\documentclass[final]{siamart220329}

\usepackage{lipsum}
\usepackage{amsfonts}
\usepackage{graphicx,subfigure}
\usepackage{epstopdf}
\usepackage{algorithmic}
\usepackage{booktabs}
\usepackage{multirow} 
\usepackage{hyperref}
\ifpdf
  \DeclareGraphicsExtensions{.eps,.pdf,.png,.jpg}
\else
  \DeclareGraphicsExtensions{.eps}
\fi

\newsiamremark{remark}{Remark}
\newsiamremark{hypothesis}{Hypothesis}
\crefname{hypothesis}{Hypothesis}{Hypotheses}
\newsiamthm{claim}{Claim}

\headers{A Neural Multigrid Helmholtz Solver}{C. Cui, K. Jiang, and S. Shu}

\title{A Neural Multigrid Solver for Helmholtz equations with high wavenumber and heterogeneous media
\thanks{Submitted to the editors DATE.}}

\author{Chen Cui\footnotemark[2]
        \and Kai Jiang\footnotemark[2]
        \and Shi Shu\thanks{Hunan Key Laboratory for Computation and Simulation in Science and Engineering,
        Key Laboratory of Intelligent Computing and Information Processing of Ministry
        of Education, School of Mathematics and Computational Science, Xiangtan University, Xiangtan, Hunan,
        China, 411105.
        (Corresponding author: Shi Shu, \email{shushi@xtu.edu.cn}).}}

\usepackage{amsopn}

\makeatletter
\newcommand*{\addFileDependency}[1]{
  \typeout{(#1)}
  \@addtofilelist{#1}
  \IfFileExists{#1}{}{\typeout{No file #1.}}
}
\makeatother

\ifpdf
\hypersetup{
  pdftitle={MG4Helm},
  pdfauthor={C. Cui, K. Jiang, and S. Shu}
}
\fi

\begin{document}

\maketitle

\begin{abstract}
     In this paper, we propose a deep learning-enhanced multigrid solver for high-frequency and heterogeneous Helmholtz equations. By applying spectral analysis, we categorize the iteration error into characteristic and non-characteristic components. 
     We eliminate the non-characteristic components by a  multigrid wave cycle, which employs carefully selected smoothers on each grid. 
     We diminish the characteristic components by a learned phase function and the approximate solution of an advection-diffusion-reaction (ADR) equation, which is solved using another multigrid V-cycle on a coarser scale, referred to as the ADR cycle.
     The resulting solver, termed Wave-ADR-NS, enables the handling of error components with varying frequencies and overcomes constraints on the number of grid points per wavelength on coarse grids. Furthermore, we provide an efficient implementation using differentiable programming, making Wave-ADR-NS an end-to-end Helmholtz solver that incorporates parameters learned through a semi-supervised training.
     Wave-ADR-NS demonstrates robust generalization capabilities for both in-distribution and out-of-distribution velocity fields of varying difficulty.
     Comparative experiments with other multigrid methods validate  its superior performance in solving heterogeneous 2D Helmholtz equations with wavenumbers exceeding 2000.
\end{abstract}



\section{Introduction}\label{sec:01}
Helmholtz equation 
\begin{equation}
 -\Delta u - k^2u = g,
 \label{eq:hel} 
\end{equation}
with wavenumber $k$ describes wave propagation in the frequency domain and has broad applications  including  acoustics, optics, electrostatics, medical and seismic imaging.
However, numerically solving it presents substantial challenges due to three key factors.
Firstly, discretizing \cref{eq:hel} requires representing each wavelength with sufficient grid points, thus resulting in a \textit{large-scale} linear system, especially noticeable with high wavenumbers. 
Secondly, the \textit{indefiniteness} of discrete system renders many numerical methods ineffective or even leads to divergence.
Thirdly, for heterogeneous Helmholtz equations, these methods available to constant wavenumber cases may encounter scalability issues. 
Many efforts have been made to address above challenges, with most of these approaches focusing on designing preconditioners for Krylov subspace methods, such as domain decomposition methods\,\cite{engquist2011sweeping,gander2019class,graham2020domain}. However, developing efficient and robust stand-alone solvers for Helmholtz equations also remains a challenging task. In this work, we aim to address these challenges by developing multigrid (MG) methods leveraging deep learning techniques.

MG methods are the state-of-the-art solvers for elliptic partial differential equations (PDEs) with linear computational complexity\,\cite{trottenberg2000multigrid,briggs2000multigrid}.
However, employing MG as a stand-alone solver for Helmholtz equation often encounters convergence issues.
The reason is attributed to the fact that standard MG methods struggle to attenuate Fourier error components $e^{i\mathbf{k}\cdot\mathbf{x}}$ with $|\mathbf{k}|\approx k$. 
These components, termed \textit{characteristic}\,\cite{brandt1997wave}, correspond to a much smaller residual since they nearly satisfy the homogeneous Helmholtz equation, rendering them nearly invisible to local operations and resulting in a slow error reduction during the  smoothing. 
On coarse grids, they are poorly approximated because their phase is altered due to insufficient grid points, rendering coarse grid correction ineffective.
To eliminate the low-frequency characteristic components, some studies restrict the size of the coarse grids and halt coarsening before these components become relatively high-frequency.
For example, an optimized discretization format has been developed to reduce the requirement for grid points per wavelength on coarse grids\,\cite{Stolk2014}.
However, this strategy compromises the scalability of MG method somewhat.
Therefore, to ensure the effectiveness of MG for solving the Helmholtz equation, it is essential to address characteristic components separately.

Currently, one of the most popular MG method for solving Helmholtz equation is the complex shifted Laplacian (CSL) preconditioner\,\cite{Erlangga2004}, which adds a complex shift to the original Helmholtz equation
$
  -\Delta u-(k^2+\beta i) u =g.
$
A large shift parameter $\beta$ displaces the eigenvalues of the CSL preconditioner distant from zero, thus restoring the effectiveness of MG methods. Obviously, a more accurate approximation needs a smaller $\beta$.
Recent researches\,\cite{gander2015applying, cocquet2017large} have shown that when $\beta$ is at most $O(k)$,  GMRES has a wavenumber-independent convergence in solving the preconditioned system provided the CSL is inverted exactly.
However, $\beta$ needs to be at least $O(k^2)$ for MG to efficiently invert the CSL preconditioner. 
There is a gap such that the CSL inverted by the standard MG V-cycle may not serve as an effective preconditioner.
Several approaches have been developed to accelerate convergence, including the utilization of different $\beta$ values on different grid levels\,\cite{Calandra2012, cools2014new} and integrating additional techniques like deflation\,\cite{Dwarka2020}. However, even with these efforts, the number of required iterations still increases  notably as the wavenumber increases.

Other approaches directly uses the information of the characteristic to design MG method, including the Wave-Ray method\,\cite{brandt1997wave, livshits2006accuracy} and characteristic-dependent near-nullspace construction method\,\cite{olson2010smoothed}.
In particular, the Wave-Ray method incorporates additional ray cycles to solve ray equations corresponding to the characteristic components within a standard MG method. 
This results in a stand-alone solver whose convergence is independent of the wavenumber for constant Helmholtz equations.
However, these methods encounter limitations when dealing with heterogeneous Helmholtz equations, as their effectiveness diminishes with pronounced variations in either amplitudes or oscillations.

Recently, there has been a growing interest in integrating deep learning techniques  with traditional iterative methods to accelerate convergence\,\cite{katrutsa2020black, cui2022fourier, kaneda2022deep}.
In the context of MG methods, current research is actively learning  smoothers\,\cite{huang2022learning, chen2022meta}, transfer operators\,\cite{greenfeld2019learning, luz2020learning}, strong threshold\,\cite{caldana2023deep, zou2023autoamg}, coarsened grids\,\cite{taghibakhshi2021optimization}, etc. 
However, most of the existing studies focus on symmetric positive definite systems, with limited attention to indefinite problems.
For the Helmholtz equation, there are a few deep learning augmented neural solvers\,\cite{azulay2022multigrid, lerer2023multigrid, drzisga2023semi, stanziola2021helmholtz}. A class of CSL-based neural Helmholtz solvers involves using UNet\,\cite{ronneberger2015u} to approximate the inverse of CSL preconditioner\,\cite{azulay2022multigrid, lerer2023multigrid}, and learning the shift parameter to speed up convergence\,\cite{drzisga2023semi}.

Motivated by aforementioned challenges and existing methods, we develop a deep learning-based MG solver, the Wave-ADR neural solver (Wave-ADR-NS), to tackle high-frequency and heterogeneous Helmholtz equations.
Our main contributions can be summarized as follows:
\begin{itemize}
     \item We design the Wave-ADR-NS by partitioning the iterative error into characteristic and non-characteristic components. 
     Non-characteristic errors are eliminated using a standard MG V-cycle, the so-called \textit{wave cycle}, employing selected smoothers on each level. 
     Characteristic errors are handled by solving an advection-diffusion-reaction (ADR) equation using another MG V-cycle on a coarse scale, termed as \textit{ADR cycle}.
     This approach addresses error components with different frequencies and overcomes limitations on the number of grid points per wavelength on coarse grids.
     \item We provide an efficient implementation using differentiable programming. This enables Wave-ADR-NS to function as an end-to-end Helmholtz solver with support for backpropagation, matrix-free computation, batch processing, and GPU acceleration.
     During the offline setup phase, parameters in Wave-ADR-NS are optimized in a semi-supervised manner.
     \item We test the effectiveness and robustness of Wave-ADR-NS on heterogeneous two-dimensional (2D) Helmholtz equations with wavenumbers up to 2000 for both in-distribution and out-of-distribution problems. Numerical experiments demonstrate that Wave-ADR-NS maintains a stable convergence rate across different velocity fields, requiring fewer iterations and less computational time compared to traditional MG preconditioners, such as the Wave-Ray method and CSL inverted by an MG V-cycle, as well as deep learning-enhanced MG preconditioners like Encoder-Solver \cite{azulay2022multigrid} and Implicit Encoder-Solver \cite{lerer2023multigrid}.
\end{itemize}

The remainder of the paper is organized as follows. \Cref{sec:02} introduces the Helmholtz equation and its discretization method. \Cref{sec:03} presents the Wave-ADR-NS and discusses its implementation. \Cref{sec:04} showcases numerical results. Finally, \Cref{sec:05} offers summary and outlook.

\section{Problem Formulation}\label{sec:02}

Consider the 2D Helmholtz equation that satisfies the Sommerfeld radiation condition far from the region of interest
\begin{equation} 
     \begin{gathered}
     -\Delta u(\mathbf{x})-\left(\frac{\omega}{c(\mathbf{x})}\right)^2 u(\mathbf{x})=g(\mathbf{ x}), \\
     \text { s.t. } \lim _{|\mathbf{x}| \rightarrow \infty}|\mathbf{x}|^{\frac{1}{2}}\left(\frac{\partial} {\partial|\mathbf{x}|}-i k\right) u(\mathbf{x})=0,  
     \end{gathered} 
     \label{eq:helm}
\end{equation}
where $k(\mathbf{x})=\omega/c(\mathbf{x})$ is the wavenumber, $\omega=2\pi f$ is the angular frequency, $f$ is the frequency, $c(\mathbf{x})$ is the wave speed, $g(\mathbf{x}) = \delta(\mathbf{x}-\mathbf{x}_0)$ is the Dirac delta function. 
When investigating solution behavior within specific regions, truncating the computational domain with appropriate boundary conditions is necessary. 
Due to the slow decay of the Helmholtz solution, Dirichlet or Neumann boundary conditions tend to reflect waves, requiring the absorbing boundaries instead.
Common choices include impedance boundary conditions\,\cite{engquist1979radiation}, sponge layers\,\cite{israeli1981approximation}, or perfectly matched layers\,\cite{berenger1994perfectly}. 
One approach is adding a damping term to the original equation
\begin{equation} 
     -\Delta u-\frac{\omega^2}{c^2}u+i\omega\gamma\frac{1}{c^{2}}u=g,\quad \mathbf{x}\in \Omega,  
     \label{eq:sponge}
\end{equation}
where $\gamma$ is the damping mask, which equals to zero
inside $\Omega$ and increasing from zero to $\omega$ within the sponge layer near the boundary. 
The thickness of this layer is typically one wavelength\,\cite{erlangga2006novel,treister2019multigrid}.

Using the second-order central finite difference method (FDM) on a uniform grid yields the discretized Helmholtz operator
\begin{equation}
    \frac{1}{h^2}\left[\begin{array}{ccc} 
     0 & -1 & 0 \\
     -1 & 4-k_{i,j}^2 h^2 & -1 \\
     0 & -1 & 0  
     \end{array}\right],
\end{equation}
where $h=1/(N+1)$ in both the $x-$ and $y-$ directions, $N$ is the number of interior nodes. 
According to the Shannon sampling principle, there should be at least $10$ grid points per wavelength, i.e., $\lambda=2\pi/k=1/f\geq 10h$, requiring $N\geq 10f$, which leads to a \textit{large-scale} linear system
\begin{equation}
    \mathbf{A}\mathbf{u}=\mathbf{g}. 
    \label{eq:linsys} 
\end{equation}

Even worse, as the wave number increases, the \textit{indefinite} property of $\mathbf{A}$ becomes more pronounced. As shown in \eqref{eq:eigs}, even for the 1D model problem, negative eigenvalues proliferate with increasing wavenumber, making many iterative methods ineffective or even divergence.

\section{Our proposed method}\label{sec:03}
In this section, we first design a MG V-cycle, termed the wave cycle, to eliminate non-characteristic error components. Although much of this content is well-known, our novel contribution lies in introducing the Chebyshev semi-iteration method with learnable parameters as a smoother for the first time. To clarify the rationale and innovation behind this method, we review the necessary analysis.
Next, we introduce a correction step, called the ADR cycle, to specifically target characteristic error components using the learned representation. Finally, we present a differentiable implementation of the proposed method, enabling efficient training.

\subsection{Wave cycle for non-characteristic error}\label{sec:31}
For clarity and without loss of generality, we design the wave cycle on the following 1D model problem
\begin{equation} 
    -u(x)^{\prime\prime}-k^2u(x)=g(x),\quad u(0)=u(1)=0,  
\end{equation}
where wave speed $c=1$. 
The coefficient matrix obtained through FDM is 
\begin{equation}\nonumber
    \mathbf{A}=(1/h^2)\text{tridiag}(-1,2,-1)-k^2\mathbf{I}.
\end{equation}
One advantage of using this model problem is that its eigenvalues
\begin{equation} 
    \begin{aligned}
    \lambda_j&=\frac{2(1-\cos j\pi h)}{h^2}-k^2=\frac{4}{h^2}\sin^2\frac{j\pi h}{2}-k^2,\quad j=1,\ldots,N,
    \end{aligned}  
    \label{eq:eigs}
\end{equation}
and eigenvectors 
\begin{equation}  
    \mathbf{v}_j=\left\{\sin ij\pi h\right\}_{i=1}^N,\quad j=1,\ldots,N, 
    \label{eq:eigvec}
\end{equation}
are known, enabling us to perform error analysis easily from the spectral viewpoint. 
It can also be observed from \eqref{eq:eigvec} that when $j$ is small (e.g., $j < N/2$), the corresponding eigenvector represents a \textit{low-frequency} mode, whereas when $j$ is large ($j \geq N/2$), the eigenvector represents a \textit{high-frequency} mode. Since the set $\{\mathbf{v}_j\}_{j=1}^N$ constitutes an orthogonal basis of $\mathbb{R}^N$, the iterative error can be expanded under this basis 
\begin{equation}
    \mathbf{e}=\sum_{j=1}^N c_j\mathbf{v}_j=\sum_{j<N/2} c_j\mathbf{v}_j+\sum_{j\geq N/2}c_j\mathbf{v}_j. 
    \label{eq:err}
\end{equation}
MG methods consist of two complementary operations: smoothing and coarse grid correction, which eliminate high-frequency and low-frequency error components respectively. Next, we analyze the errors of these two operations separately.

\subsubsection{Error of the smoother}
The general iteration scheme of a smoother is given by
\begin{equation}
     \mathbf{u}^{(q+1)}=\mathbf{u}^{(q)}+\mathbf{B}(\mathbf{g}-\mathbf{A}\mathbf{u}^{(q)}),  
\end{equation}
where $\mathbf{B}$ is an approximation of $\mathbf{A}^{-1}$ and is computationally efficient. 

We employ the damped Jacobi method as an example to conduct error analysis. Similar analysis can be applied to other smoothers as well.
The error propagation matrix of the damped Jacobi method is
\begin{equation} 
\mathbf{S}_\omega=\mathbf{I}-\omega \mathbf{D}^{-1}\mathbf{A},
\label{eq:sjac}
\end{equation} 
where $\mathbf{D}$ is the diagonal of $\mathbf{A}$.
From the eigenvalues of $\mathbf{A}$ in \cref{eq:eigs}, a direct calculation obtains the eigenvalues of $\mathbf{S}_\omega$ 
\begin{equation}
\mu_j=1-\omega\left(1-\frac{2\cos j\pi h}{2-k^2h^2}\right),\quad j=1,\ldots, N.  
\end{equation}
Assume that the initial error can be  expanded as in \cref{eq:err}. After $q$ relaxation steps, the error becomes
\begin{equation} 
\mathbf{e}^{(q)}=\mathbf{S}_\omega^q \sum_{j=1}^N c_j\mathbf{v}_j=\sum_{j=1}^N c_j\mu_j^q \mathbf{v}_j.  
\end{equation}  
Obviously, when $\mu_j$ is close to zero, the corresponding error component $\mathbf{v}_j$ decays quickly. However, when $\mu_j$ is close to one or even larger than one, $\mathbf{v}_j$ decays slowly or even amplifies.

\begin{figure}[!htb]
     \centering
     \subfigure[$k=8\pi,N=47$]{\label{subfig:eig1}\includegraphics[width=0.29\textwidth]{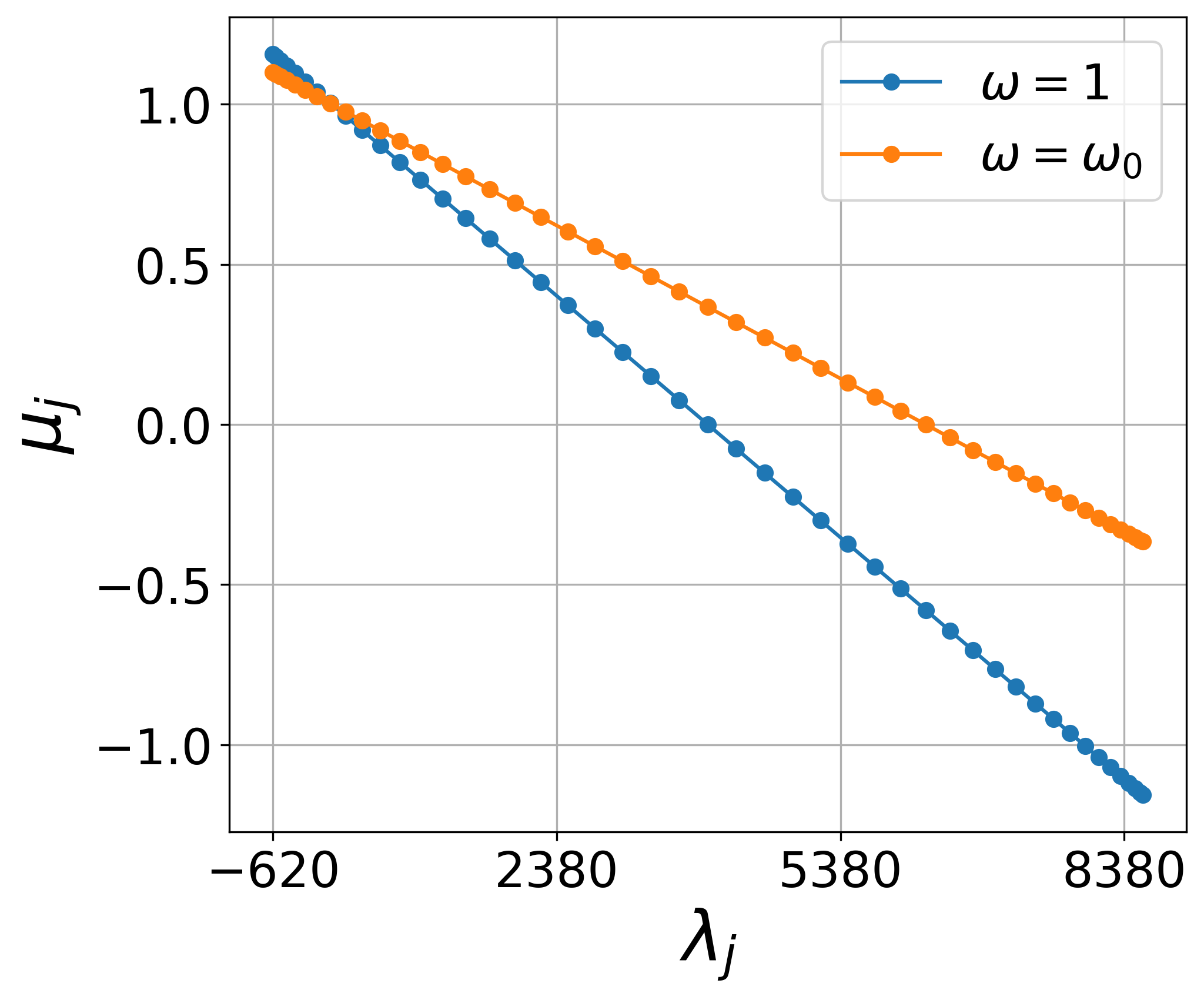}}
     \subfigure[$k=8\pi,N=23$]{\includegraphics[width=0.28\textwidth]{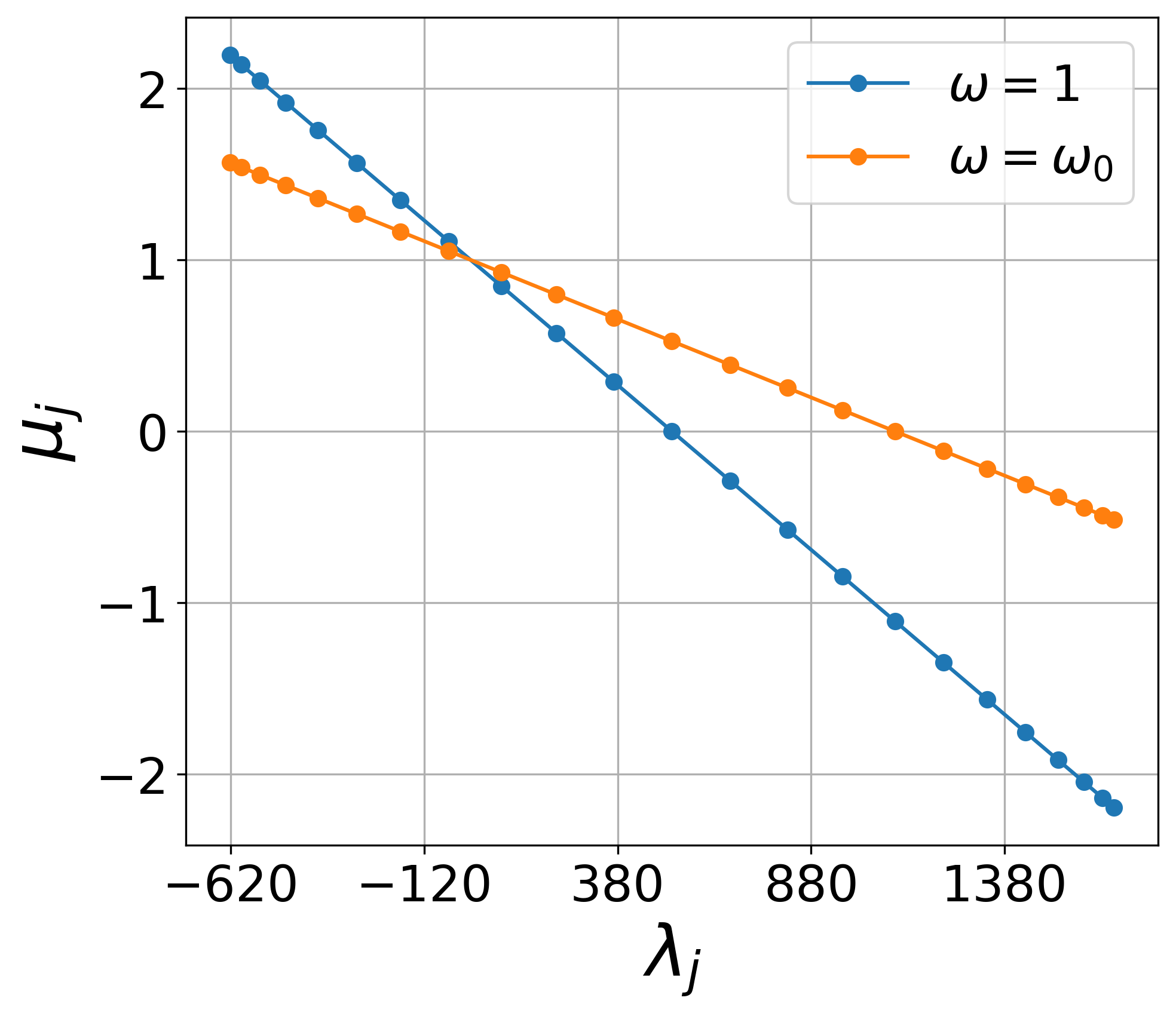}}
     \subfigure[$k=8\pi,N=11$]{\includegraphics[width=0.3\textwidth]{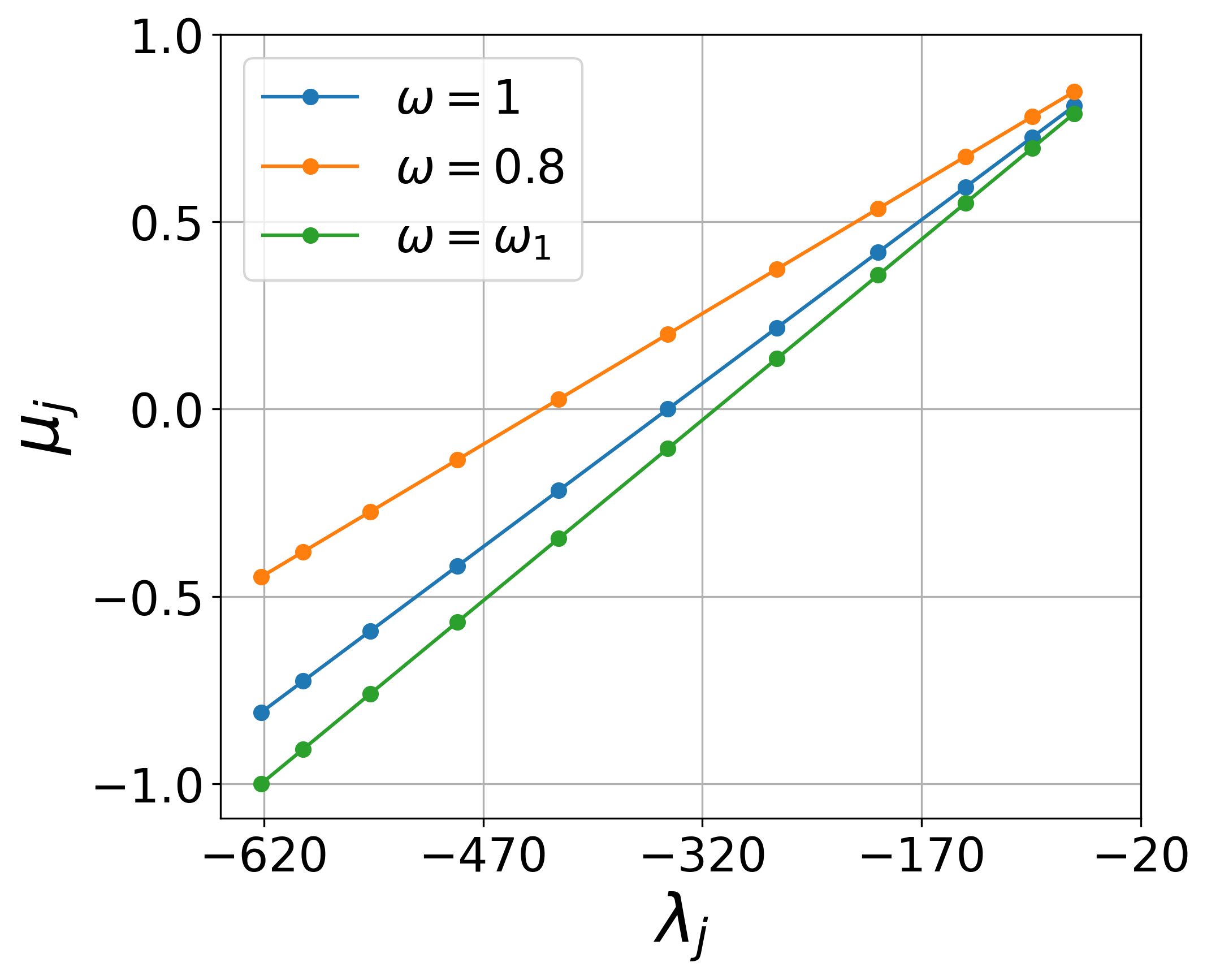}}
     \caption{Eigenvalues of the error propagation matrix in  \eqref{eq:sjac} for the damped Jacobi method on nested mesh grids.}
     \label{fig:jacobi_eigs}
\end{figure}
\cref{subfig:eig1} displays the eigenvalues of $\mathbf{S}_\omega$  when $f=4,N=47$. On the finest grid, one can adjust $\omega$ to dampen the most oscillatory mode.
The optimal choice is $\omega_0=\frac{2-k^2h^2}{3-k^2h^2}$ \cite{elman2001multigrid,ernst2012difficult}. 
As the grid coarsens, the amplification factor for smooth errors increases, making the damped Jacobi method being not applicable on this level.
Further coarsening grids,  the coarse operator becoming negative definite, and the damped Jacobi method converges again when $\omega \in (0, \omega_1)$, where $\omega_1=\frac{2-k^2h^2}{2\sin^2(\pi h/2)-k^2h^2/2}$. 
However, on this coarse level, the damped Jacobi method fails to play a smoothing role as it cannot eliminate high-frequency error components.

Currently, there are mainly two types of methods that provide a smoothing effect on coarse grids for the Helmholtz equation: Krylov subspace methods like GMRES\,\cite{elman2001multigrid}, and iterative methods appliedon the normal equation\,\cite{brandt1997wave}. 
To avoid the difficulty of determining optimal GMRES steps\,\cite{elman2001multigrid}, we use the iterative method on the normal equation. Concretely, we use the Chebyshev semi-iteration method\,\cite{adams2003parallel,richefort2023toward}
\begin{equation}  
     \mathbf{u}^{(q+1)}=\mathbf{u}^{(q)}+\beta_q\left(\mathbf{A^*g}-\mathbf{A^*\mathbf{A} u}^{(q)}\right),
     \label{eq:cheb}
\end{equation}  
where $\mathbf{A^*}$ denotes the conjugate transpose of $\mathbf{A}$. The corresponding error propagation matrix is $\mathbf{S}_{\beta_q}=\mathbf{I}-\beta_q \mathbf{A^*\mathbf{A}}.$
After applying the iteration $q$ times, the error becomes
\begin{equation}
     \mathbf{e}^{(q)}=\left(\mathbf{I}-\beta_{q-1} \mathbf{A^*\mathbf{A}} \right) \mathbf{e}^{(q-1)}
                    =\ldots :=p(\mathbf{A^*}{\mathbf{A}}) \mathbf{e}^{(0)},
\end{equation}  
where $p(\mathbf{A^*\mathbf{A}})$ is a matrix polynomial of order $q$ with $p(0)=\mathbf{I}$.
$\beta_i$ are chosen to minimize the maximal modulus of $p$ within the interval $\left[\lambda_{\max}/\alpha, \lambda_{\max}\right]$, $\alpha$ is determined such that $\lambda_{\max}/\alpha$ approximates the median of all eigenvalues of $\boldsymbol{A^*\boldsymbol{A}}$ to compress high-frequency error components.

\begin{figure}[!htbp]
     \centering  
     \subfigure[$k=8\pi,N=47, \alpha=4.6$]{\includegraphics[width=0.43\textwidth]{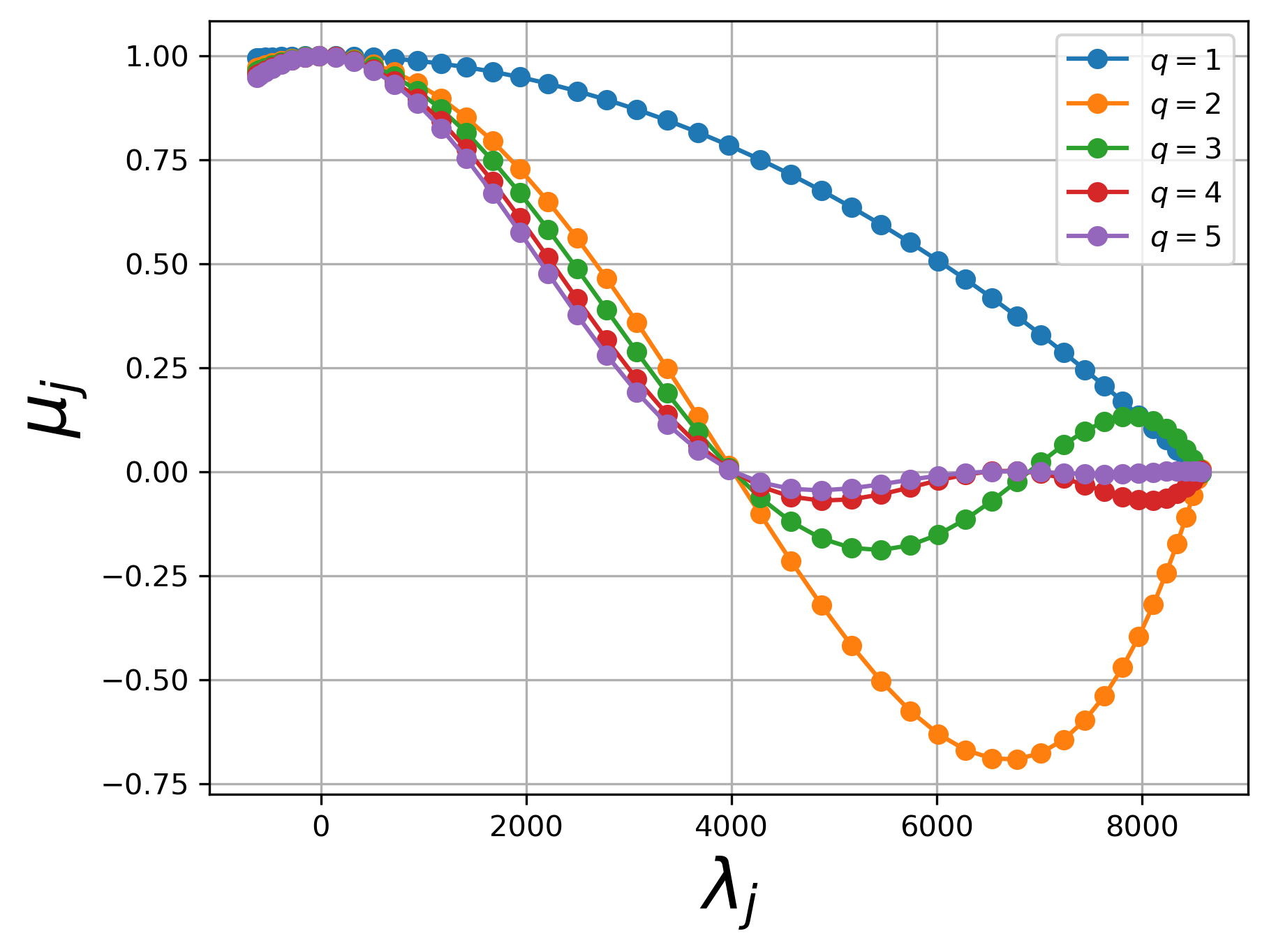}}\quad
     \subfigure[$k=8\pi,N=23, \alpha=7.1$]{\includegraphics[width=0.43\textwidth]{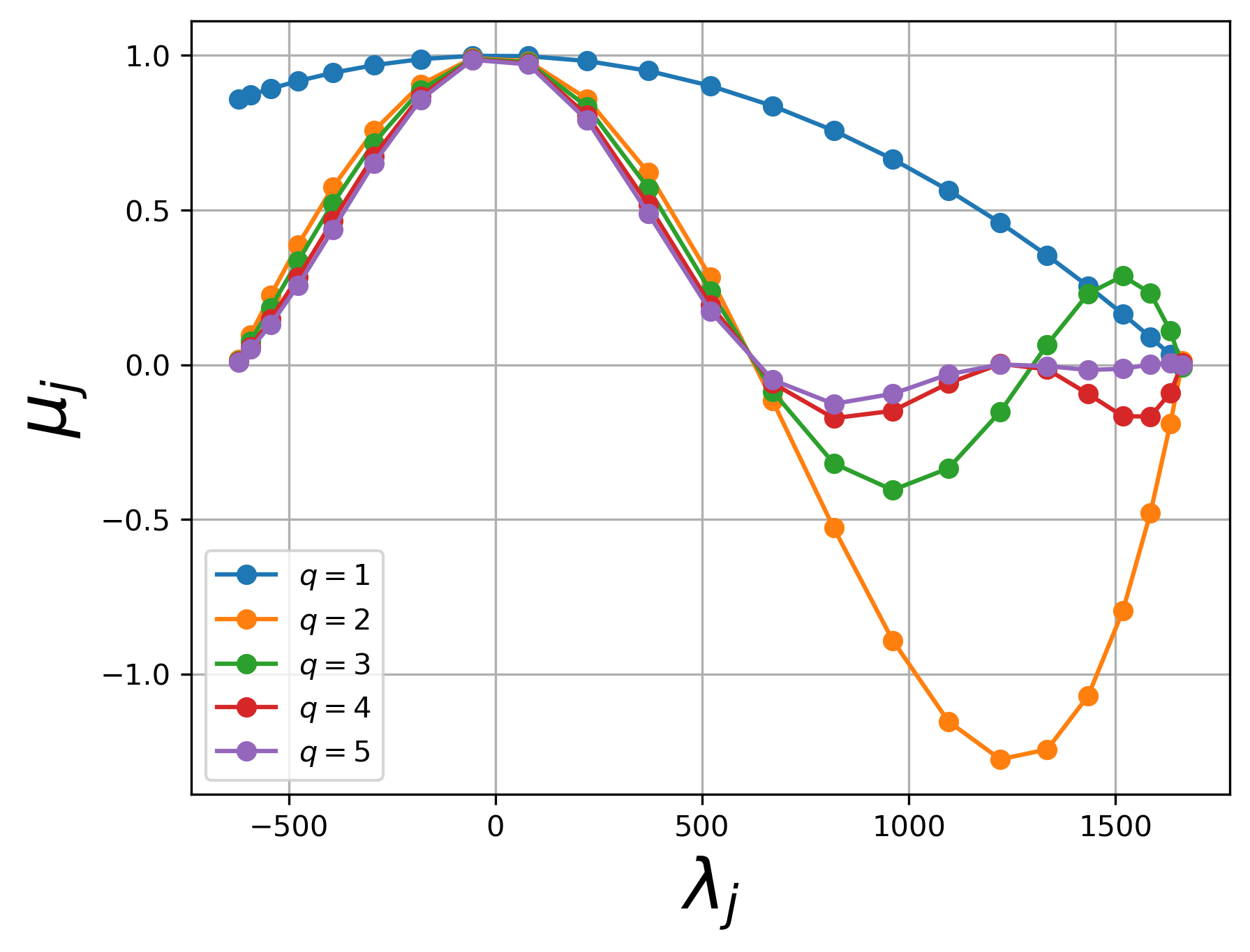}}\\
     \subfigure[$k=8\pi,N=11, \alpha=3.3$]{\includegraphics[width=0.43\textwidth]{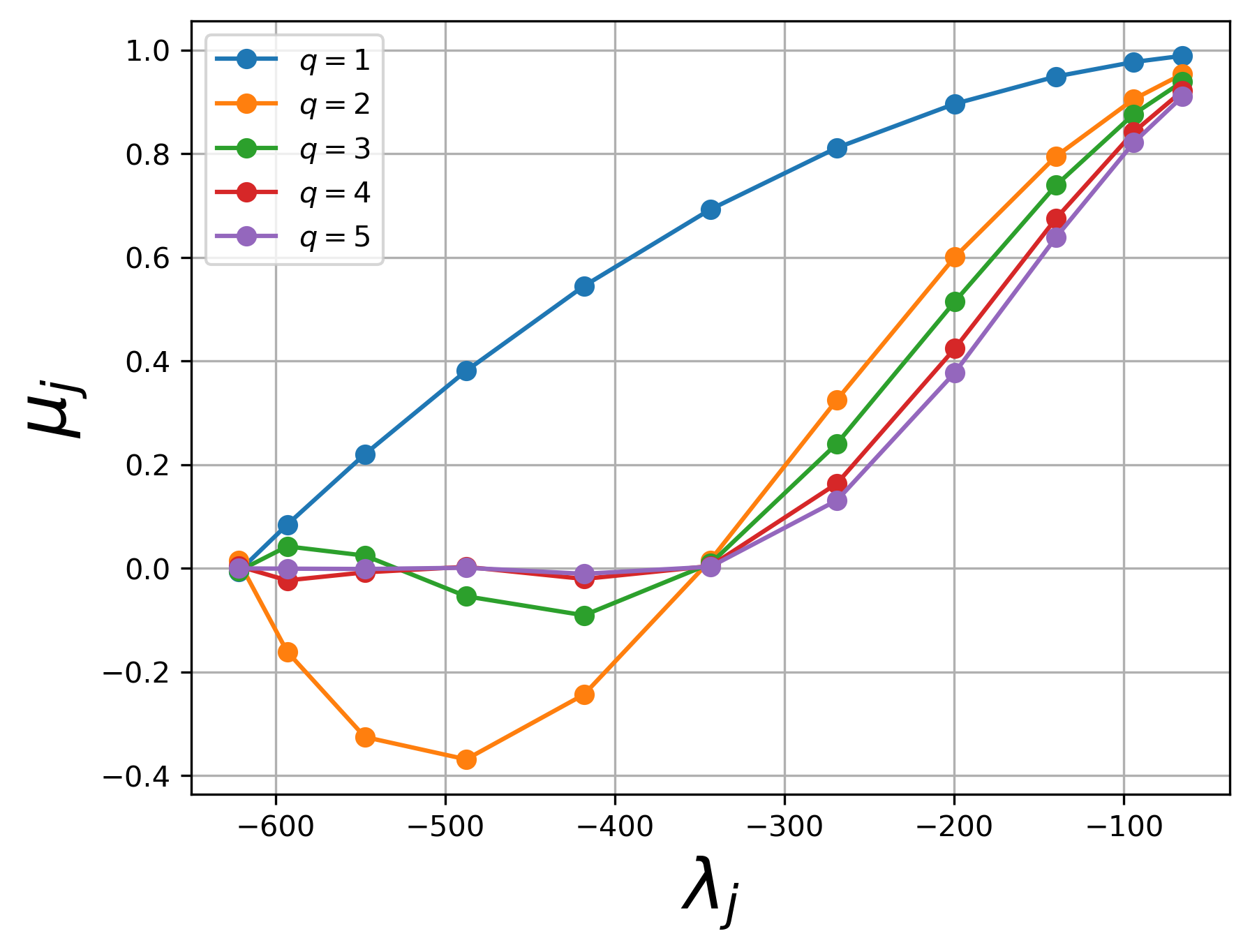}}  \quad
     \subfigure[$k=8\pi,N=5, \alpha=1.2$]{\includegraphics[width=0.43\textwidth]{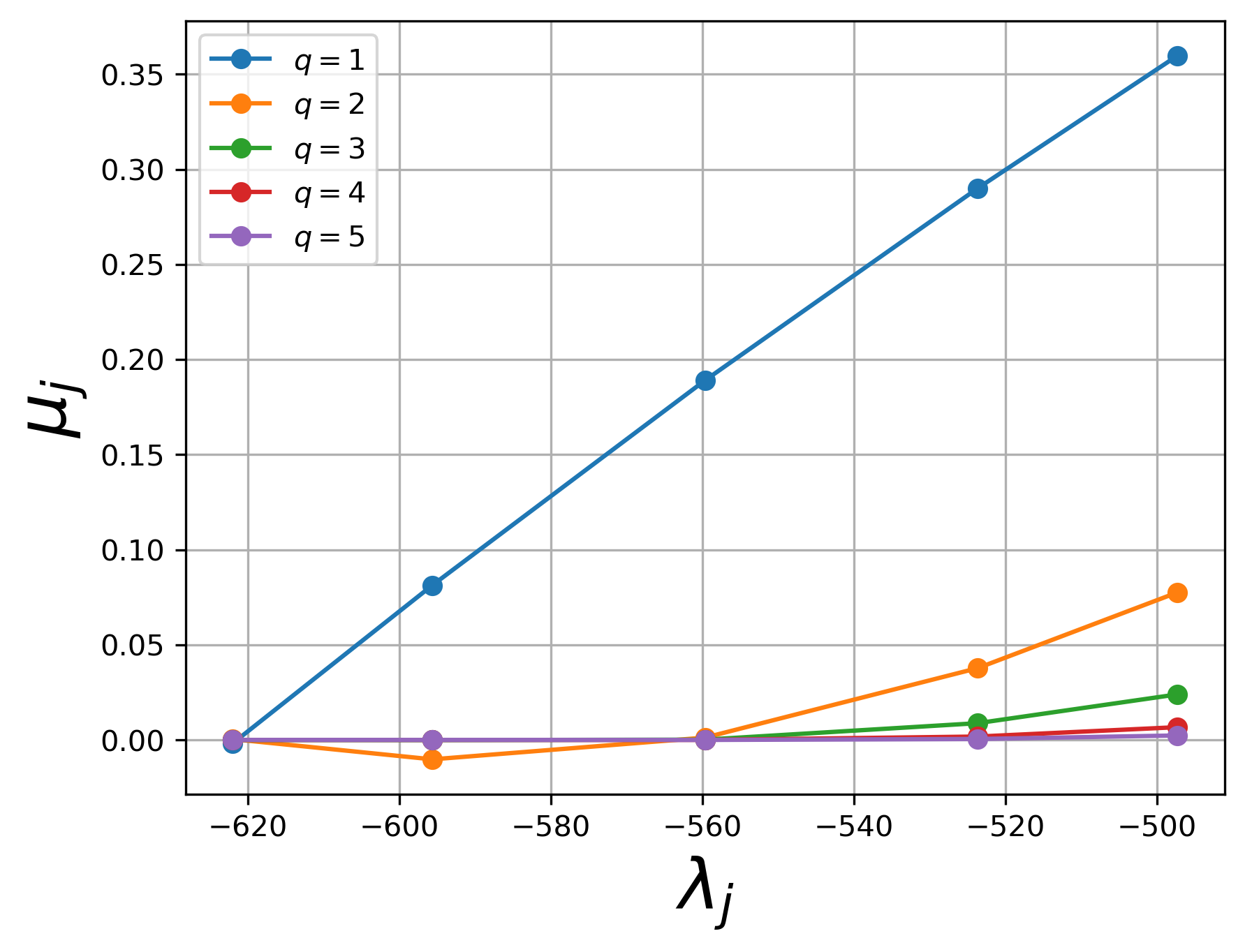}}
     \caption{Eigenvalues of the error propagation matrix of the Chebyshev semi-iteration method on nested mesh grids.}  
     \label{fig:cheby_eigs}
\end{figure}
\cref{fig:cheby_eigs} shows the eigenvalues of $\mathbf{S}_{\beta_q}$. As $q$ increases, the Chebyshev semi-iteration method exhibits an improved smoothing effect, except on the third level for the following reason. 
On the third level ($kh \approx 2$), the largest eigenvalue of $\mathbf{A}$ close to zero, which corresponds to the smallest eigenvalues of $\mathbf{A^*A}$. 
Since the Chebyshev semi-iterative method targets errors associated with large eigenvalues, the error corresponding to near-zero eigenvalues is the least likely to be eliminated.
Moreover, $\alpha$ varies across levels, making its calculation a significant challenge. To address this, we employ a neural network to learn these parameters.

\subsubsection{Error of the coarse grid correction} 
Let $\mathbf{A}^h$ denote the coefficient matrix on the fine grid, and $\mathbf{A}^H$ represent the coarse grid counterpart. The error propagation matrix of the coarse grid correction operator in the two-grid method is given by
\begin{equation}\nonumber
\mathbf{E}=\mathbf{I}-\mathbf{I}_{H}^{h}(\mathbf{A}^{H})^{-1} \mathbf{I}_{h}^{H}\mathbf{A}^h,
\end{equation}
where $\mathbf{I}_{h}^{H}$ is the full weighting operator
\begin{equation} \nonumber
\left[\mathbf{I}_h^H\mathbf{u}^h\right]_i:=\frac{1}{4}\left(\left[\mathbf{u}^h\right]_{2i -1}+2\left[\mathbf{u}^h\right]_{2i}+\left[\mathbf{u}^h\right]_{2i+1}\right),\quad i= 1,\ldots,(N-1)/2,  
\end{equation}  
and the interpolation operator $\mathbf{I}_H^h=2(\mathbf{I}_h^H)^{\top}$.

From the classical two-grid analysis \cite{elman2001multigrid}, we know that the compressibility of $\mathbf{E}$ for low-frequency error components $\mathbf{v}_j$ ( $j < N/2$) is
\begin{equation} 
     \mathbf{Ev}_j \approx \mathbf{v}_j - \frac{\lambda^h_j}{\lambda^H_j}\mathbf{I}_H^h\mathbf{I}_h^H\mathbf{v}_j = \left(1 - \frac{\lambda^h_j}{\lambda^H_j}\right)\mathbf{v}_j.
\end{equation}
where $\lambda^h_j$ and $\lambda^H_j$ are the eigenvalues of $\mathbf{A}^h$ and $\mathbf{A}^H$ coresponding to $\mathbf{v}_j$, respectively.
Thus, the compressibility depends on the ratio $\lambda^h_j/\lambda^H_j$. When $\lambda^h_j$ and $\lambda^H_j$ are equal, coarse grid correction is perfect. However, the trouble arises when one eigenvalue is close to zero while the other is not, or if they have opposite signs. 
From the eigenvalues expression \eqref{eq:eigs}, it can be seen that when $j\approx 2f$, the eigenvalues are close to zero. 
During the coarsening process, some eigenvalues may undergo a slight leftward shift, potentially causing positive eigenvalues to become negative on levels where the largest eigenvalue is close to zero.
The eigenvectors associated with these eigenvalues correspond to components that are challenging to eliminate or may even be amplified by coarse grid correction.

\subsubsection{Summary} 
From the preceding analysis, we can derive a MG V-cycle for the Helmholtz equation that employs the damped Jacobi method on the finest grid, the Chebyshev semi-iteration method with learnable parameters on the coarser grids, and no iterations on the third level ($kh \approx 2$). We refer to this V-cycle as the \textit{wave cycle}, following the Wave-Ray method \cite{brandt1997wave}.

However, even with these carefully chosen smoothers, the wave cycle still faces challenges in solving the Helmholtz equation:
\begin{enumerate}
\item Both smoothing and coarse grid correction struggle to eliminate error components with frequency $f$;
\item Error components associated with eigenvalues that change sign during coarsening (frequencies near $f$) are also difficult to eliminate.
\end{enumerate}

Next, we introduce the ADR cycle to address these problematic components.

\subsection{ADR cycle for characteristic error}
\begin{figure}[!htb]
     \centering
     \subfigure[$u$]{\label{subfig:helmholtz}\includegraphics[width=0.35\textwidth]{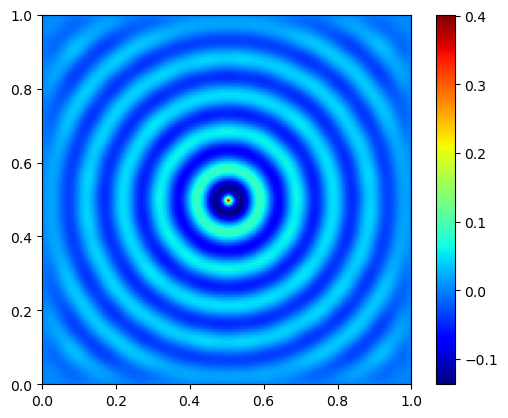}}\quad\quad
     \subfigure[$|\mathcal{F}(e)|$]{\label{subfig:smoothererror}\includegraphics[width=0.35\textwidth]{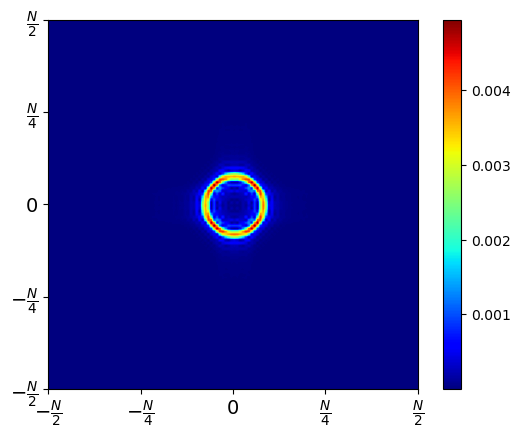}}
     \caption{(a) Solution of Helmholtz equation in homogeneous media with a point source and absorbing boundary condition.
     (b) The module of the error in Fourier space between one wave cycle iterative solution and the reference solution.}
\end{figure}

We use the wave cycle derived above to solve the 2D Helmholtz equation with an absorbing boundary condition. 
\cref{subfig:helmholtz} shows the numerical solution obtained using FDM and a direct linear solver for $f = 10$ and $c = 1$, with the source point located at $\mathbf{x}_0 = (0.5, 0.5)$. 
\cref{subfig:smoothererror} illustrates the  error in Fourier space between the wave cycle iterative solution and the reference solution from \cref{subfig:helmholtz}. As expected from the preceding analysis, the error concentrates near the circle with radius $f$.

\begin{figure}[!htb]
     \centering
    \includegraphics[width=0.35\textwidth]{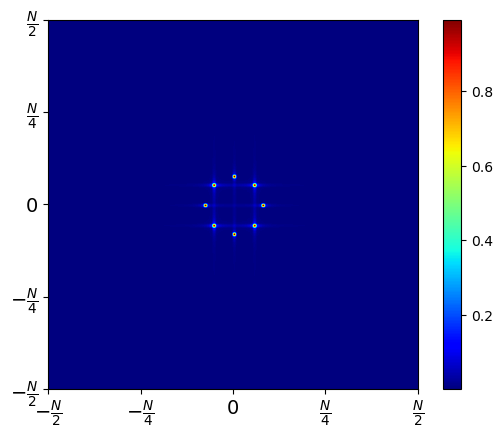}
     \caption{The module of Fourier transform of basis functions used in Wave-Ray method.}
     \label{subfig:rayhat}
\end{figure}

The Wave-Ray method\,\cite{brandt1997wave,livshits2006accuracy} is based on the fact that the characteristic error $v(\mathbf{x})$ can be expressed as 
$$
     v(\mathbf{x})=\sum_{m=1}^M v_m(\mathbf{x}) = \sum_{m=1}^M a_{m}(\mathbf{x}) e^{i\omega\left(k_1^{m} x_1+k_2^{m} x_2\right)},
$$
where
$$
     \left(k_{1}^m, k_{2}^m\right)=\left(\cos \theta_m , \sin \theta_m\right):=\vec{k}^m, \quad \theta_m = (m-1) \frac{2 \pi}{M}, \quad m=1, \dots, M,
$$
and $a_{m}$ are smooth ray functions.
Assuming the residual corresponding to the component $v_m(\mathbf{x})$ is $r_m(\mathbf{x})$, substituting these into \cref{eq:sponge} (where $v_m(\mathbf{x})$ is the solution and $r_m(\mathbf{x})$ is the right-hand side) gives the following equation for $a_{m}$
$$
     \Delta a_{m} + 2i\omega\vec{k}^m \cdot \nabla a_{m} + i\omega\gamma s^2 a_{m} = \hat{r}_m,
$$
where $\hat{r}_m = r(\mathbf{x})e^{-i\omega\left(k_1^{m} x_1+k_2^{m} x_2\right)}$. Solving these ray equations provides the correction value $v(\mathbf{x})$ for the wave cycle.
\cref{subfig:rayhat} demonstrates the distribution of the used plane wave basis functions $e^{i\omega\left(k_1^{m} x+k_2^{m} y\right)}$ for $M=8$. 
It can be observed that a sufficiently large $M$ is required to approximate the characteristic components. As the wavenumber increases, the radius of this circle also expands, necessitating a larger value of $M$.
Moreover, for the heterogeneous Helmholtz equation, the distribution of the characteristic no longer resembles a ring (see \cref{fig:7}), rendering the use of such basis functions inefficient.

To diminish the characteristic error, we employ neural networks to learn its representation. To achieve this, we begin by modeling the expression of the characteristic error. Given that it approximates the solution of the homogeneous Helmholtz equation, we aim to find the basis function capable of representing this error within the kernel space of the Helmholtz operator $ \Delta + k^2 $.
Concretly, we generalize the plane wave to variable wavenumber cases by the form  $e^{-i\omega\tau(\mathbf{x})}$, where $\tau$ is a general phase function. 
Suppose that $e^{-i\omega\tau(\mathbf{x})}$ belongs to the kernel space of the Helmholtz operator, that is
\begin{equation}
     |\nabla \tau(\mathbf{x})|^2-1/c^2(\mathbf{x})=\frac{-i}{w}\Delta\tau(\mathbf{x}).
     \label{eq:kernel}
\end{equation}
As the right-hand side of \cref{eq:kernel} is small for large wavenumbers, we let $\tau(\mathbf{x})$ satisfy the following eikonal equation
\begin{equation}
     |\nabla \tau(\mathbf{x})|^2=1/c^2(\mathbf{x}):=s^2(\mathbf{x}),\quad \tau(\mathbf{x }_0)=0,
     \label{eq:eikeq}
\end{equation}
where $s(\mathbf{x}) = 1/c(\mathbf{x})$ is the \textit{slowness} model. 
The resulting $e^{-i\omega\tau(\mathbf{x})}$ is anticipated to approximate the kernel space of the Helmholtz operator, facilitating effective representation of the characteristic.

\begin{figure}[!htbp]
     \centering
     \subfigure[$\tau=\|\mathbf{x}-\mathbf{x}_{0}\|_{2}$]{\label{subfig:eikonal}\includegraphics[width=0.27\textwidth]{ 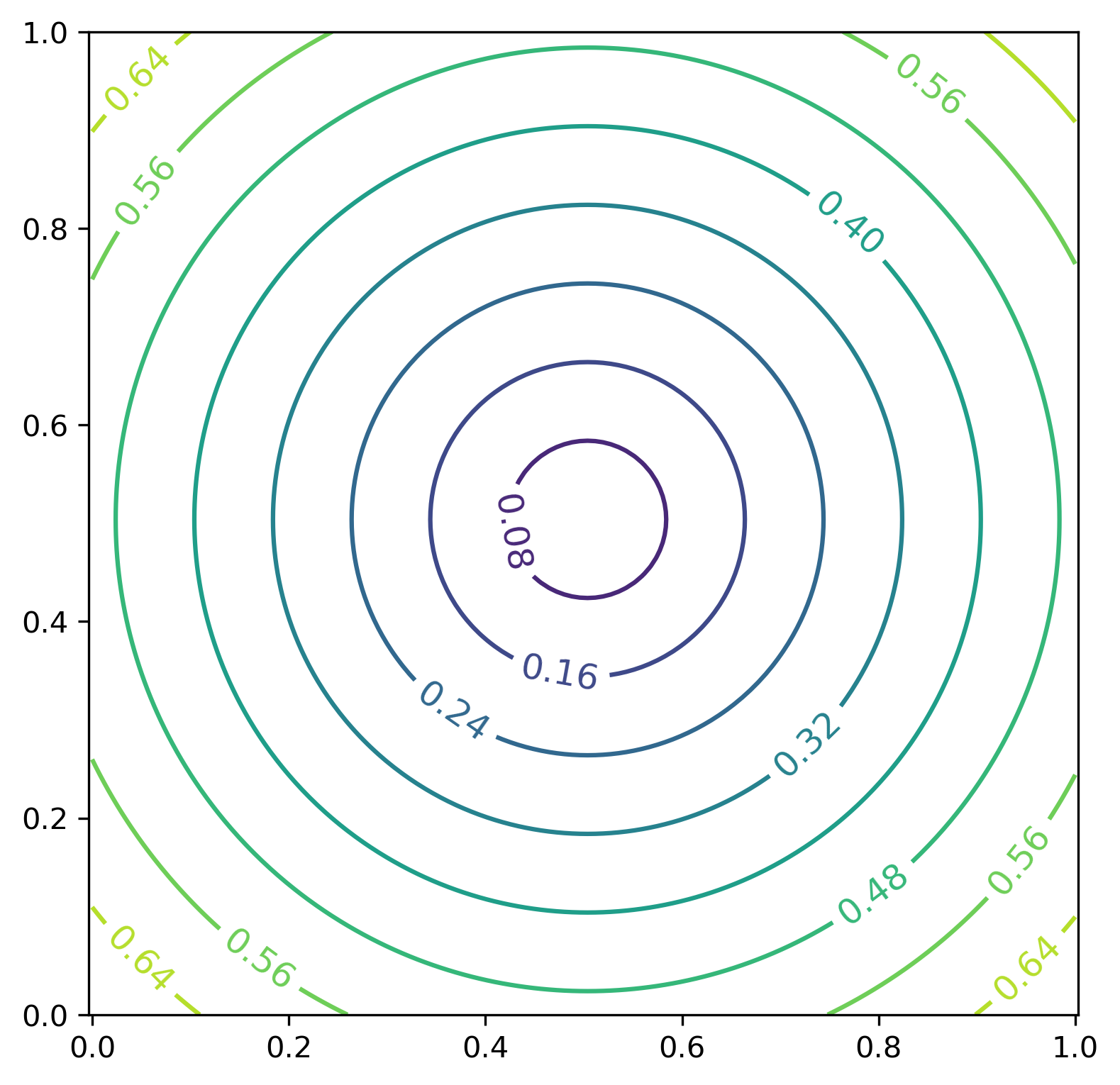}}  
     \subfigure[$\nabla \tau$]{\label{subfig:constant_tau}\includegraphics[width=0.27\textwidth]{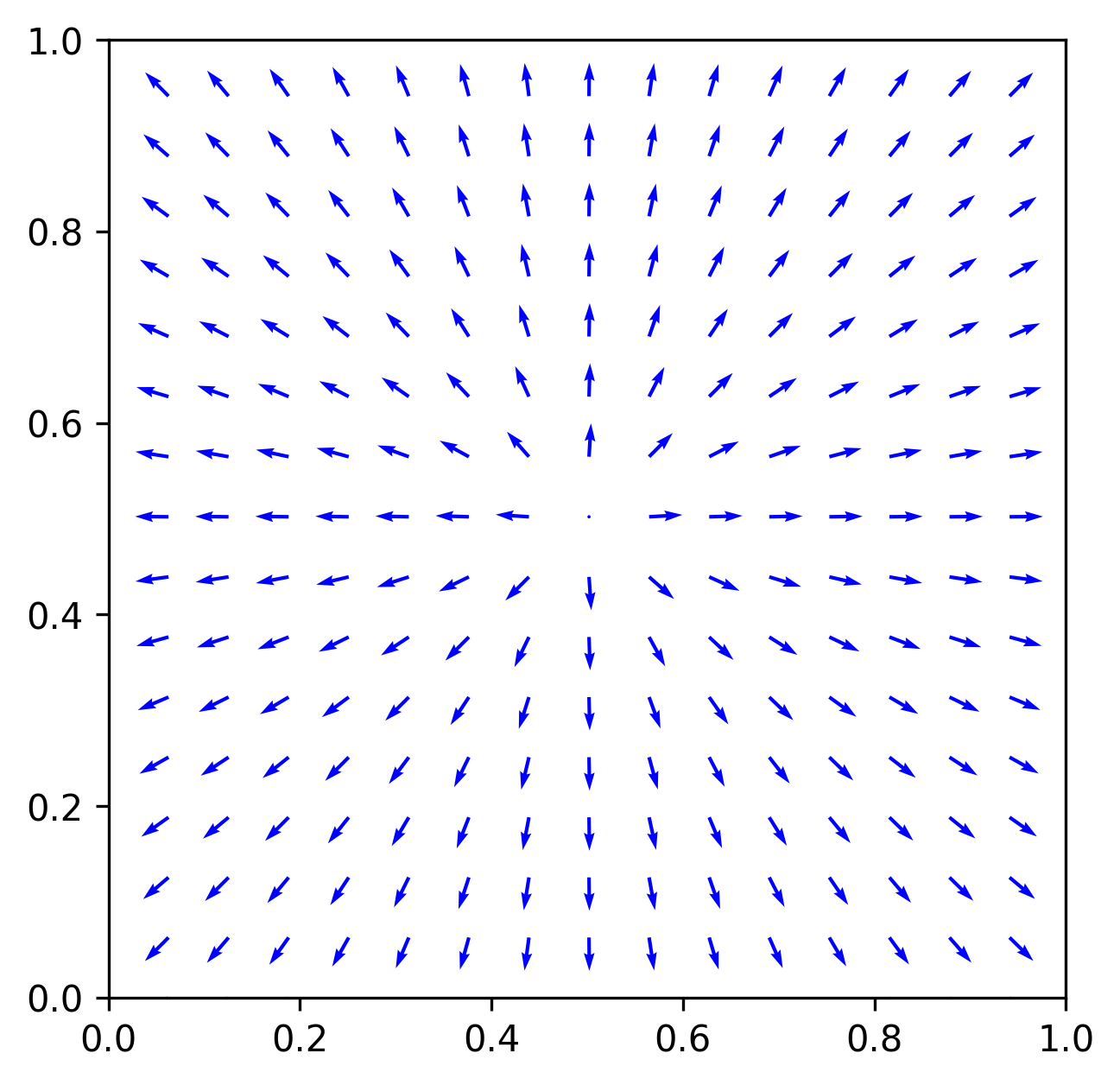}}
     \subfigure[$|\mathcal{F}(e^{-i\omega \tau})|$]{\label{subfig:pwhat}\includegraphics[width=0.31\textwidth]{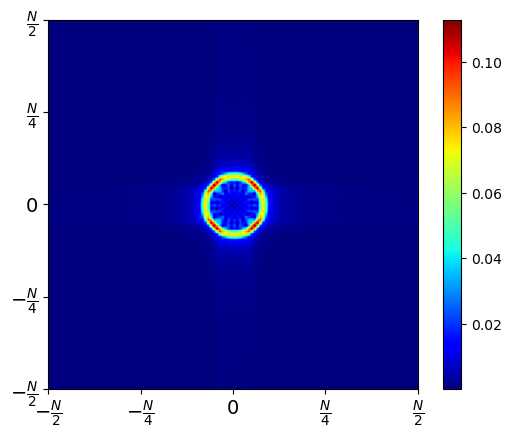}} 
     \caption{(a) Solution of the eikonal equation with uniform media. (b) Gradient field of the eikonal solution. (c) The module of Fourier transform of $e^{-i\omega \tau}$.}    
\end{figure}
\cref{subfig:eikonal} shows the solution of the eikonal \cref{eq:eikeq} for $s=1$
\begin{equation}
     \tau=\|\mathbf{x}-\mathbf{x}_{0}\|_{2}. 
     \label{eq:eikonal_cons}
\end{equation}  
\cref{subfig:constant_tau} depicts the gradient field of $\tau$, which can reflect the direction of wave propagation. 
\cref{subfig:pwhat} illustrates the Fourier transform of $e^{-i\omega\tau}$, closely matching the characteristic components.  

Therefore, we substitute the $M$ linear phase functions used in the Wave-Ray method with $\tau$, which is obtained from solving the eikonal equation.
We use the expression  
\begin{equation}  
     v(\mathbf{x})=a(\mathbf{x}) e^{-i \omega \tau(\mathbf{x})},   
     \label{eq:Rytov}  
\end{equation}  
to approximate the characteristic components, where the undetermined amplitude function $ a(\mathbf{x}) $ is expected to be independent of $ \omega $.
By substituting it into \cref{eq:sponge},
we obtain an advection-diffusion-reaction (ADR) equation on $a$
\begin{equation}  
     -\Delta a + 2i\omega\nabla\tau\cdot\nabla a + i\omega(\Delta \tau)a + \omega^2(|\nabla\tau|^2-s^2)a + i\omega\gamma s^2a =\hat{r},
     \label{eq:ADR}
\end{equation}
with $\hat{r}=re^{i \omega \tau}$, $r$ is the residual corresponding to $v$.
Solving the ADR equation and then correcting the error \cref{eq:Rytov} for the wave cycle completes one Wave-ADR cycle. 

\cref{alg:cycle} outlines the computation process for a single Wave-ADR cycle. Compared to a standard MG V-cycle, the Wave-ADR incorporates ADR correction in lines 23-25 to eliminate characteristic errors. 
The chosen second level for correction strikes a balance between effectiveness and computational efficiency: solving the ADR equation on the finest level is prohibitively expensive, applying the correction on the third level is inefficient, and characteristic errors become indistinguishable on coarser levels.
Since only a fixed number of matrix-vector multiplications are performed on each level, the computational complexity of a single Wave-ADR cycle remains $O(N)$.
Next, we discuss how neural networks can learn the parameters $\alpha$ and $\tau$ for the iteration, and how to design an efficient ADR cycle to solve the ADR equation.
\begin{algorithm}
     \caption{Wave-ADR cycle}
     \label{alg:cycle}
     \begin{algorithmic}[1]
     \REQUIRE{$\mathbf{A}$, $\mathbf{g}$, current state $\mathbf{u}$, level $l$, max-level $L$, PARAMETERS $\alpha, \tau$}
     \STATE \textit{Pre-smoothing}
     \IF{$l == 1$}
        \STATE Apply Jacobi iteration: $\mathbf{u} \gets \text{Jacobi}(\mathbf{A}, \mathbf{g}, \mathbf{u}, \omega_0, \text{iter}=1)$
     \ELSIF{$l == L$}
        \STATE Apply Chebyshev semi-iteration: $\mathbf{u} \gets \text{Chebysemi}(\mathbf{A}, \mathbf{g}, \mathbf{u}, \alpha_{L}, \text{iter}=10)$
        \RETURN $\mathbf{u}$
     \ELSIF{$l == 3$}
        \STATE Keep current $\mathbf{u}$
     \ELSE
        \STATE Apply Chebyshev semi-iteration: $\mathbf{u} \gets \text{Chebysemi}(\mathbf{A}, \mathbf{g}, \mathbf{u}, \alpha_{l}, \text{iter}=5)$
     \ENDIF
     \STATE \textit{Coarse grid correction}
     \STATE Compute residual: $\mathbf{r} \gets \mathbf{g} - \mathbf{A} \mathbf{u}$
     \STATE Restrict residual to coarse grid: $\mathbf{r}_c \gets \text{down}(\mathbf{r})$ 
     \STATE Call V-cycle recursively on coarse grid, using re-discretized operator $\mathbf{A}_H$ and zero initial guess $\mathbf{0}$: $\mathbf{e}_c \gets \text{wave\_adr}(\mathbf{A}_H, \mathbf{r}_c, \mathbf{0}, l+1, L, \alpha, \tau)$
     \STATE Interpolate error to fine grid: $\mathbf{e} \gets \text{up}(\mathbf{e}_c)$ 
     \STATE Correct solution: $\mathbf{u} \gets \mathbf{u} + \mathbf{e}$
     \STATE \textit{Post-smoothing}
     \IF{$l == 1$}
        \STATE Apply Jacobi iteration: $\mathbf{u} \gets \text{Jacobi}(\mathbf{A}, \mathbf{g}, \mathbf{u}, \omega, \text{iter}=1)$
     \ELSIF{$l == 2$ ($kh\approx 1$)}
        \STATE Apply Chebyshev semi-iteration: $\mathbf{u} \gets \text{Chebysemi}(\mathbf{A}, \mathbf{g}, \mathbf{u}, \alpha_{l}, \text{iter}=5)$
        \STATE Compute residual: $\mathbf{r} \gets \mathbf{g} - \mathbf{A} \mathbf{u}$
        \STATE Apply ADR correction: $\mathbf{e} \gets \text{adr\_correction}(\mathbf{r}, \tau)$ \COMMENT{Solved by ADR cycle}
        \STATE Correct solution: $\mathbf{u} \gets \mathbf{u} + \mathbf{e}$
     \ELSE
        \STATE Apply Chebyshev semi-iteration: $\mathbf{u} \gets \text{Chebysemi}(\mathbf{A}, \mathbf{g}, \mathbf{u}, \alpha_{l}, \text{iter}=5)$
     \ENDIF
     
     \RETURN Updated state $\mathbf{u}$
     \end{algorithmic}
\end{algorithm}

\subsubsection{Neural networks for learning $\alpha$ and $\tau$}

Recall that a relative optimal choice of $\alpha$ aims to make $\lambda_{\max}/\alpha$ approximate the median of the eigenvalues of $\boldsymbol{A^*\boldsymbol{A}}$, ensuring that the Chebyshev semi-iteration provides effective smoothing. Thus, $\alpha$ is tied to the spectral density of $\boldsymbol{A^*\boldsymbol{A}}$, and consequently to the parameters of the Helmholtz equation. 
The parameters of interest are $\omega$ and $s(\mathbf{x})$. To predict the smoothing parameter $\alpha_l$ for each layer $l$ (where $l = 2,4,\dots,L$), we introduce the convolutional neural network (CNN) that takes $\omega$, $s(\mathbf{x})$, and $N$ as inputs. We chose a CNN due to the simplicity of the task, which involves predicting a single parameter. For this, the classical CNN architecture, composed of convolutional layers, fully connected layers, and ReLU activations, is sufficient.

As previously discussed, an appropriate $\tau$, representing the characteristic error, can be obtained by solving the eikonal equation. However, its numerical computation is also challenging. In this work, we solve eikonal equations using neural operators, which offers two key advantages. 
First, neural operators provide better computational efficiency than traditional methods and are easily implementable on GPUs. Second, the learned $\tau$ is expected to perform better for Wave-ADR-NS than $\tau$ obtained directly from solving the eikonal equation.

Notice that the eikonal solution is independent of the wavenumber, our approach involves training at low resolution and testing at high resolution. To achieve this, we use the Fourier neural operator (FNO)\,\cite{li2020fourier,benitez2023fine}, leveraging its \textit{discretization-invariant} property \cite{kovachki2021neural}. 
It is also important to note that the eikonal solution exhibits non-smoothness near the source, which presents challenges during training.
Following the methodology outlined in\,\cite{treister2019multigrid}, we decompose $\tau$ into $\tau_0\tau_1$, where $\tau_0$ in \cref{eq:eikonal_cons} analytically solves the eikonal equation when $s=1$. 
Subsequently, our objective is to solve the factored eikonal equation for $\tau_1$
\begin{equation}
     \left|\tau_0 \nabla \tau_1+\tau_1 \nabla \tau_0\right|^2=s^2. 
\end{equation}
Since the non-smooth contribution concentrates in $\tau_0$, $\tau_1$ remains relatively smooth near the source.

The input to the FNO consists of four channels: the slowness $s(x,y)$, the real component of the point source $g_r$, and their coordinates $x, y$. The output is the factored travel time $\tau_1(x,y)$
\begin{equation}
     (s, g_r, x, y) \rightarrow \tau_1.  
\end{equation}
We then use FDM (second-order interior, first-order on the boundary) to obtain the derivatives of $\tau$
$$\begin{aligned}\nabla\tau&=\tau_0\nabla\tau_1+\tau_1\nabla\tau_0\\\Delta\tau&=\tau_1\Delta\tau_0+2\nabla\tau_0\cdot\nabla\tau_1+\tau_0\Delta\tau_1.\end{aligned}$$

\subsubsection{ADR solver}
Since the ADR equation \eqref{eq:ADR} satisfies the Péclet condition on the level $k h \approx 1$
\begin{equation}
     2\omega h \max(|\tau_x|,|\tau_y|)\leq 2,  
\end{equation}
the central difference discretization can be utilized without an additional stabilizing term.
Nevertheless, when employing MG method to solve the discretized system, the ADR equations on coarse-grids fail to satisfy the Péclet condition, resulting in diminishing \textit{h-ellipticity} and impeding the effectiveness of smoothers\,\cite{trottenberg2000multigrid}. 
Indeed, when employing the central difference discretization method, the difficulty of solving the ADR equation mirrors that of the original Helmholtz equation \cite{gander2024analysis}.
Hence, we employ first-order upwind discretization method instead 
\begin{equation}
    \left(\frac{\partial \tau}{\partial x} \cdot \frac{\partial a}{\partial x}\right)_j \approx \begin{cases}
            \left(\frac{\partial \tau}{\partial x}\right)_j \cdot \frac{a_j-a_{j-1}}{ h}, & \text{if } \left(\frac{\partial \tau}{\partial x}\right)_j>0 \\
            \left(\frac{\partial \tau}{\partial x}\right)_j \cdot \frac{a_{j+1}-a_{j}}{h}, & \text{if } \left(\frac{\partial \tau}{\partial x}\right)_j<0.
            \end{cases}    
        \label{eq:upwind}
\end{equation}

\begin{figure}[!htb]
     \centering
     \subfigure[]{\label{subfig:discremethod}\includegraphics[width=0.4\textwidth]{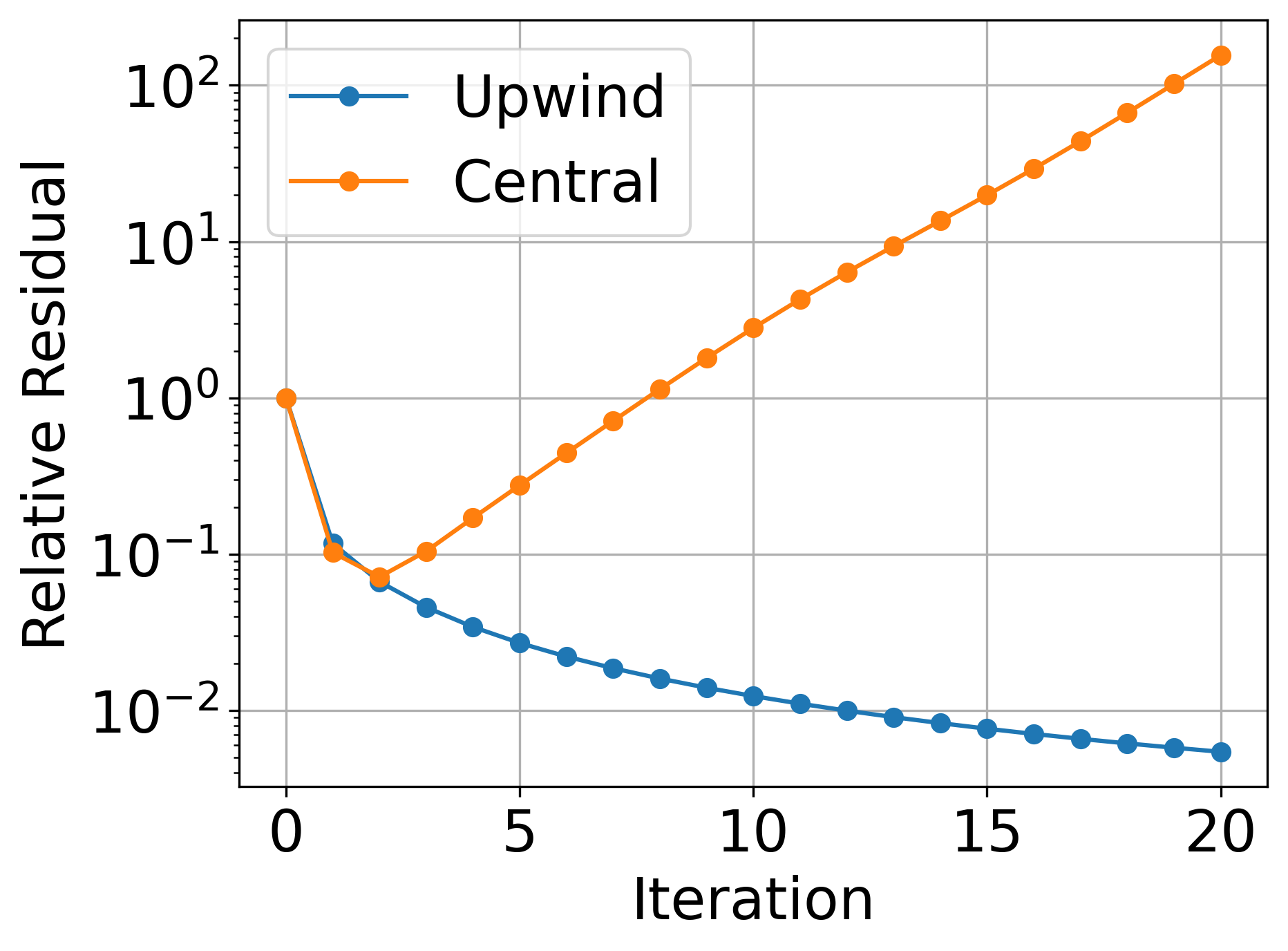}}\quad\quad
     \subfigure[]{\label{subfig:difftol}\includegraphics[width=0.4\textwidth]{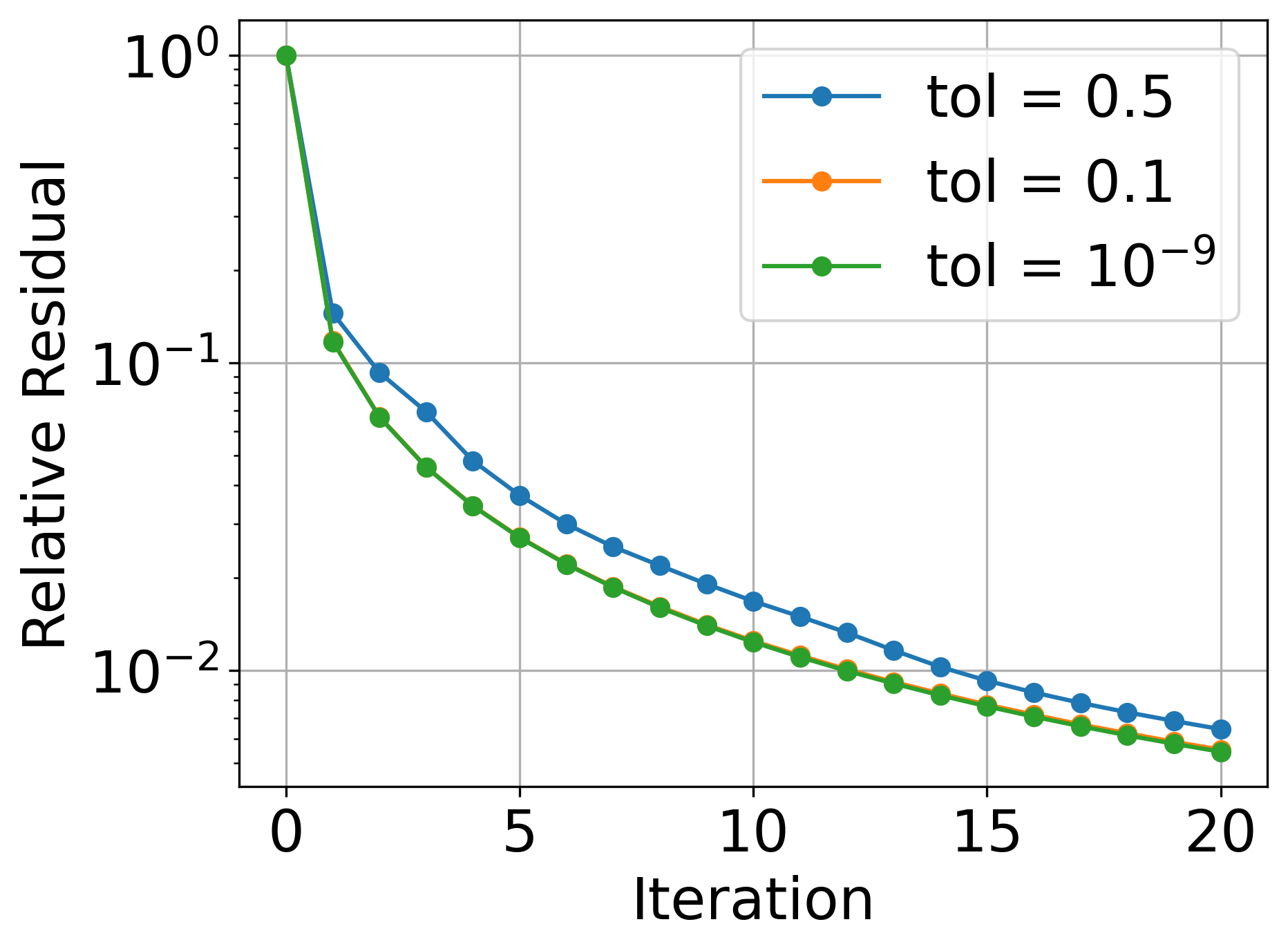}}
     \caption{(a) The effect of different discretization methods for ADR equation on the convergence rate of Wave-ADR-NS. (b) The influence of different solving accuracy for ADR upwind discrete system on the convergence rate of Wave-ADR-NS.}
\end{figure}

Another advantage of the first-order upwind discretization method is its better stability compared to the central difference method when used for correction.
\cref{subfig:discremethod} illustrates the convergence history of Wave-ADR-NS for the constant Helmholtz equation with $N=128$ and $\omega=20\pi$, employing different ADR discretization methods. The linear ADR system is solved using a direct method.
It can be observed that  Wave-ADR-NS diverges when the central difference method is used without stability.

Next, we analyze the influence of different solution accuracy for ADR discrete system on the convergence rate of Wave-ADR-NS. 
\Cref{subfig:difftol} illustrates the convergence history of Wave-ADR-NS under various ADR solution accuracy. These results demonstrate that once the solution accuracy of the ADR linear system surpasses $0.1$, further increasing accuracy hardly improves the performance of Wave-ADR-NS. Therefore, we employ an ADR solver with the accuracy of $0.1$.
Furthermore, results from ablation experiments demonstrate that a single MG V-cycle with the GMRES(3) smoother can reduce the relative residual of the ADR linear system to below $0.1$. Therefore, we utilize it as the ADR solver.

\subsection{Training procedures and implementation}\label{sec:wavesolver}
We have introduced the iteration algorithm and network structure. Now we generate the training data and give the loss function to complete the learning of parameters $\alpha$ and $\tau$ in the Wave-ADR-NS. 
We aim to adopt an unsupervised training approach, where the training data consists of triplets $(\omega_i, s_i, g_i)$. The slowness model $s_i \in \mathbb{R}^{N \times N}$ is sampled from existing datasets, with the size $N$ determined by different angular frequencies $\omega_i$ to satisfy $\omega h \approx 0.5$, as detailed in the next section. Here, $g_i$ represents a point source centered at $\Omega$. 
For each triplet $(\omega_i, s_i, g_i)$, we first use the CNN and FNO to obtain $\alpha$ and $\tau$, respectively, and then perform Wave-ADR cycle $K$ times. The relative residual at this stage serves as the loss for training
\begin{equation}
     \mathcal{L} = \frac{1}{N_{\mathrm{batch}}} \sum_{i=1}^{N_{\mathrm{batch}}} \frac{||\mathbf{g}_{i} - \mathbf{A}_{i} \mathbf{u}^{(K)}_{i}||_{2}^{2}}{||\mathbf{g}_{i}||_{2}^{2}}.
     \label{eq:loss}
\end{equation}
Here, $\mathbf{A}_{i}$ is determined by $(\omega_i, s_i)$, $\mathbf{g}_{i}$ represents a discrete point source located at the center of $\Omega$, and $\mathbf{u}^{(K)}_{i}$ is the solution obtained at the $K$-th Wave-ADR iteration, $K=3$ is used in our experiments. The training process is summarized in \cref{alg:train}.
\begin{algorithm}
     \caption{Wave-ADR-NS offline training}
     \label{alg:train}
     \begin{algorithmic}[1]
     \REQUIRE{Training data $\{(\omega_i, s_i, g_i)\}_{i=1}^{N_{\mathrm{train}}}$, batch size $N_{\mathrm{batch}}$, learning rate $\eta$, number of epochs $N_{\mathrm{epochs}}$}
     \FOR{epoch = 1, ..., $N_{\mathrm{epochs}}$}
          \STATE Shuffle training data
          \FOR{each batch $\{(\omega_i, s_i, g_i)\}_{i=1}^{N_{\mathrm{batch}}}$}
               \FOR{each training sample $(\omega_i, s_i, g_i)$ in the batch} 
                    \STATE Forward pass through CNN: $\alpha_i \gets \text{CNN}(\omega_i, s_i)$
                    \STATE Forward pass through FNO: $\tau_i \gets \text{FNO}(s_i, g_i)$
                    \STATE Initialize solution: $\mathbf{u}^{(0)} = \mathbf{0}$ 
                    \FOR{$j = 1, \cdots K$}
                         \STATE Call Wave-ADR cycle (\cref{alg:cycle}) with inputs $(\mathbf{A}_i, \mathbf{g}_i, \mathbf{u}^{(j-1)}, l=1, L, \alpha_i, \tau_i)$ to compute $\mathbf{u}^{(j)}$
                    \ENDFOR
                    \STATE Compute relative residual (loss): 
                    $\text{loss}_i = \frac{\|\mathbf{g}_i - \mathbf{A}_i \mathbf{u}^{(K)} \|_2^2}{\|\mathbf{g}_i\|_2^2}$   
               \ENDFOR
               \STATE Compute total loss for the batch: 
               $
               \mathcal{L} = \frac{1}{N_{\mathrm{batch}}} \sum_{i=1}^{N_{\mathrm{batch}}} \text{loss}_i
               $
               \STATE Backpropagate to compute gradients: $\nabla_{\text{CNN}}, \nabla_{\text{FNO}} \gets \nabla_{\mathcal{L}}$
               \STATE Update parameters of CNN and FNO using Adam optimizer \cite{kingma2014adam}
          \ENDFOR
     \ENDFOR
     \RETURN Trained models CNN and FNO
     \end{algorithmic}
\end{algorithm}

\begin{figure}[!htb]
     \centering
     \includegraphics[width=\textwidth]{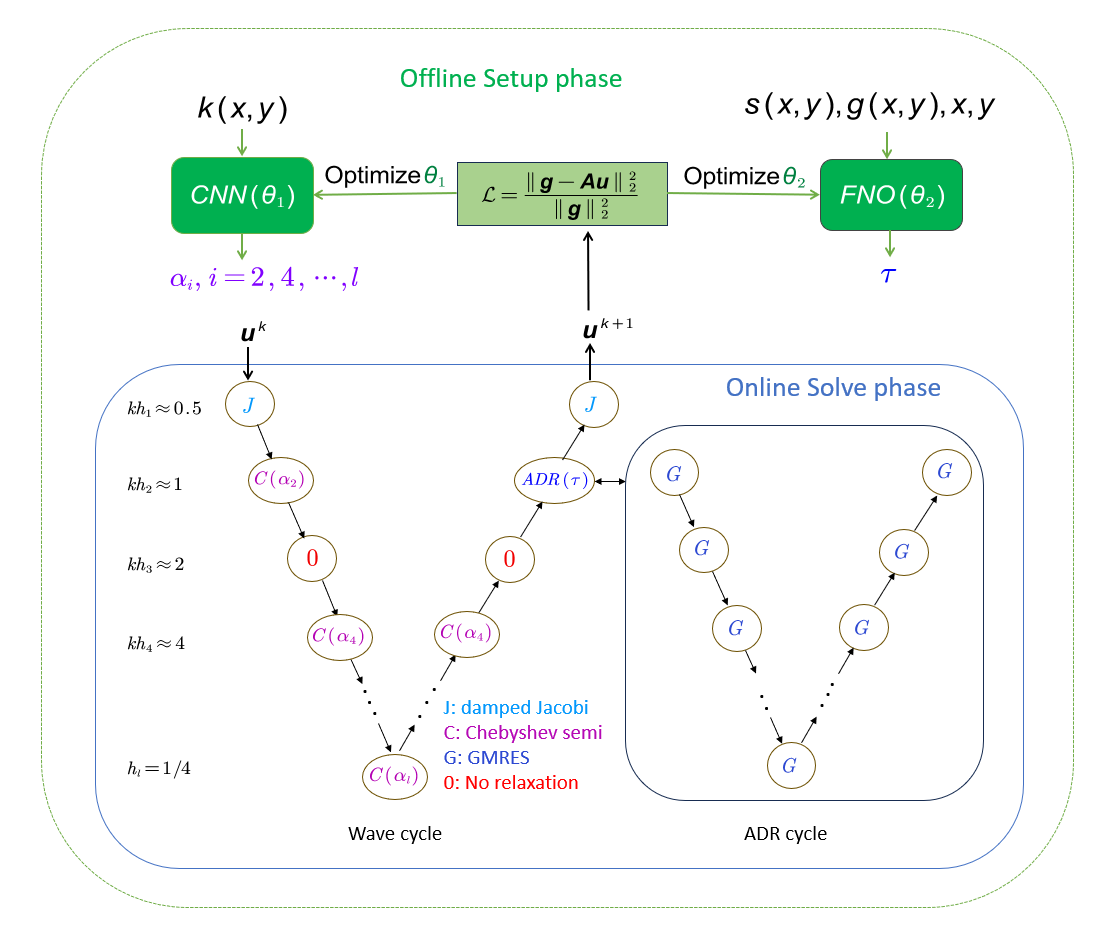}
     \caption{Diagram of the Wave-ADR-NS: The dashed green box represents the setup phase, during which two neural networks, CNN and FNO, are trained. The blue box denotes the solve phase, where different relaxation methods are executed at each level to address errors with varying frequencies. At the second level, the ADR cycle is employed to eliminate the near-null space error components.}
     \label{fig:cycle}
\end{figure} 
\Cref{fig:cycle} summarizes the diagram of the Wave-ADR-NS, which comprises two main phases: the \textit{setup phase} and the \textit{solve phase}. 
During the offline setup phase, a CNN and a Fourier neural operator (FNO)\,\cite{li2020fourier,benitez2023fine} are employed to learn the parameters $\alpha$ and $\tau$ in the solver, respectively.
The online solve phase consists of two V-cycles: the wave cycle and ADR cycle, which aim to eliminate non-characteristic and characteristic error components, respectively.

\begin{figure}[!htb]
     \centering
     \includegraphics[width=\textwidth]{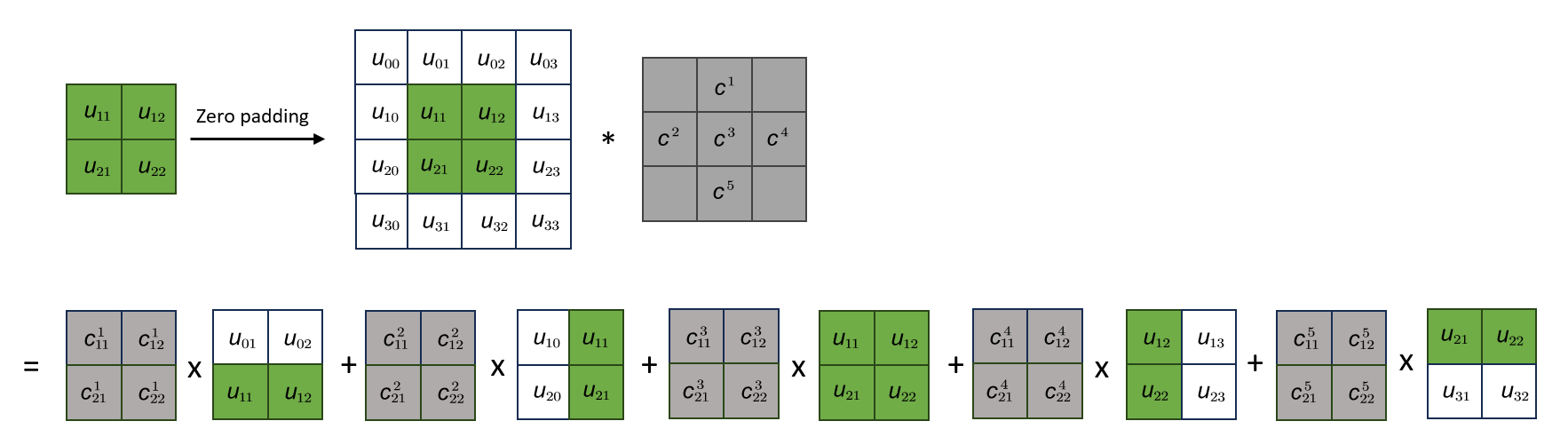}
     \caption{Matrix-vector product with variable stencil.}
     \label{fig:conv}
\end{figure}

To avoid storing the sparse matrix $\mathbf{A}$, we implement matrix-vector multiplication in convolutional form using the PyTorch deep learning framework \cite{Paszke2019PyTorchAI}. The second-order FD stencil for the Laplace operator in the Helmholtz and ADR equations is given by
\begin{equation}\nonumber
    \frac{1}{h^2}\left[\begin{array}{ccc}
      & -1 &  \\
    -1 & 4 & -1 \\
      & -1 & 
    \end{array}\right],
\end{equation}
and convolutions with this kernel and a stride of one are performed. 
The stencil for the first-order upwind discretization of the convection term $\nabla\tau \cdot \nabla a$ is
\begin{equation}\nonumber
    \frac{1}{2h}\left[\begin{array}{ccc}
      & -\tau_y-|\tau_y| & \\
    -\tau_x-|\tau_x| & 2(|\tau_x|+|\tau_y|) & \tau_x-|\tau_x| \\
      & \tau_y-|\tau_y| &
    \end{array}\right] := \left[\begin{array}{ccc}
     & c^1 &  \\
   c^2 & c^3 & c^4 \\
     & c^5 & 
   \end{array}\right],
\end{equation}
where $\nabla\tau=(\tau_x, \tau_y)$.
As each entry $c^j$ depends on the coordinate position, it cannot be implemented using the usual convolution. 
Instead, we implement the matrix-vector product by: (1) multiplying each $c^j$ with the appropriate vector padding and slice, and (2) summing the five channel outputs, as shown in \cref{fig:conv}.

The MG restriction and interpolation utilize standard full-weighting and bilinear interpolation operators
\begin{equation}\nonumber
    K_{I_h^H}=\frac{1}{16}\left[\begin{array}{lll}
    1 & 2 & 1 \\
    2 & 4 & 2 \\
    1 & 2 & 1
    \end{array}\right], \quad K_{I_H^h}=\frac{1}{4}\left[\begin{array}{lll}
    1 & 2 & 1 \\
    2 & 4 & 2 \\
    1 & 2 & 1
    \end{array}\right],
\end{equation}
convolution and deconvolution operations using these kernels with stride two are performed, respectively.

Inner iterative methods such as the damped Jacobi method, Chebyshev semi-iteration, and GMRES rely solely on matrix-vector products, allowing for straightforward implementation. This ensures Wave-ADR functions as a differentiable solver, enabling proper backpropagation during the training process.
Unlike recently developed nonlinear preconditioners \cite{azulay2022multigrid,kopanivcakova2024deeponet,rudikov2024neural} or hybrid iterative methods \cite{zhang2022hybrid,hu2024hybrid}, which often decouple training and testing phases and may introduce uncertainties like divergence when applying the preconditioner in testing, Wave-ADR-NS creates a \textit{closed loop} between training and testing. This allows the training loss to directly reflect the convergence rate.
In summary, Wave-ADR-NS enables matrix-free computation, batch processing, and GPU acceleration, making it well-suited for gradient-based optimization while significantly reducing storage and memory demands compared to sparse matrix representations. Code and data for the following experiments are available at \url{https://github.com/cuichen1996/Wave-ADR-NS}.

\section{Numerical experiments}\label{sec:04}
In this section, we first train Wave-ADR-NS on velocity fields transformed from the CIFAR-10 dataset \cite{krizhevsky2009learning}, which has also been used in other deep learning-augmented MG preconditioners \cite{azulay2022multigrid,lerer2023multigrid}. 
Next, we illustrate the learned the smoother parameters $\alpha$ and the phase function $\tau$.
Finally, we evaluate the in-distribution generalization and out-of-distribution transfer capabilities of Wave-ADR-NS, comparing it with other MG methods. 
All experiments are conducted on an Nvidia A100-SXM4-80GB GPU.

\subsection{Dataset}\label{sec:kappa}
The CIFAR-10 dataset serves as a widely used benchmark for image classification, comprising 50,000 images with $32\times 32$ pixels distributed across $10$ classes. 
This diverse dataset enables a comprehensive testing of our method across various scenarios.
To convert natural images into slowness models $s(x,y)$ for Helmholtz equation, we take the following three steps following \cite{azulay2022multigrid}
\begin{enumerate}
     \item resize the images to match the grid size through bilinear interpolation;
     \item apply Gaussian smoothing on the image;
     \item normalize the values to fall within the range of $[0.25, 1]$.
\end{enumerate}

\begin{figure}[!htb]
     \centering
     \subfigure[$N=128$]{
     \includegraphics[width=0.25\textwidth]{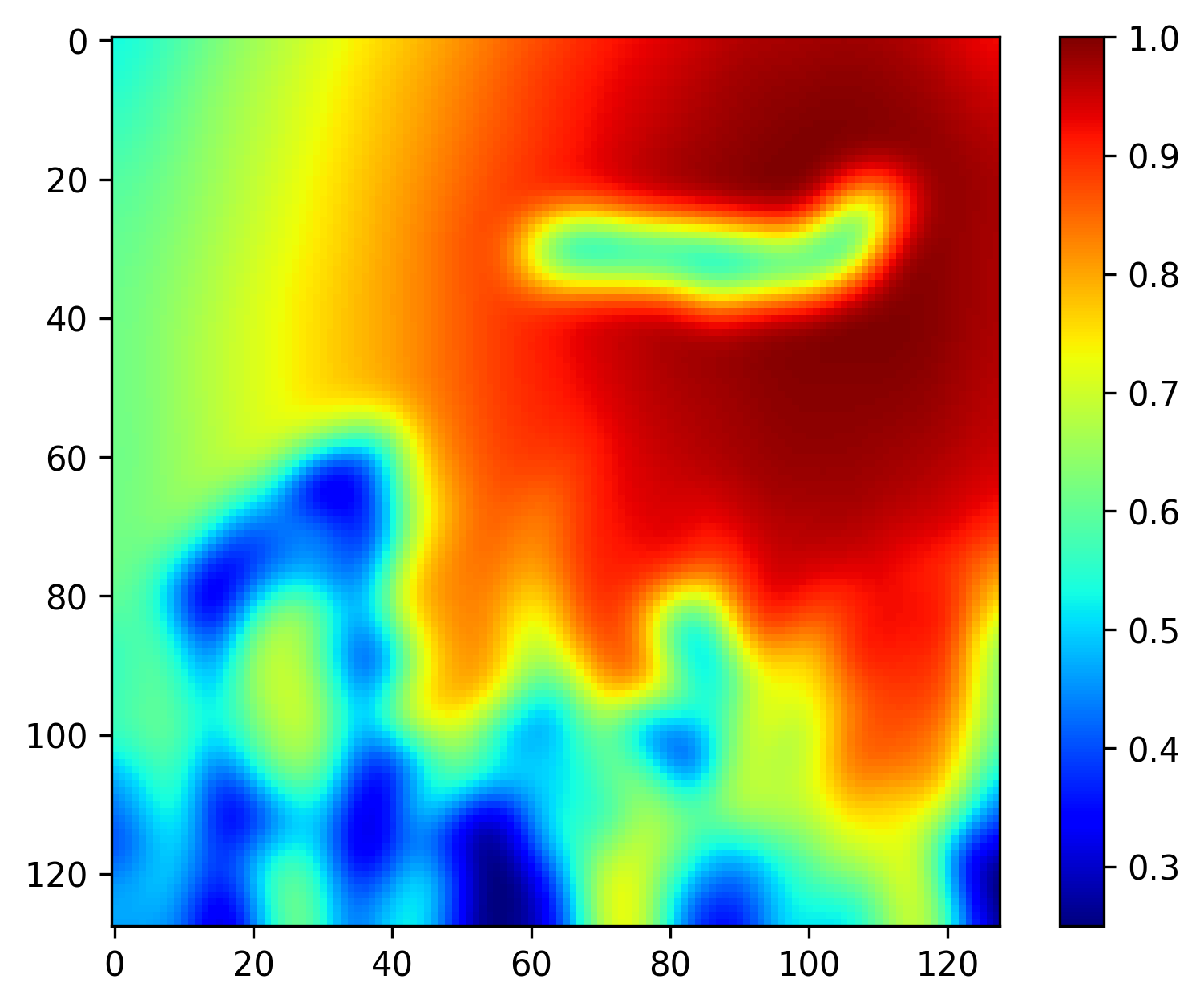}}
     \subfigure[$N=256$]{
     \includegraphics[width=0.25\textwidth]{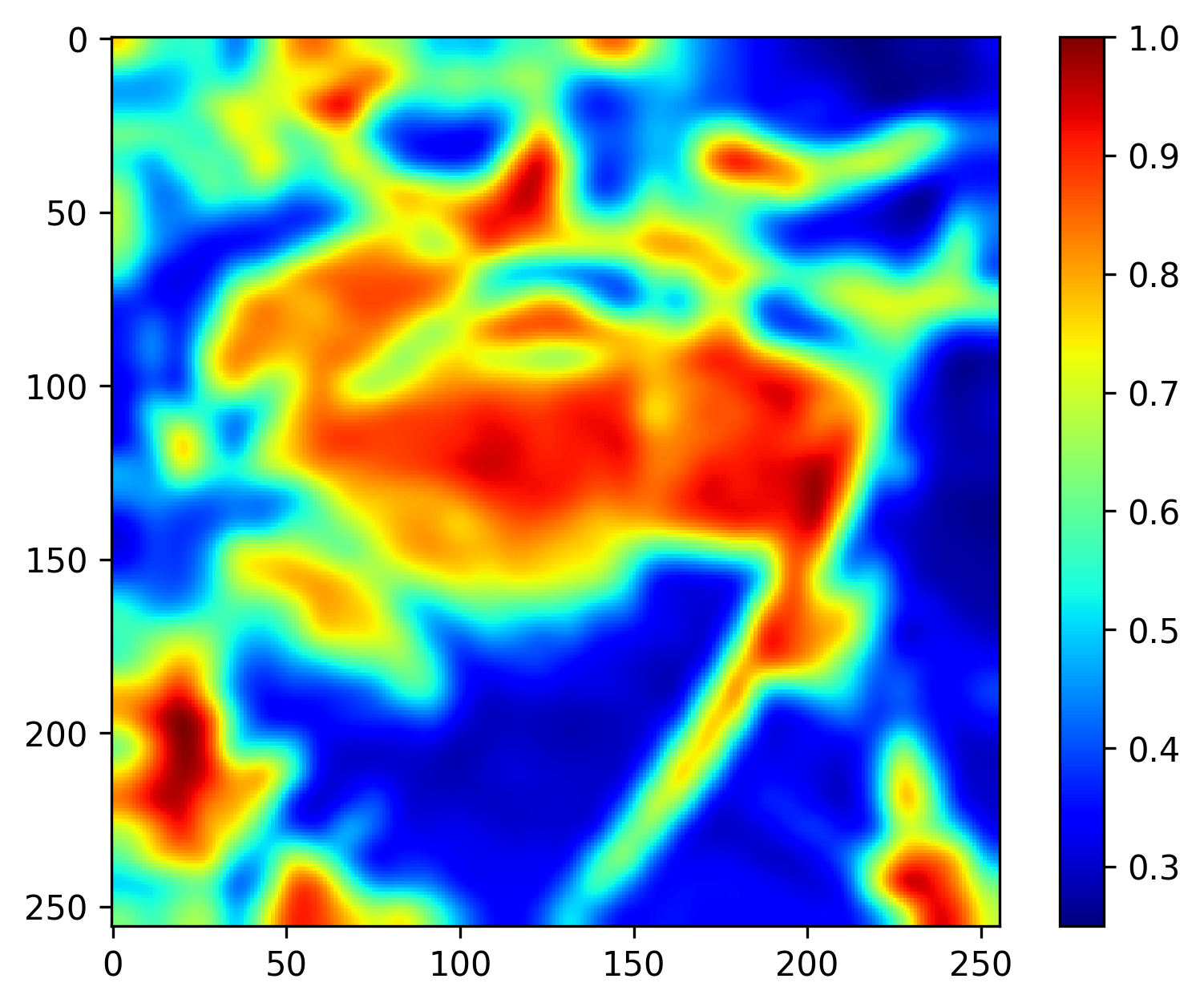}}
     \subfigure[$N=512$]{
     \includegraphics[width=0.25\textwidth]{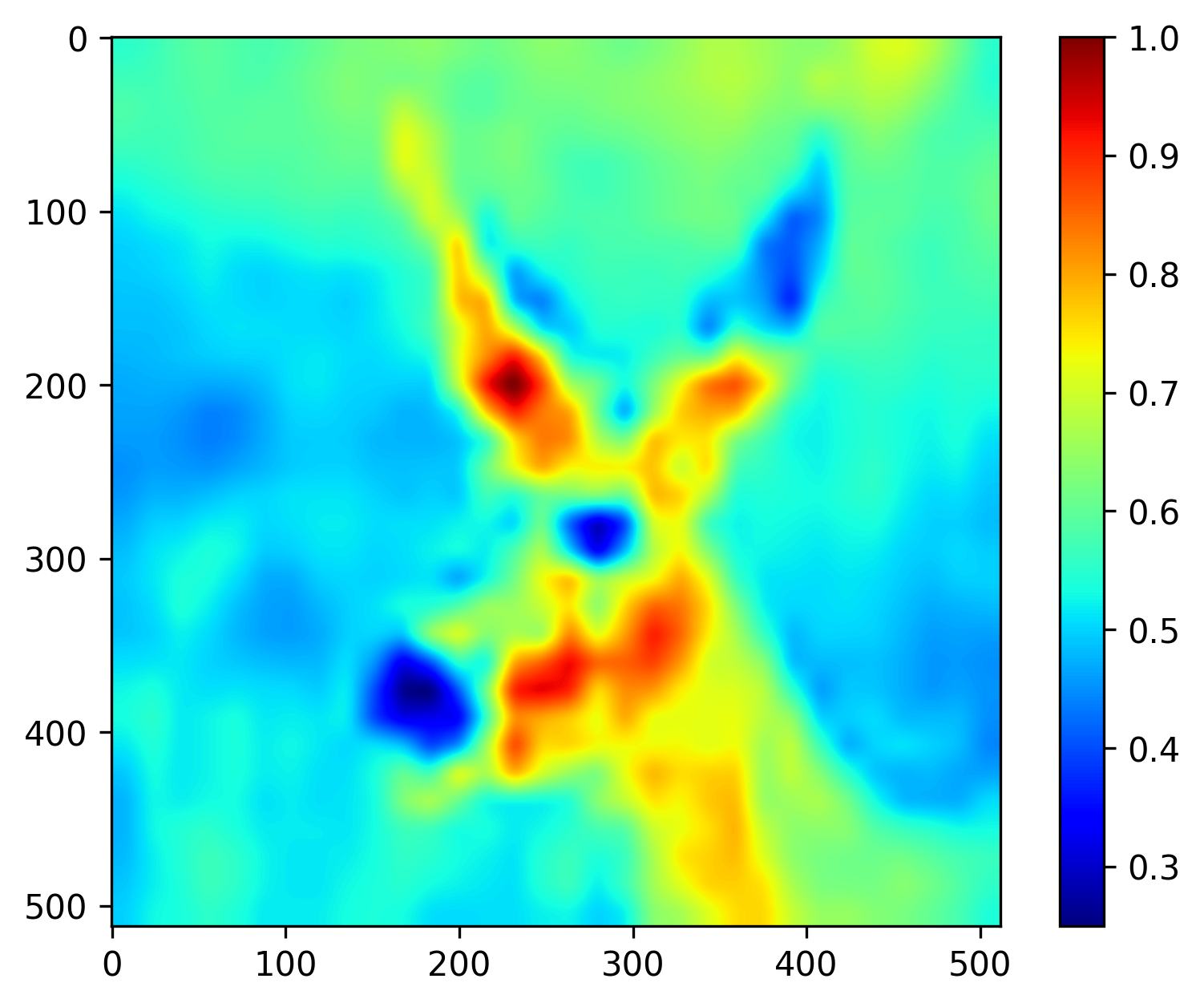}}
     \caption{Examples of slowness models converted from the CIFAR-10 dataset with different resolutions}
     \label{fig:example}
\end{figure}

\cref{fig:example} displays three randomly generated slowness models from the CIFAR-10 dataset, with grid sizes of $128\times 128$, $256\times 256$, $512\times 512$, respectively. 
In our experiments, we select a total of $16000$, $10000$, and $6000$ slowness models of sizes $128\times 128$, $256\times 256$, and $512\times 512$ as the training set, following the approach outlined in\,\cite{lerer2023multigrid}.
The testing set is re-selected to include slowness model with grid size of $4096\times 4096$.

\begin{figure}[!htbp]
     \centering
     \subfigure[Slowness model example]{\label{fig:kappa_test}\includegraphics[width=0.3\textwidth]{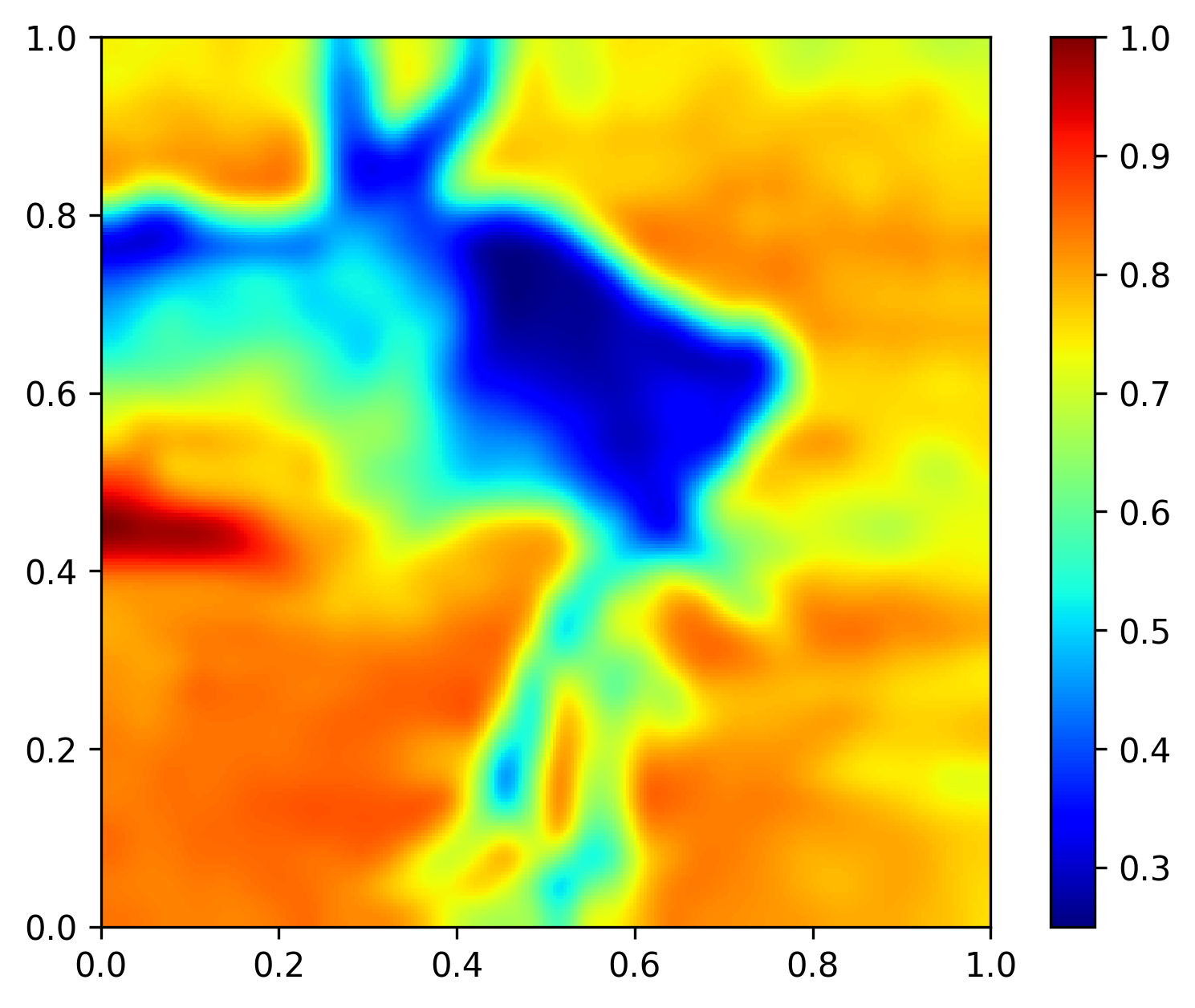}}\quad
     \subfigure[Helmholtz solution]{\label{subfig:helmsol}\includegraphics[width=0.3\textwidth]{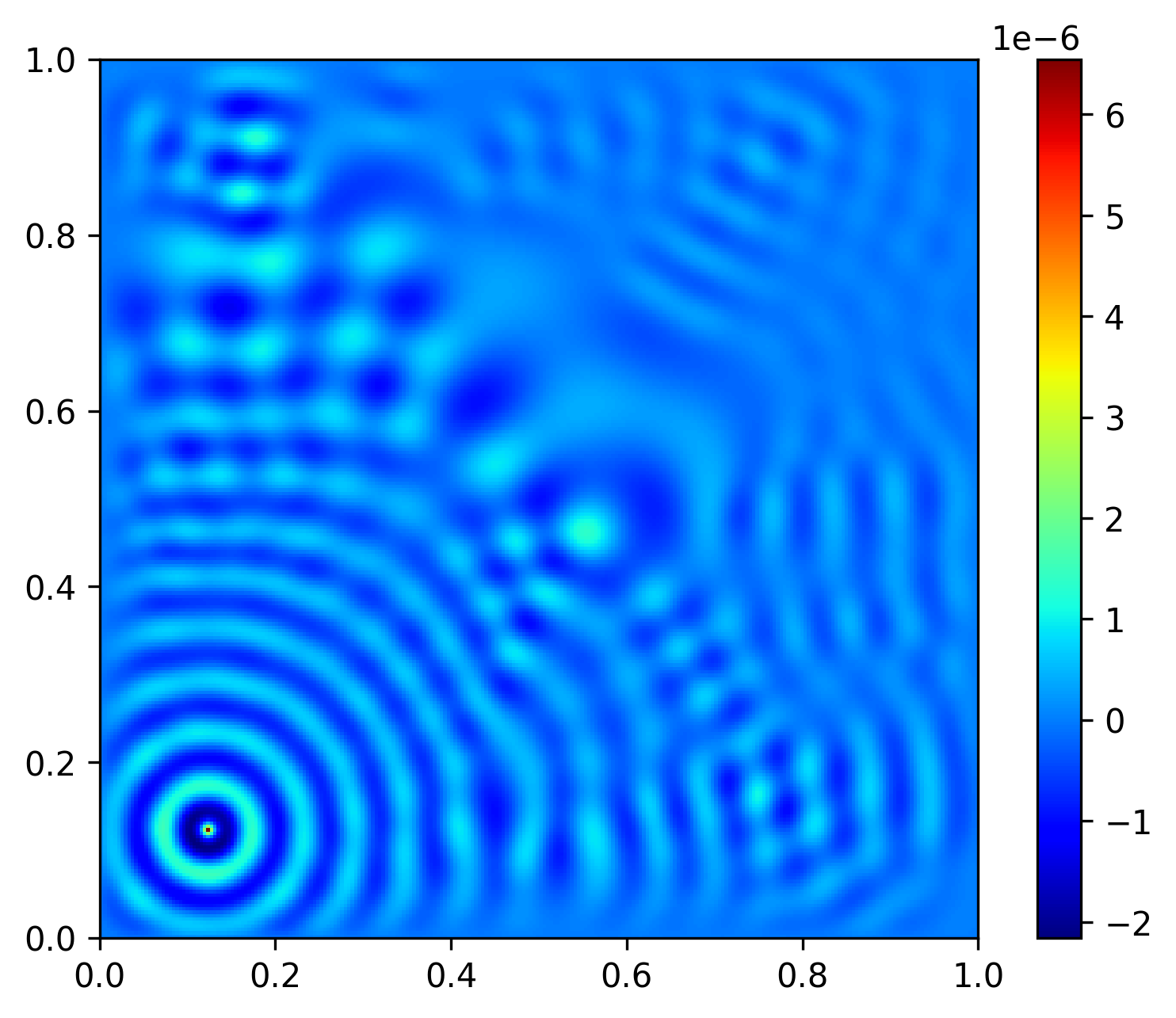}}
     \subfigure[$|\mathcal{F}(e)|$]{\label{subfig:what}\includegraphics[width=0.3\textwidth]{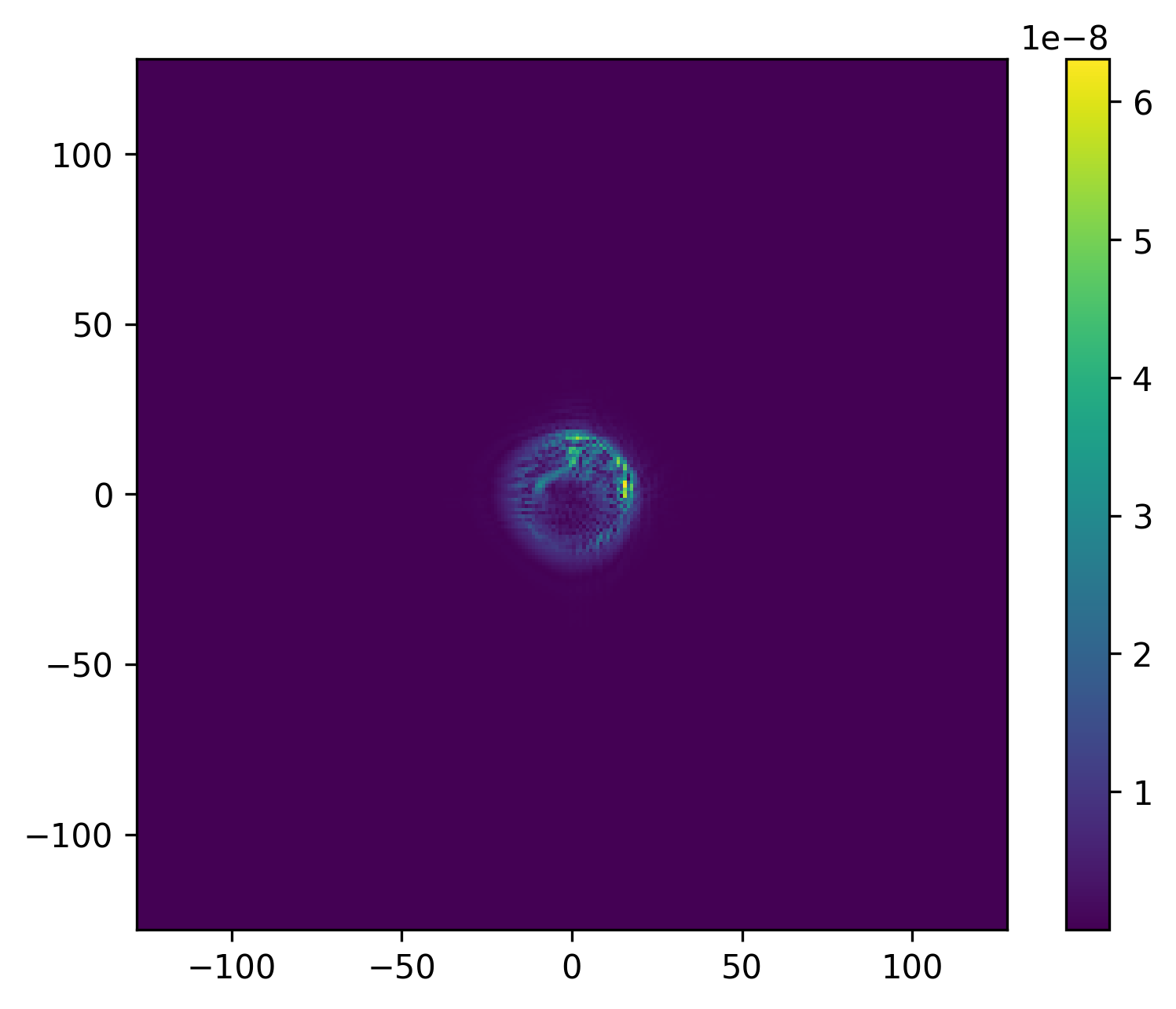}}\\
     \subfigure[$\tau$]{\label{subfig:tau}\includegraphics[width=0.27\textwidth]{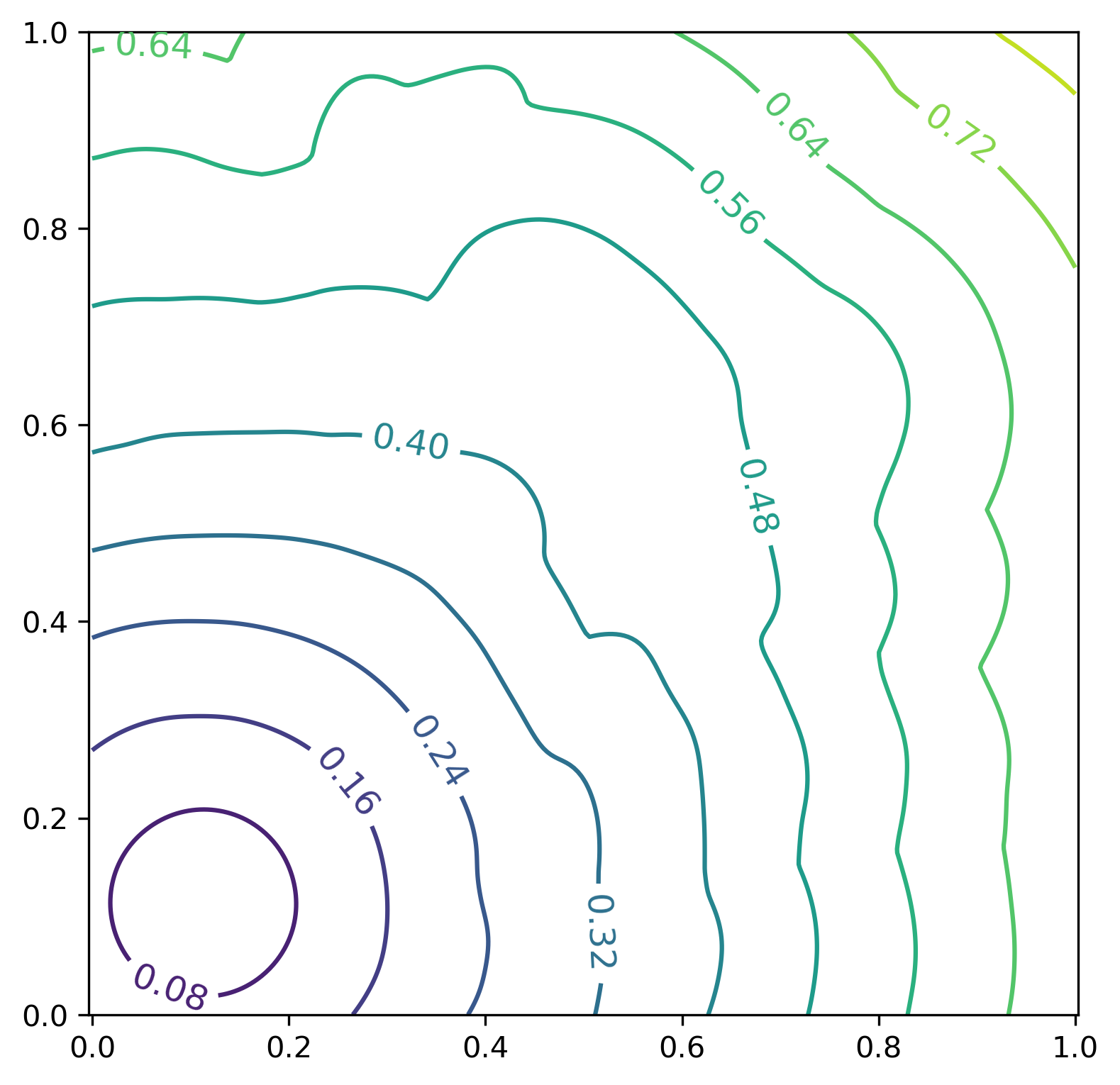}}\quad\quad
     \subfigure[$\nabla \tau$]{\label{subfig:gradtau}\includegraphics[width=0.27\textwidth]{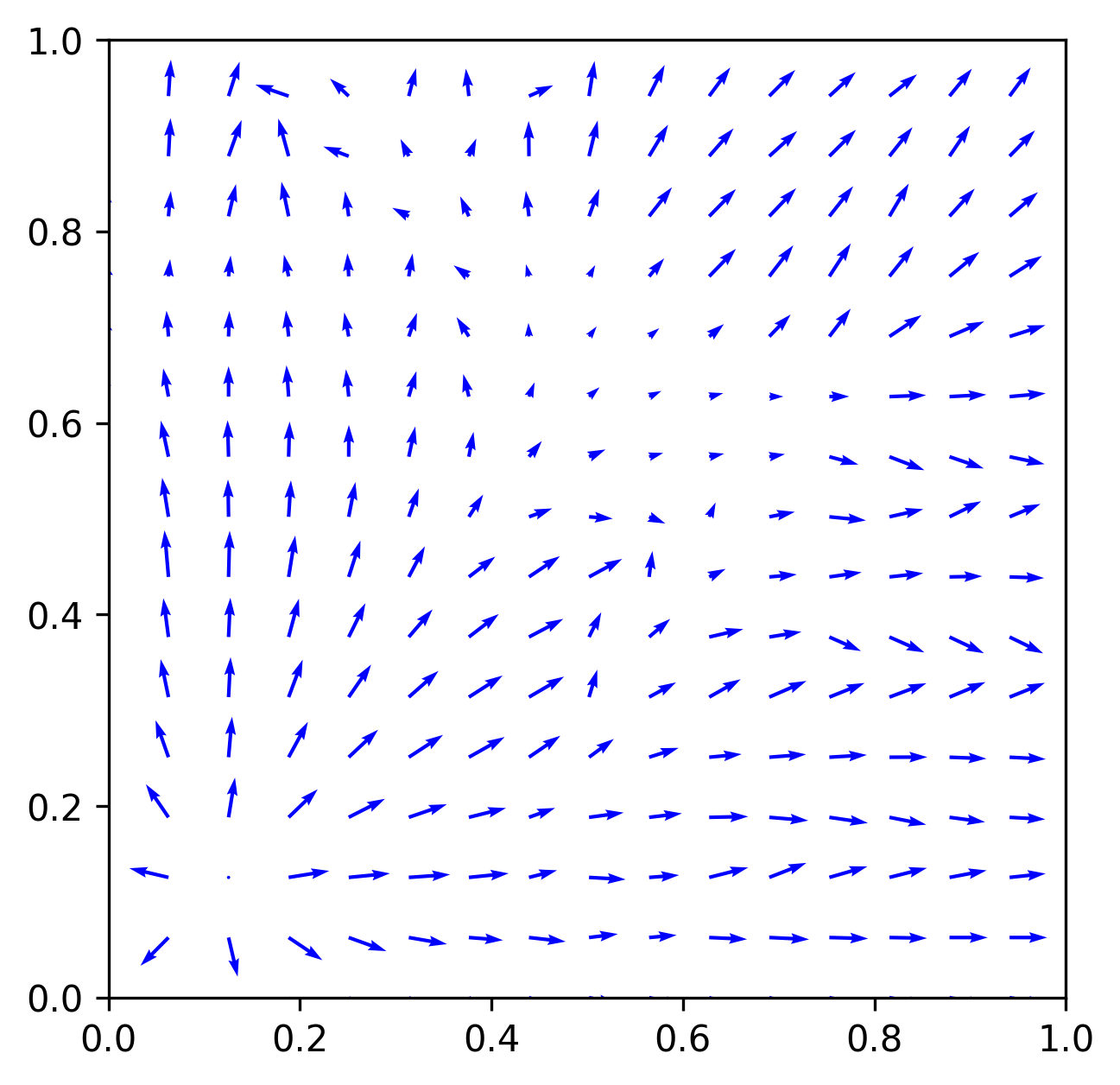}}\quad
     \subfigure[$|\mathcal{F}(e^{-i\omega \tau})|$]{\label{subfig:taufft}\includegraphics[width=0.32\textwidth]{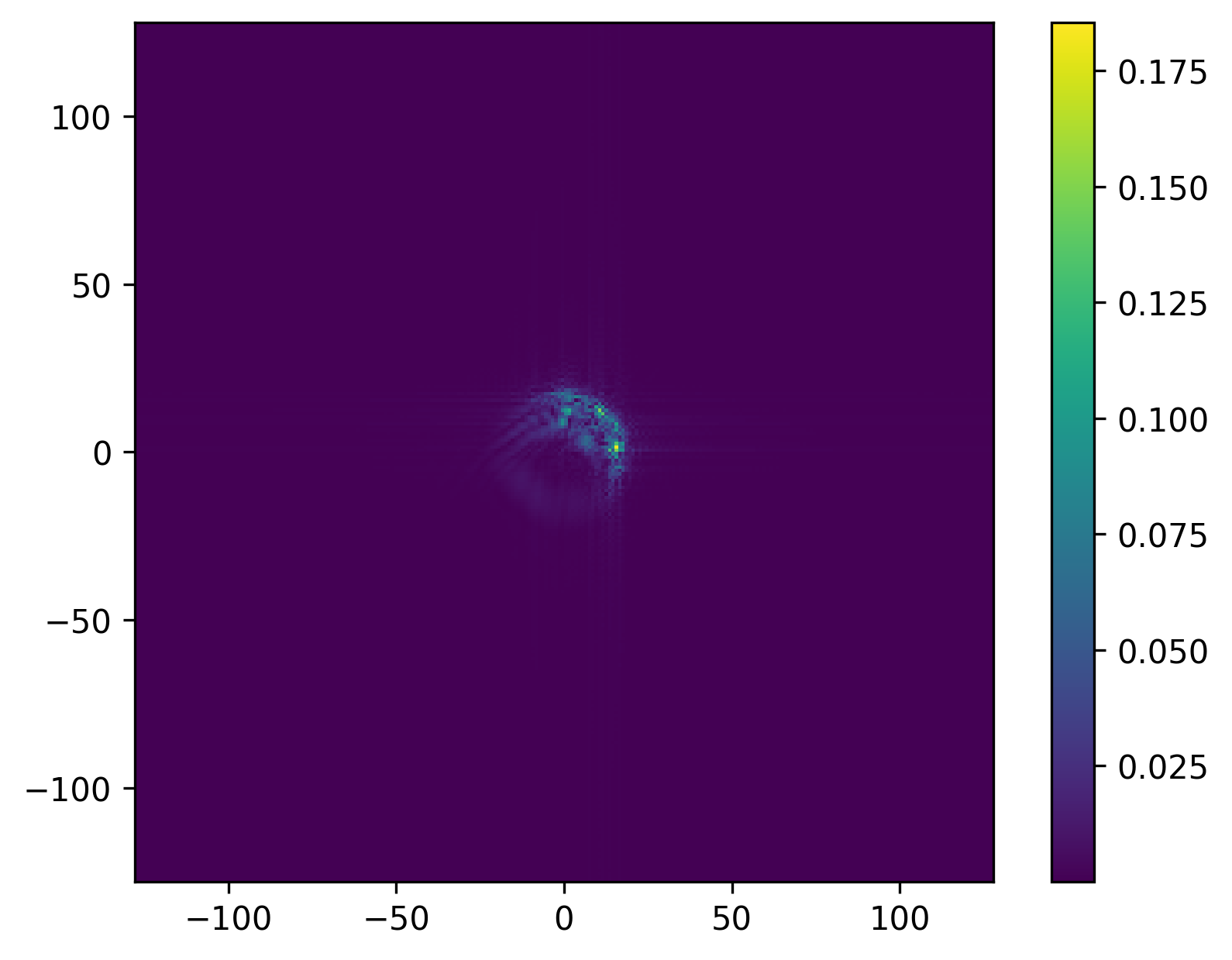}}
     \caption{(a) Example of a slowness model from CIFAR-10 dataset, (b) Solution of the heterogeneous Helmholtz equation with a point source when $\omega=40\pi$ and $N=256$, (c) The module of Fourier transform of the error after one wave cycle iteration, (d) Solution of  the eikonal equation, (e) Gradient field of $\tau$, (f) The module of Fourier transform of $e^{-i\omega\tau}$.}
     \label{fig:7}
\end{figure}

In \cref{fig:7}, a slowness model is presented as an example to illustrate the solution of the corresponding heterogeneous Helmholtz equation (see \cref{subfig:helmsol}), the error in Fourier space between one wave cycle iterative solution and the reference solution (see \cref{subfig:what}), the phase function $\tau$ (see \cref{subfig:tau}) with its gradient field (see \cref{subfig:gradtau}) obtained by solving the eikonal equation, and the distribution of $e^{-i\omega\tau}$ in Fourier space (see \cref{subfig:taufft}).
It can be observed that the gradient of $\tau$ can reflect the propagation direction of the wavefield, and $e^{-i\omega\tau}$ represents the distribution of the characteristic components, even for variable wavenumber problems.

\subsection{Learning parameters $\alpha$ and $\tau$}\label{sec:422}
In this subsection,  we validate the effect of the learned $\alpha$ and $\tau$.

\begin{figure}[!htbp]
     \centering
     \subfigure[$N=128, \omega=20\pi$]{\includegraphics[width=0.3\textwidth]{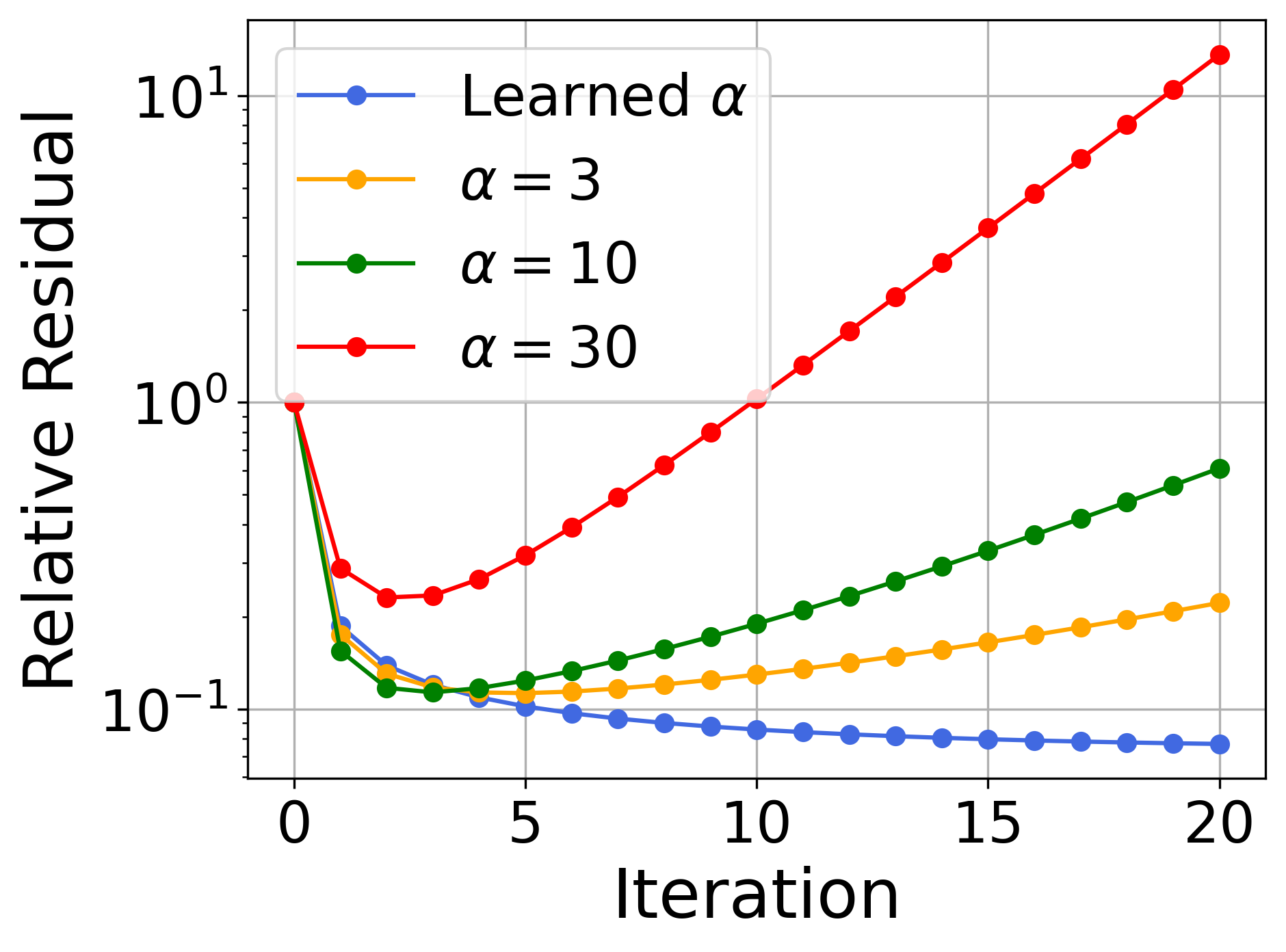}}
     \subfigure[$N=256, \omega=40\pi$]{\includegraphics[width=0.3\textwidth]{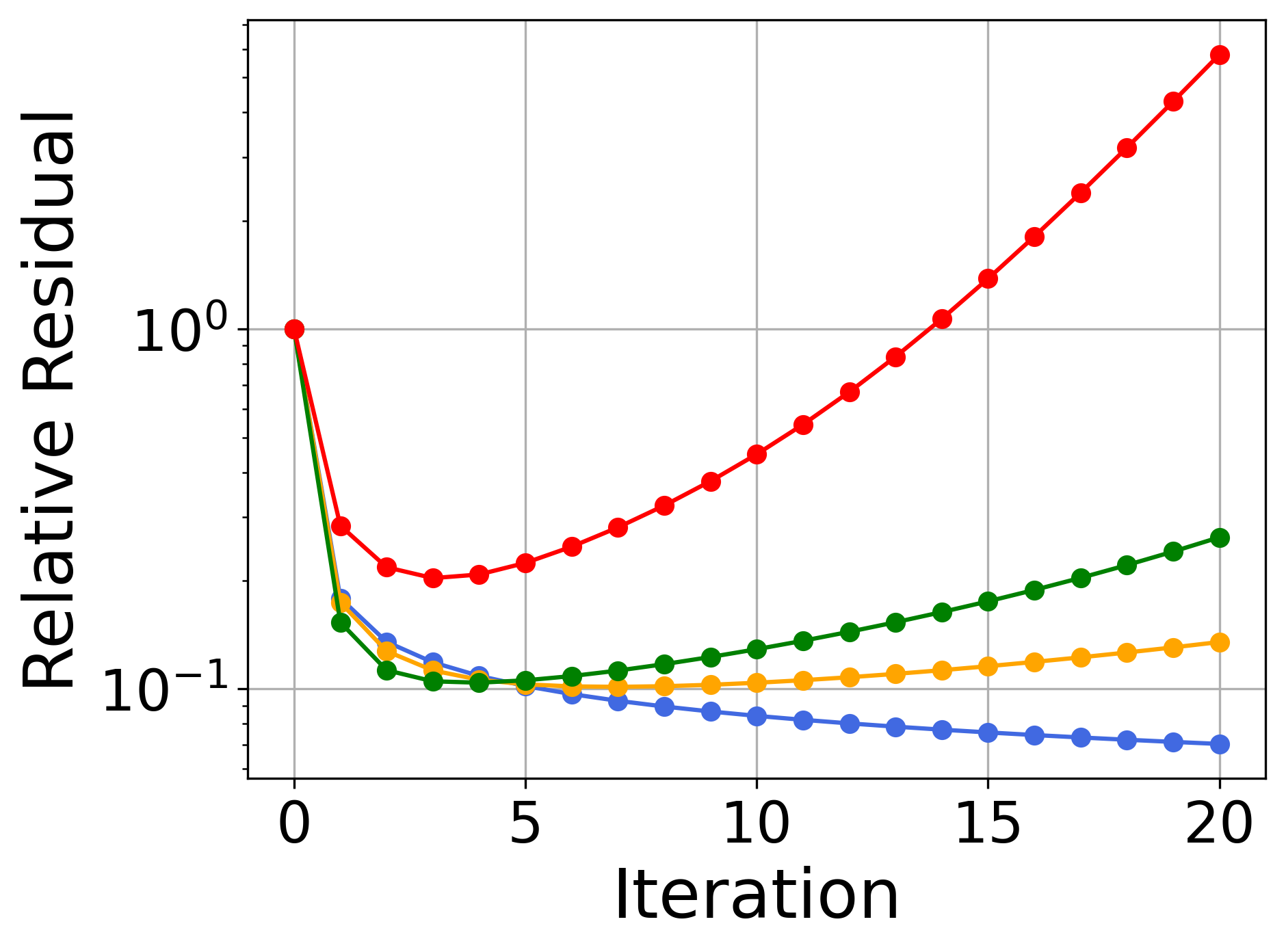}}
     \subfigure[$N=512, \omega=80\pi$]{\includegraphics[width=0.3\textwidth]{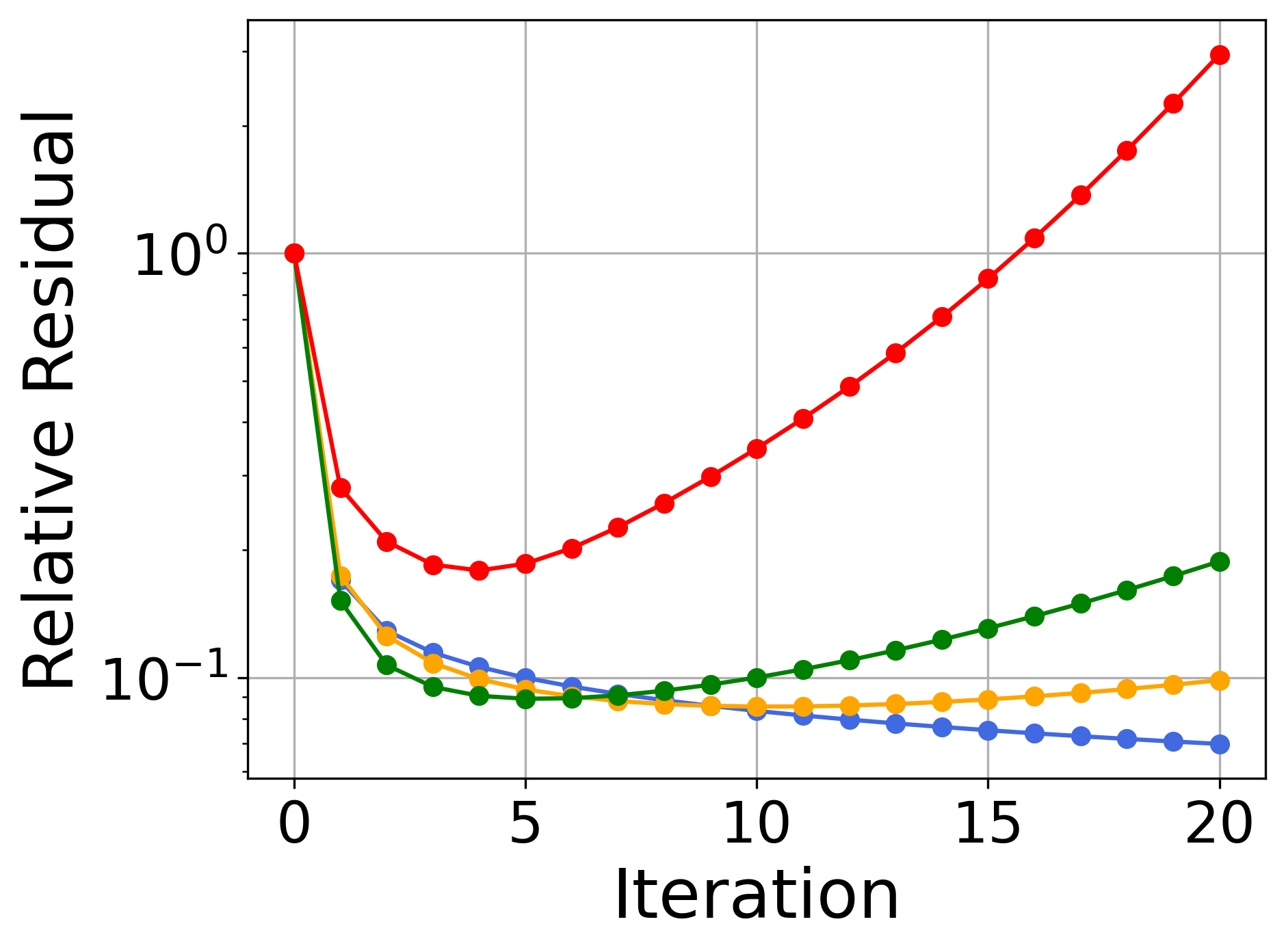}}
     \caption{Convergence history of wave cycle with different values of $\alpha$ when solving the point source Helmholtz equation with different wavenumbers.}
     \label{fig:alpha}
\end{figure}

We have introduced the importance of selecting the parameter $\alpha$ of the Chebyshev semi-iteration method in \Cref{sec:03}. 
In prior applications, $\alpha$ commonly adopts default parameters, typically $3$, $10$, or $30$\,\cite{adams2003parallel}.
In Wave-ADR-NS, the parameter $\alpha$ is alternatively determined by the CNN introduced before.
A comparison of the performance between the learned $\alpha$ and the default value is conducted.
\cref{fig:alpha} illustrates the convergence history of the wave solver for solving the point source Helmholtz equation at different frequencies with the slowness model shown in \cref{fig:kappa_test}. It is evident that the learned $\alpha$ outperforms these default values.
 
\begin{figure}[!htb]
     \centering
     \includegraphics[width=0.4\textwidth]{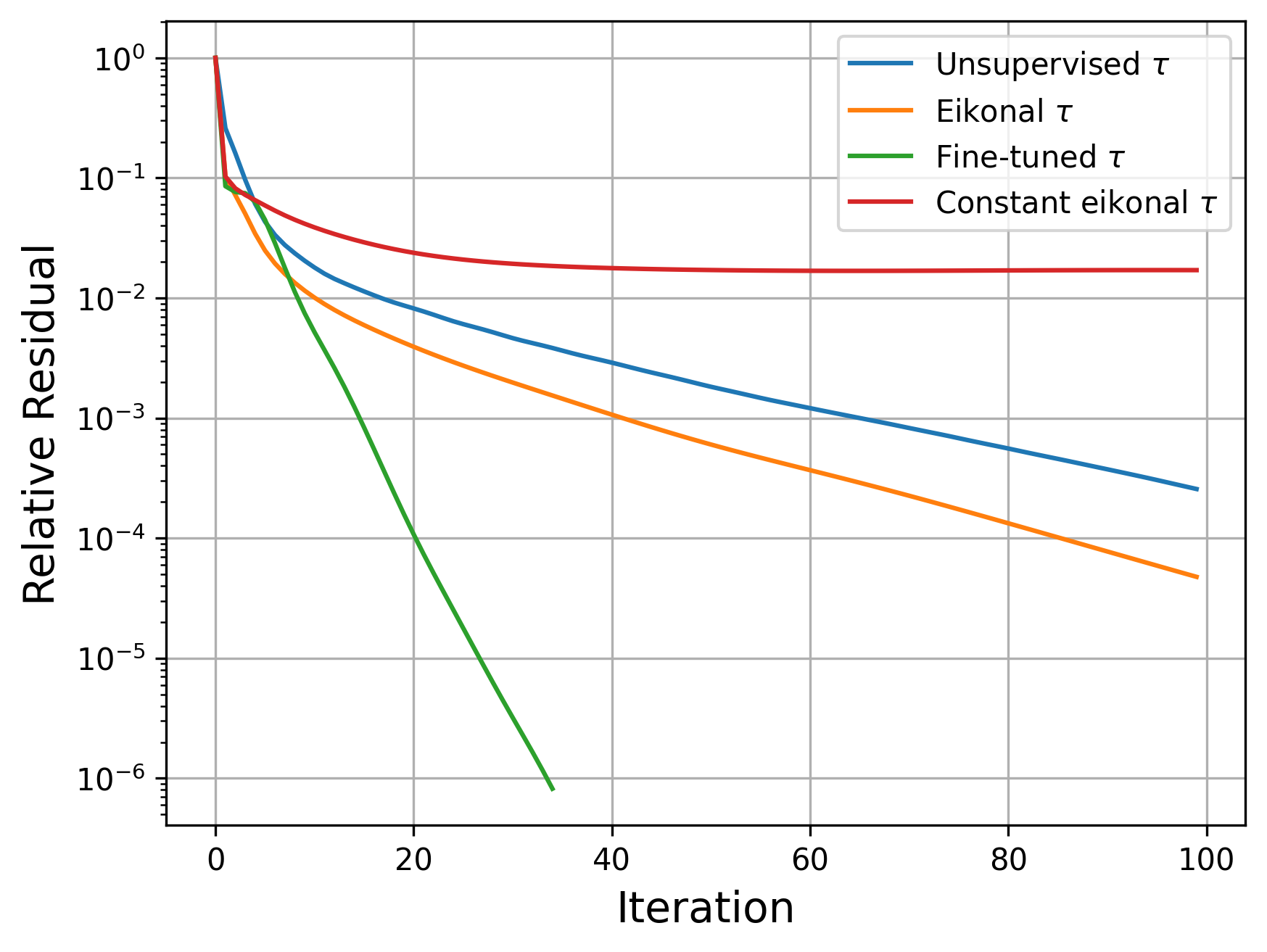}
     \label{fig:fine_tune}
     \caption{Convergence history of Wave-ADR-NS for different $\tau$, where the constant eikonal $\tau$ is defined in \cref{eq:eikonal_cons}.}
\end{figure}

Now, we analyze the impact of using different $\tau$ on the performance of  Wave-ADR-NS. 
\cref{fig:fine_tune} illustrates the convergence history of Wave-ADR-NS with different $\tau$ when solving the Helmholtz equation with $N=128$ and $\omega=20\pi$, utilizing the slowness model depicted in Figure \ref{fig:kappa_test}. 
It can be observed that while the solver with $\tau$ from unsupervised learning achieves convergence, its convergence rate is comparatively slower than that obtained from solving the eikonal equation. 
This drives us to develop a new method to select a better  $\tau$.

We first train the FNO to learn $\tau$ obtained by solving the eikonal equation. We utilize factored fast marching method\,\cite{treister2016fast} to compute $\tau_1$ as supervised data.
\Cref{tab:noresults} illustrates the relative error of FNO tested on slowness models of various resolutions and compares it with UNet \cite{azulay2022multigrid}.
It can be observed that FNO exhibits robust generalization capabilities for larger-scale test cases compared with UNet.
Then, we transfer the parameters of the learned FNO to Wave-ADR-NS and fine-tune them (with 10 training epochs) together with the parameters $\alpha$ using the loss function \cref{eq:loss} to obtain the final Wave-ADR-NS in such a semi-supervised manner. 
\Cref{fig:fine_tune} illustrates the convergence history of Wave-ADR-NS with the fine-tuned $\tau$. It is observed that the fine-tuned $\tau$ exhibits much faster convergence than others.
\begin{table}[!htb]
     \centering
     \footnotesize{
     \caption{Average relative test error of FNO and UNet across different slowness models.}
     \label{tab:noresults}
     \begin{tabular}{ccccccc}
     \toprule
     $N$    & 128    & 256    & 512    & 1024   & 2048   & 4096   \\ \midrule
     FNO  & \textbf{0.2401} & \textbf{0.1972} & \textbf{0.2370} & \textbf{0.2907} & \textbf{0.3379} & \textbf{0.3285} \\ \midrule
     UNet & 1.2118 & 0.7144 & 0.3322 & 0.4960 & 0.6219 & 1.1132 \\ \bottomrule
     \end{tabular}
}
\end{table}

\begin{figure}[!htbp]
     \centering
     \subfigure[$K=1$]{\includegraphics[width=0.23\textwidth]{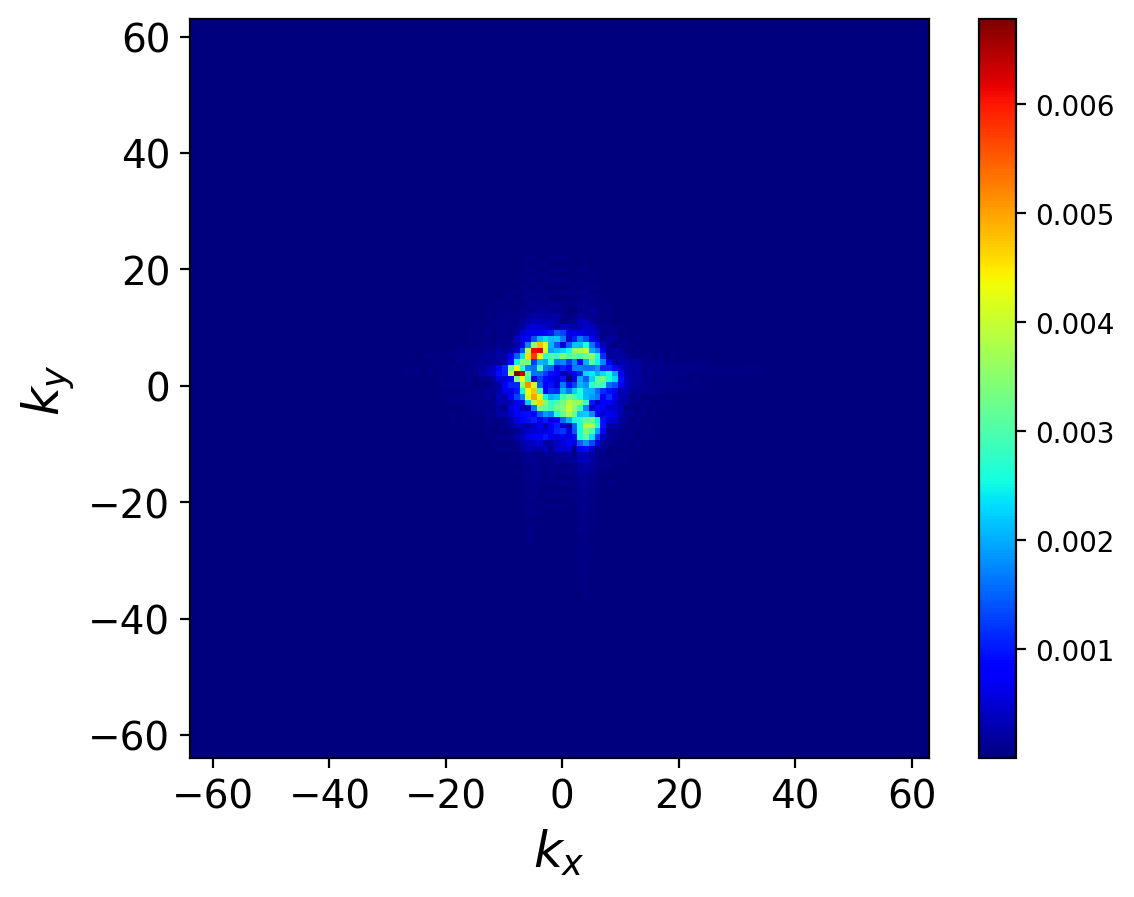}}
     \subfigure[$K=2$]{\includegraphics[width=0.23\textwidth]{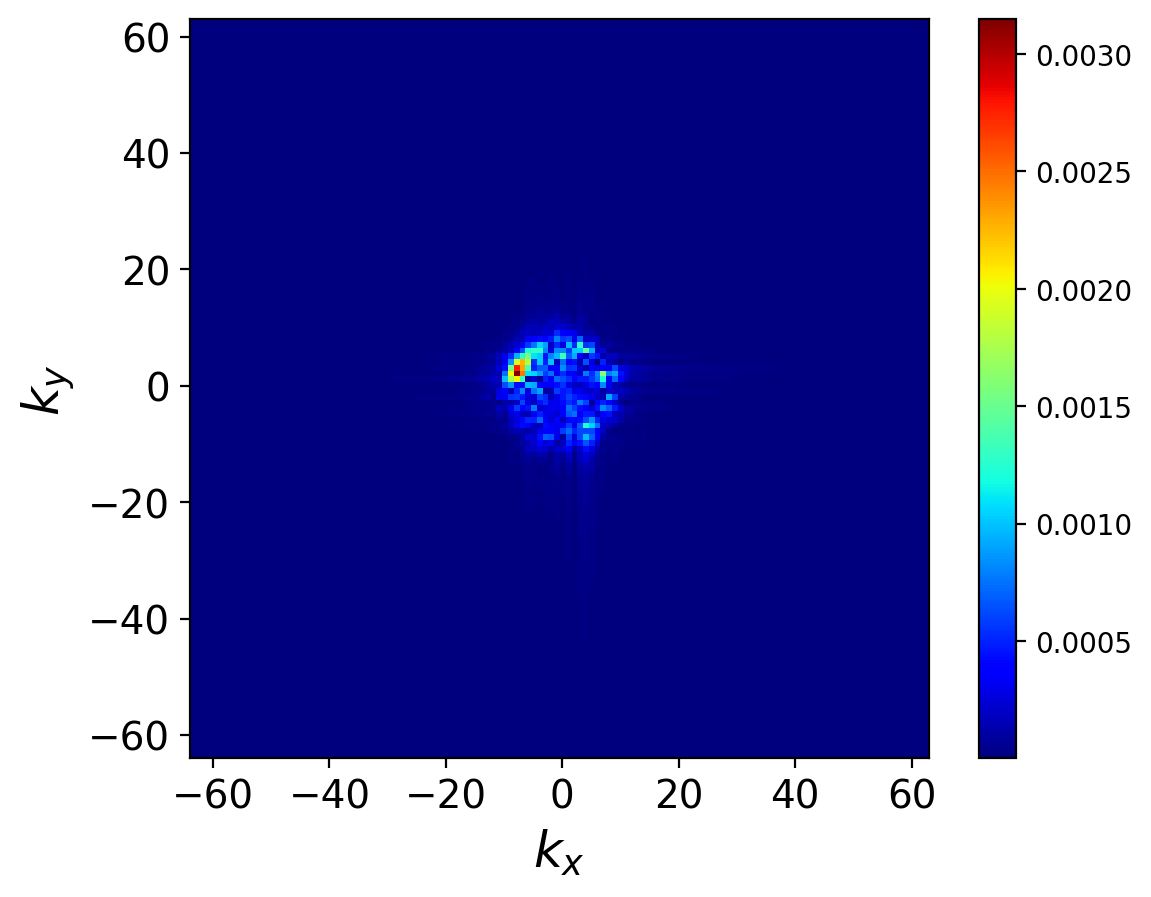}}
     \subfigure[$K=3$]{\includegraphics[width=0.23\textwidth]{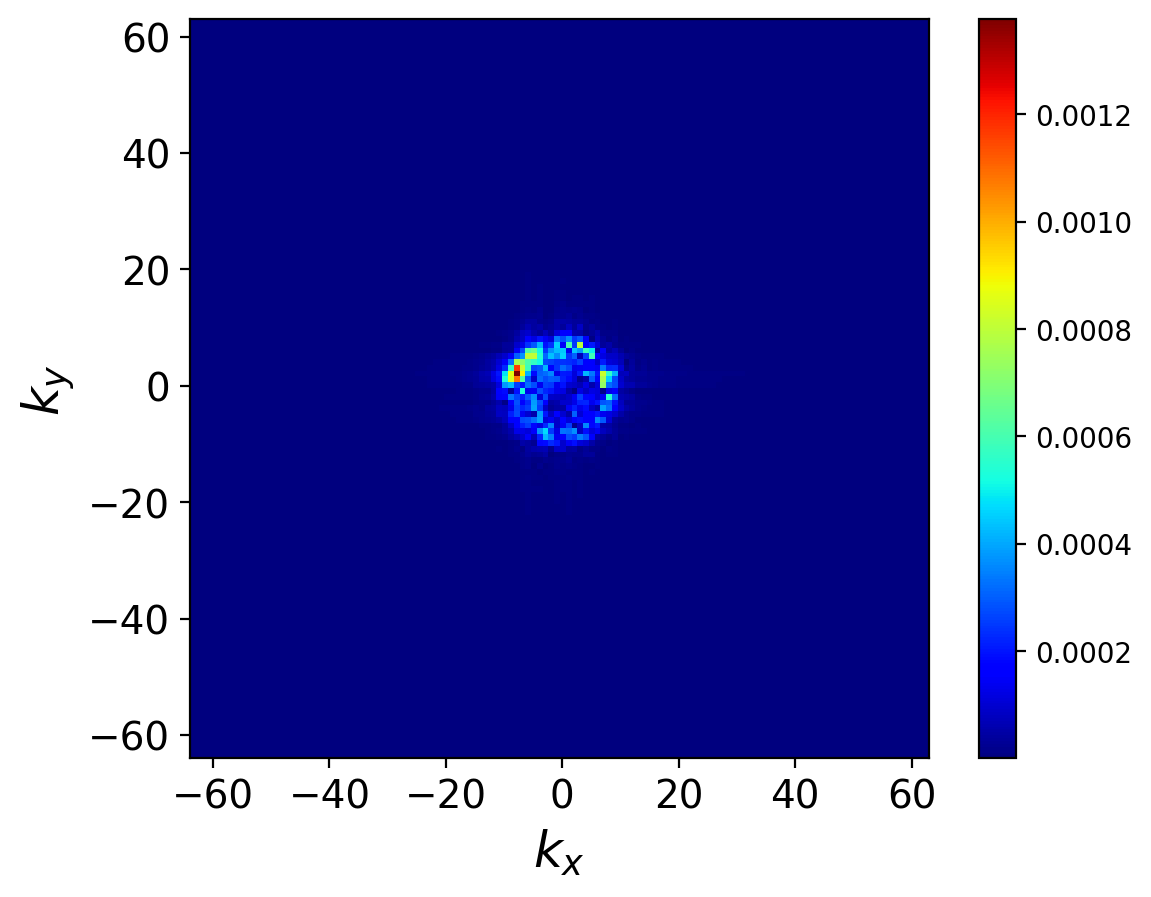}}
     \subfigure[$K=4$]{\includegraphics[width=0.23\textwidth]{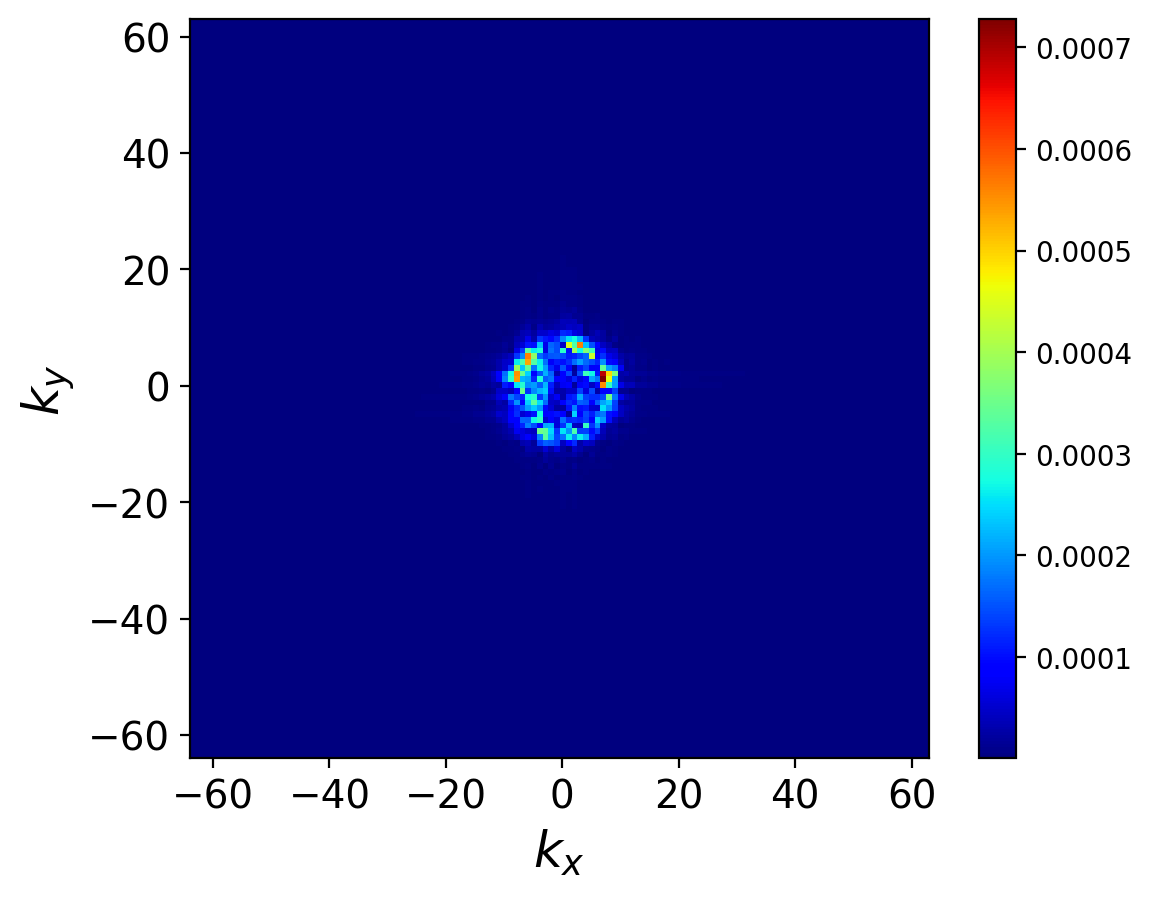}}\\
     \subfigure[$K=1$]{\includegraphics[width=0.23\textwidth]{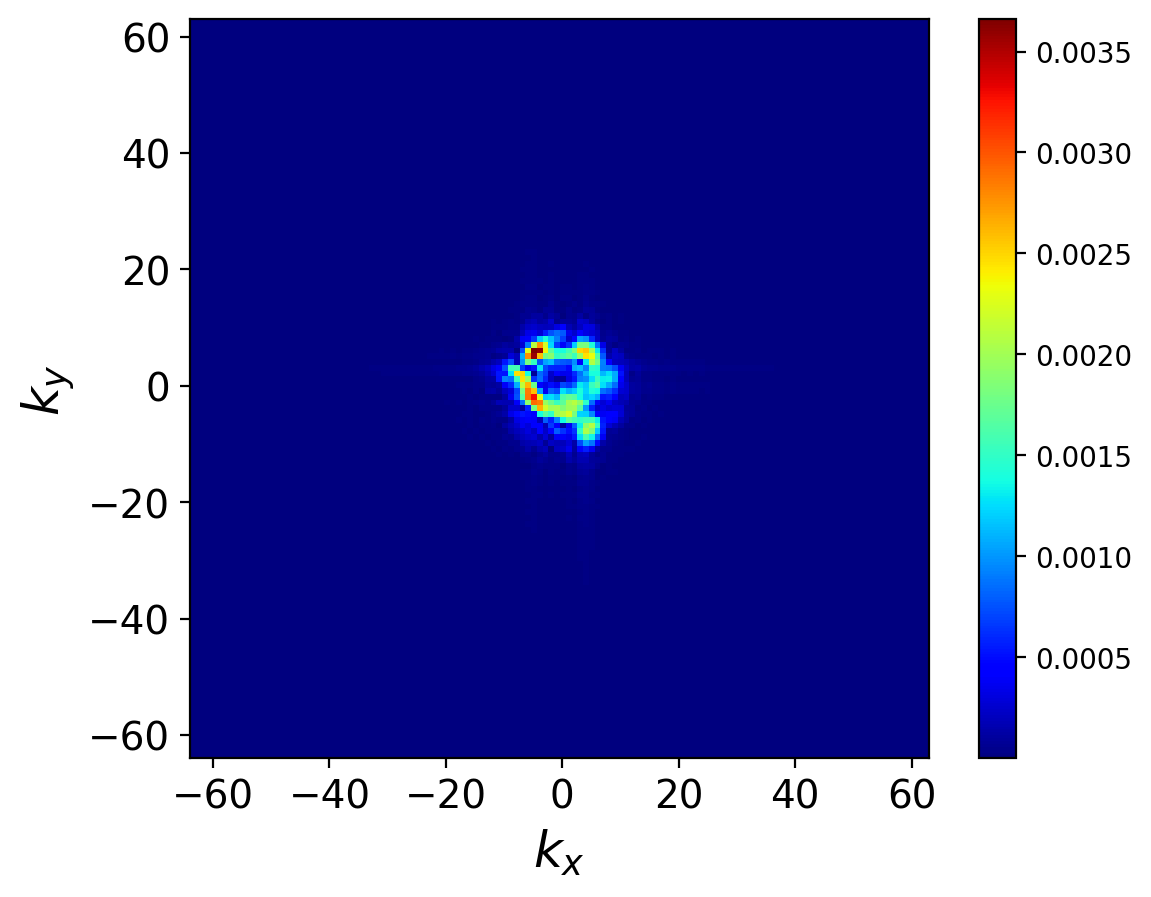}}
     \subfigure[$K=2$]{\includegraphics[width=0.23\textwidth]{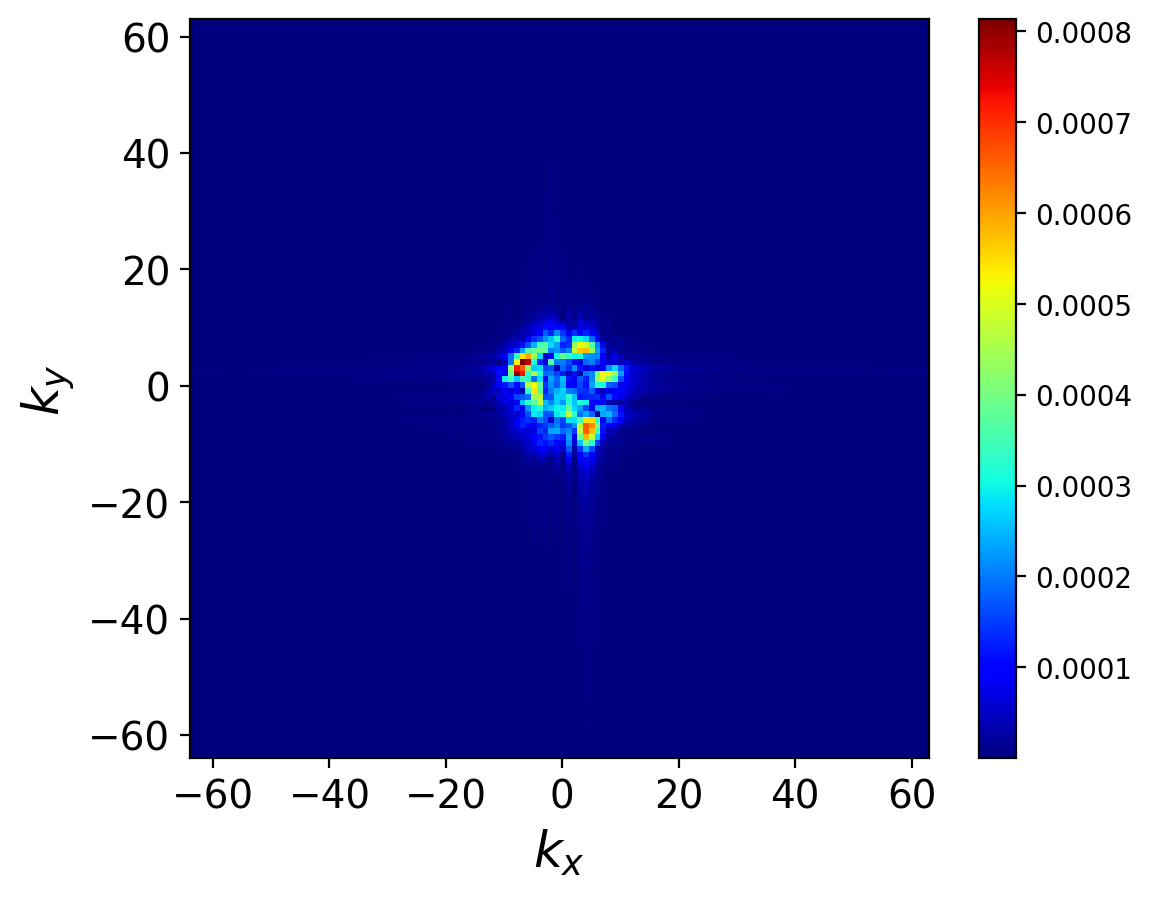}}
     \subfigure[$K=3$]{\includegraphics[width=0.23\textwidth]{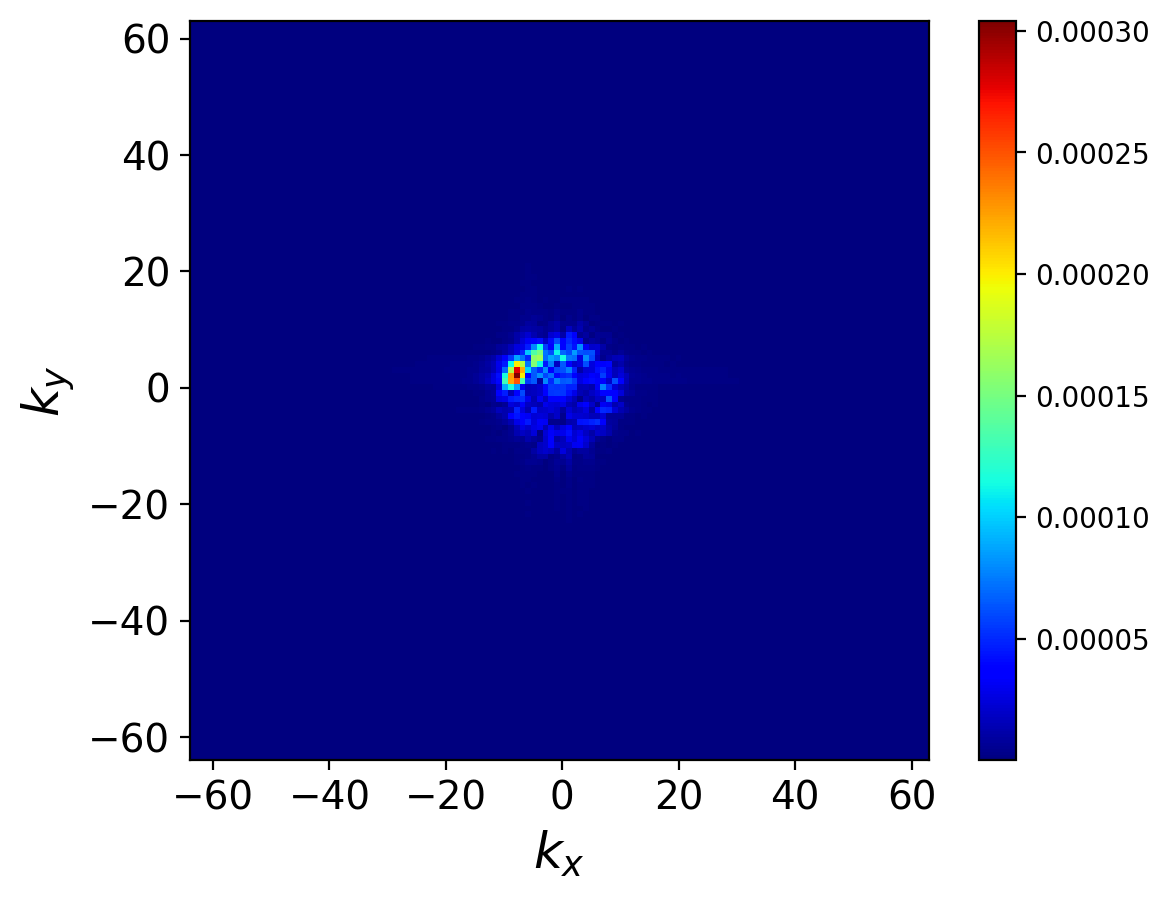}}
     \subfigure[$K=4$]{\includegraphics[width=0.22\textwidth]{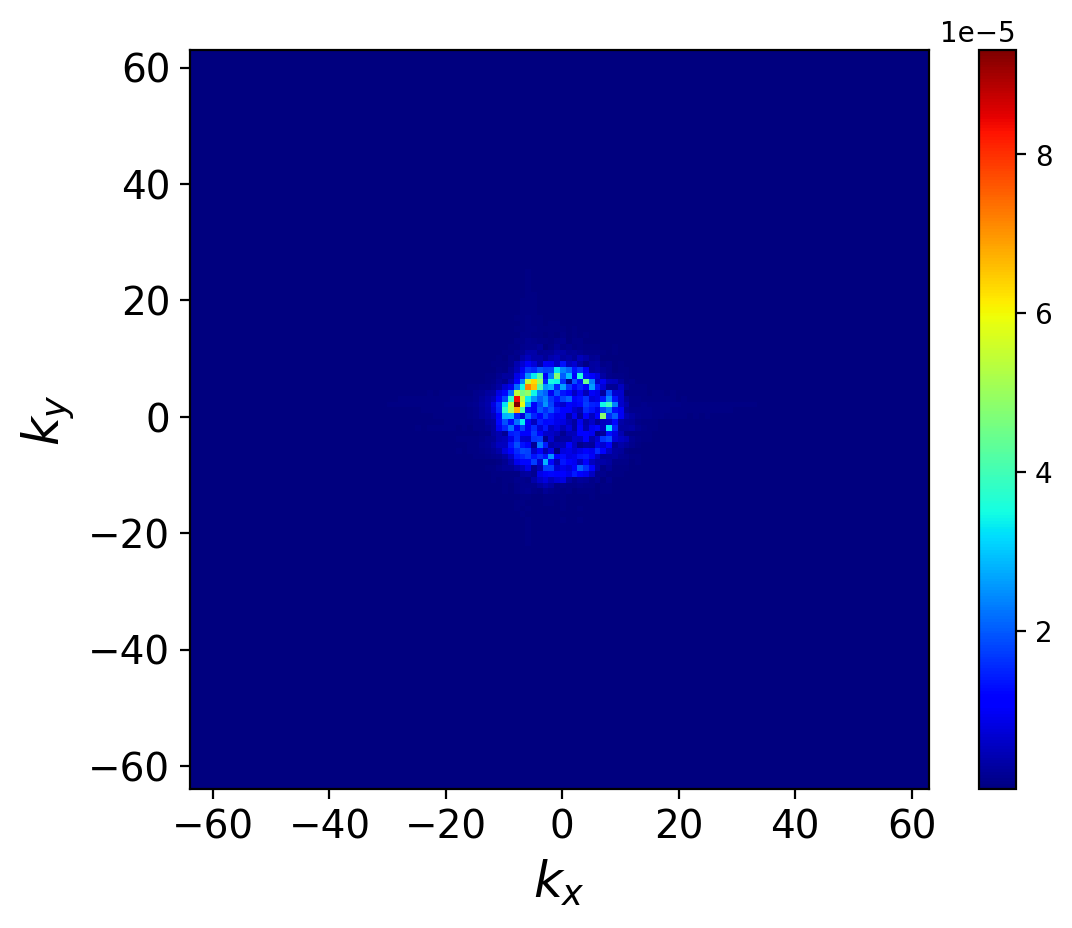}}
     \caption{Comparison of the module of errors in Fourier space of Wave-ADR-NS at different iteration steps when using $\tau$ obtained by solving the eikonal equation. Top: Iterative error before correction; Bottom: Error correction obtained through the ADR cycle.}
     \label{fig:eik_error}
 \end{figure}

\begin{figure}[!htbp]
     \centering
     \subfigure[$K=1$]{\includegraphics[width=0.23\textwidth]{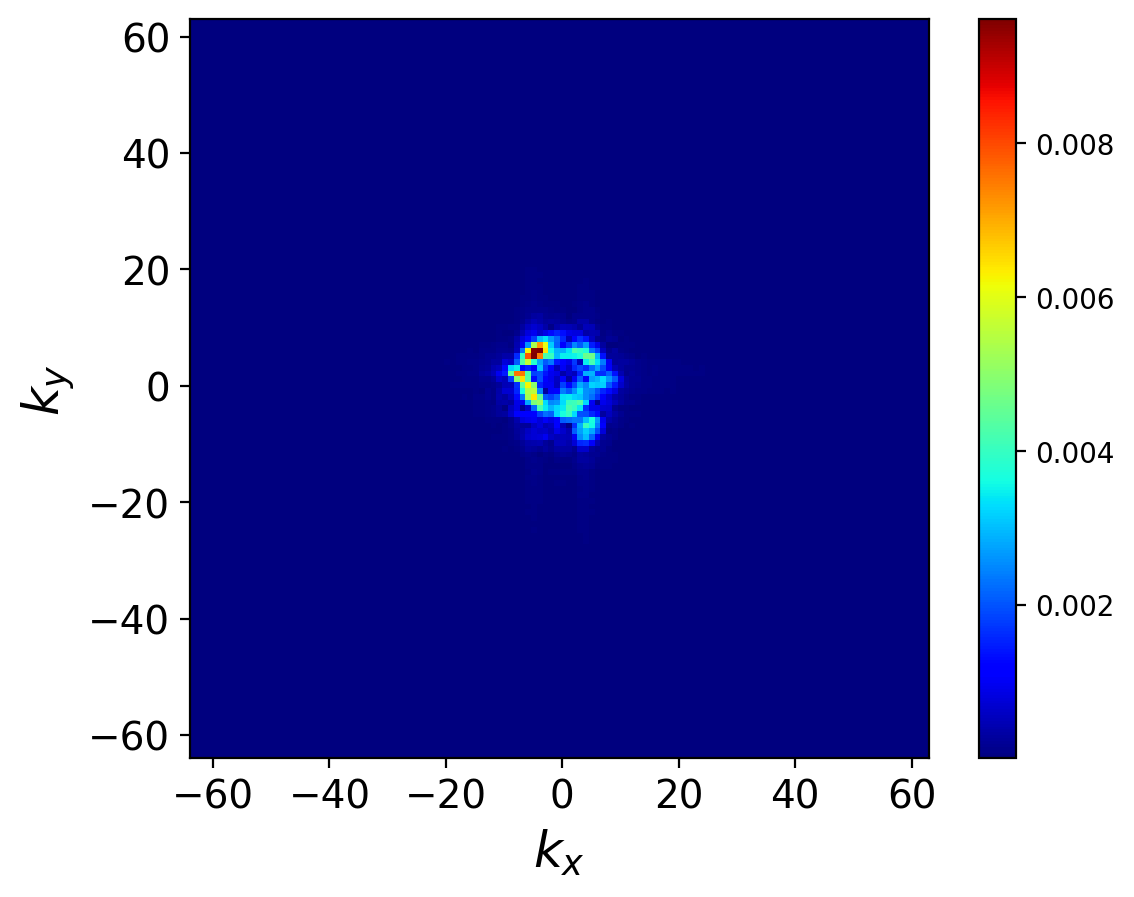}}
     \subfigure[$K=2$]{\includegraphics[width=0.23\textwidth]{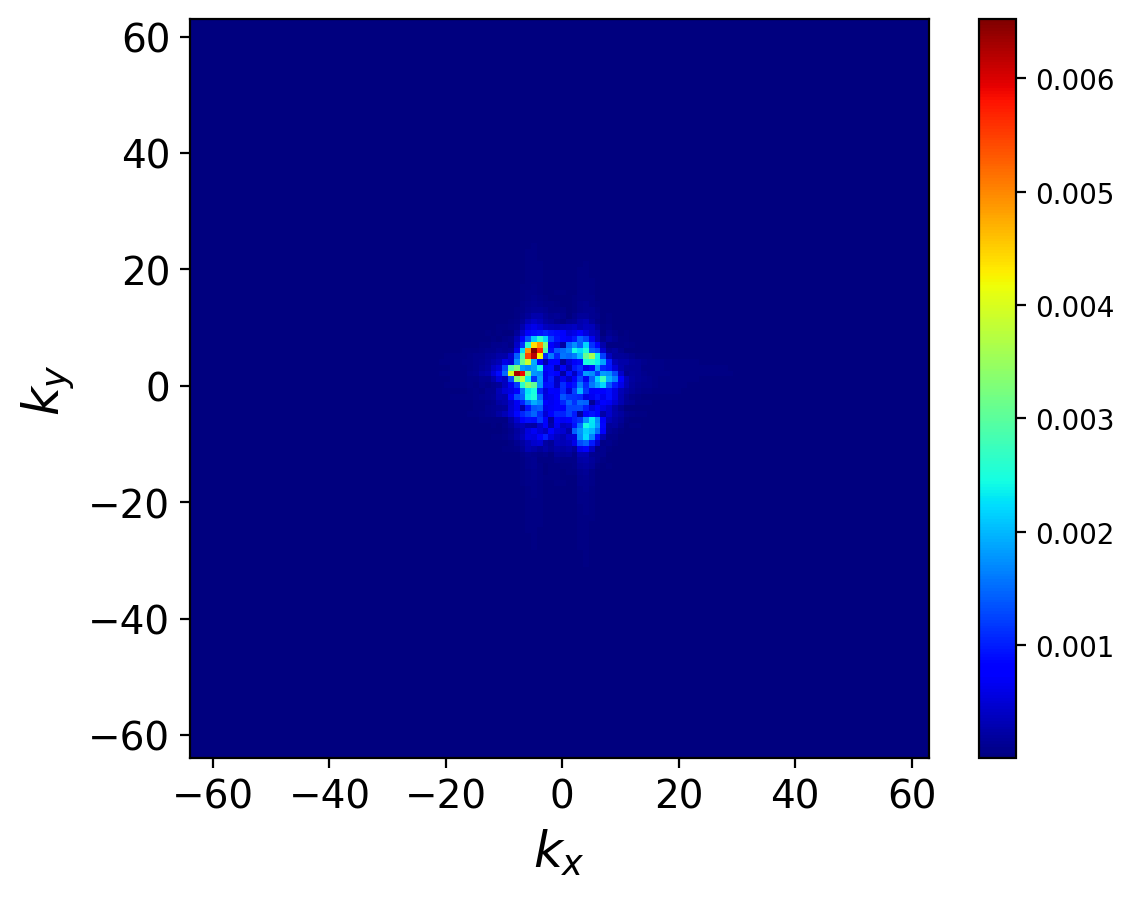}}
     \subfigure[$K=3$]{\includegraphics[width=0.23\textwidth]{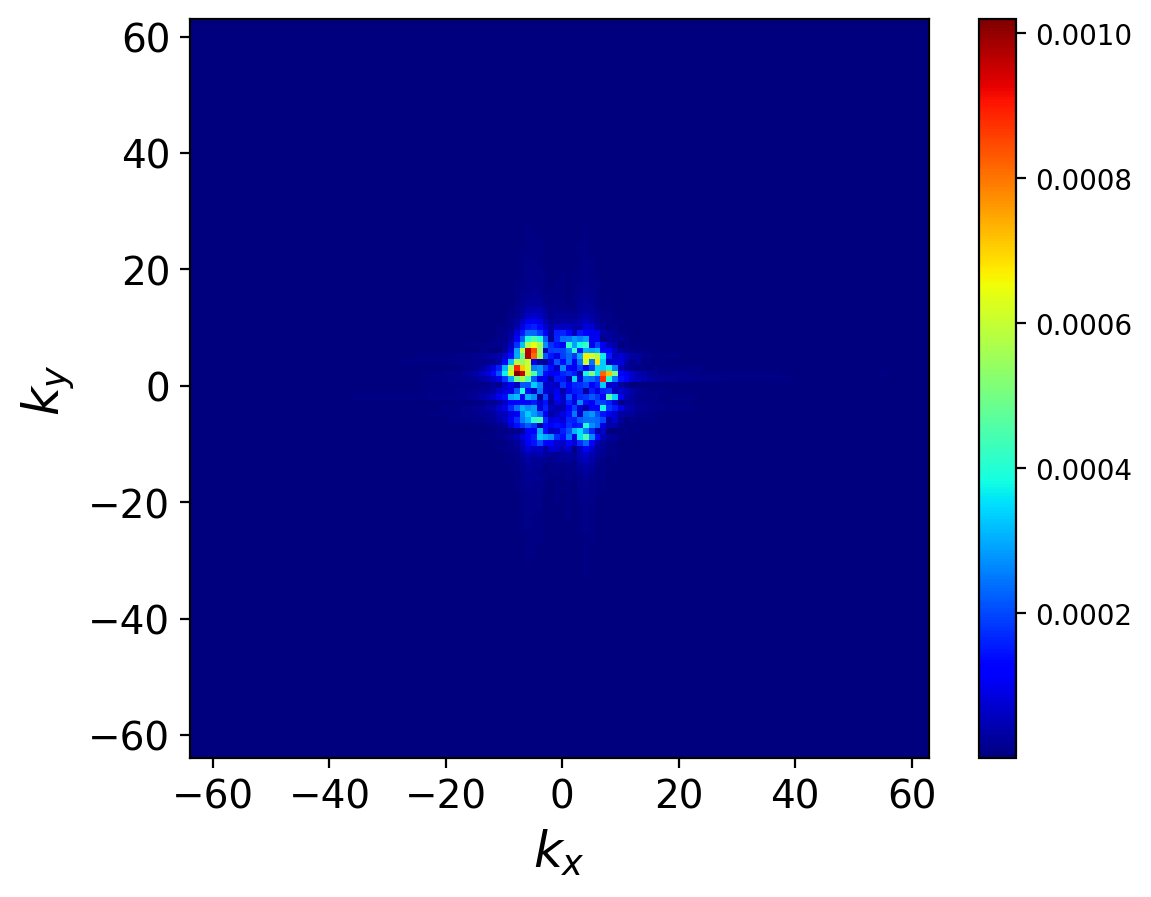}}
     \subfigure[$K=4$]{\includegraphics[width=0.23\textwidth]{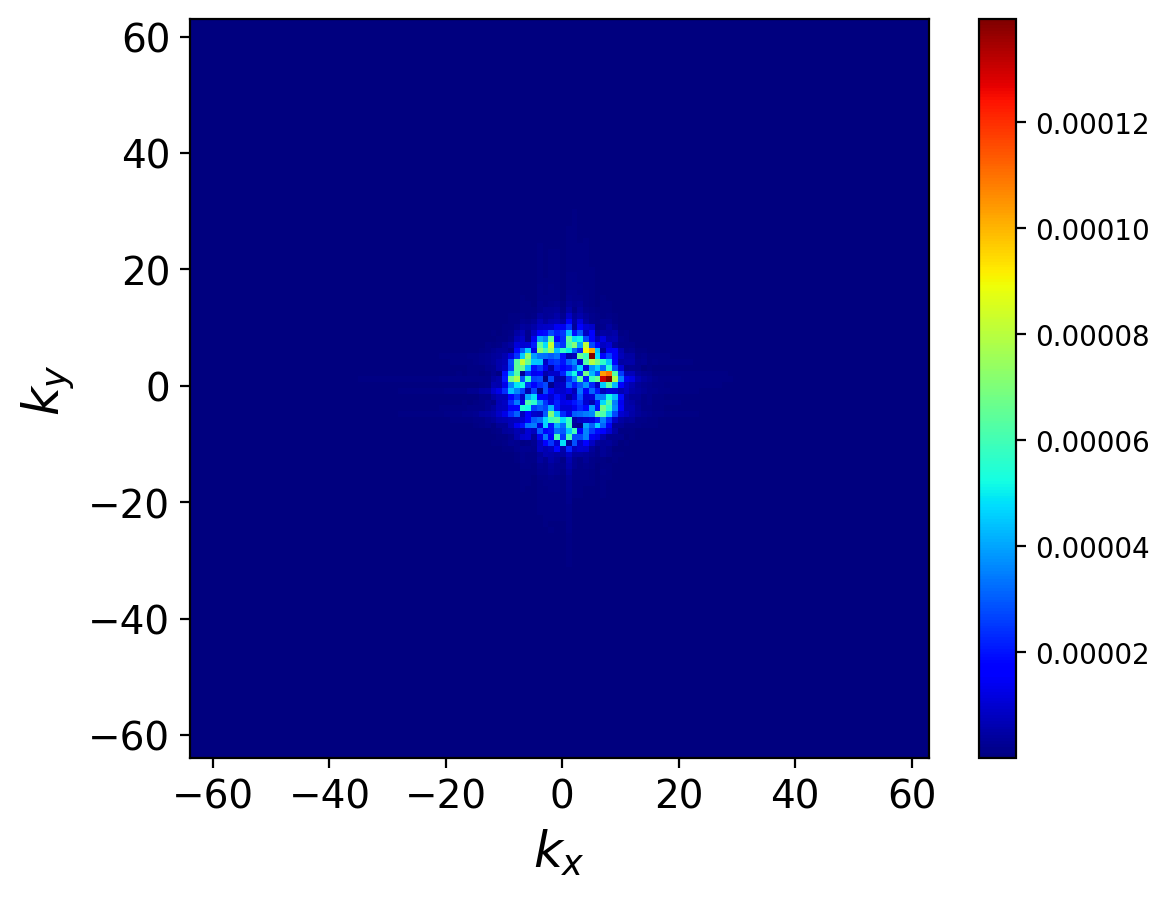}}\\
     \subfigure[$K=1$]{\includegraphics[width=0.23\textwidth]{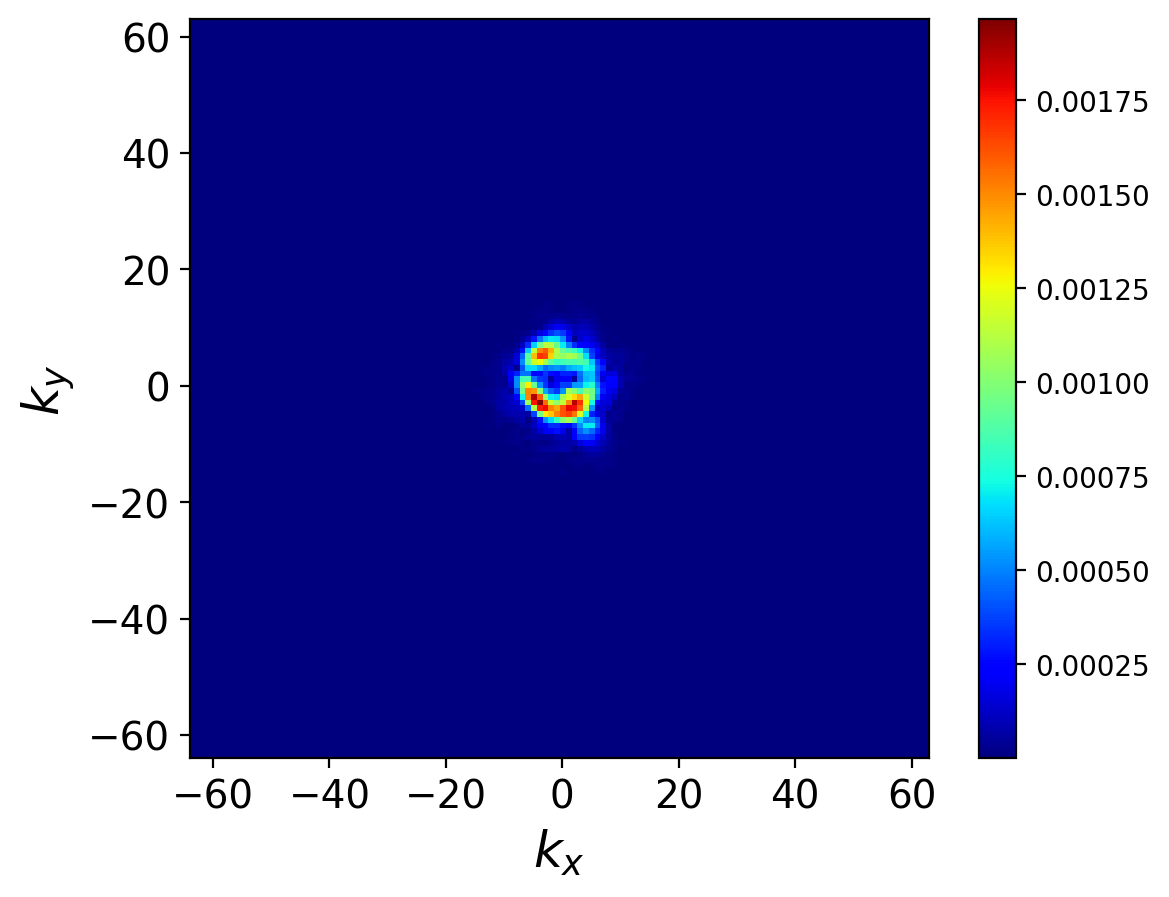}}
     \subfigure[$K=2$]{\includegraphics[width=0.23\textwidth]{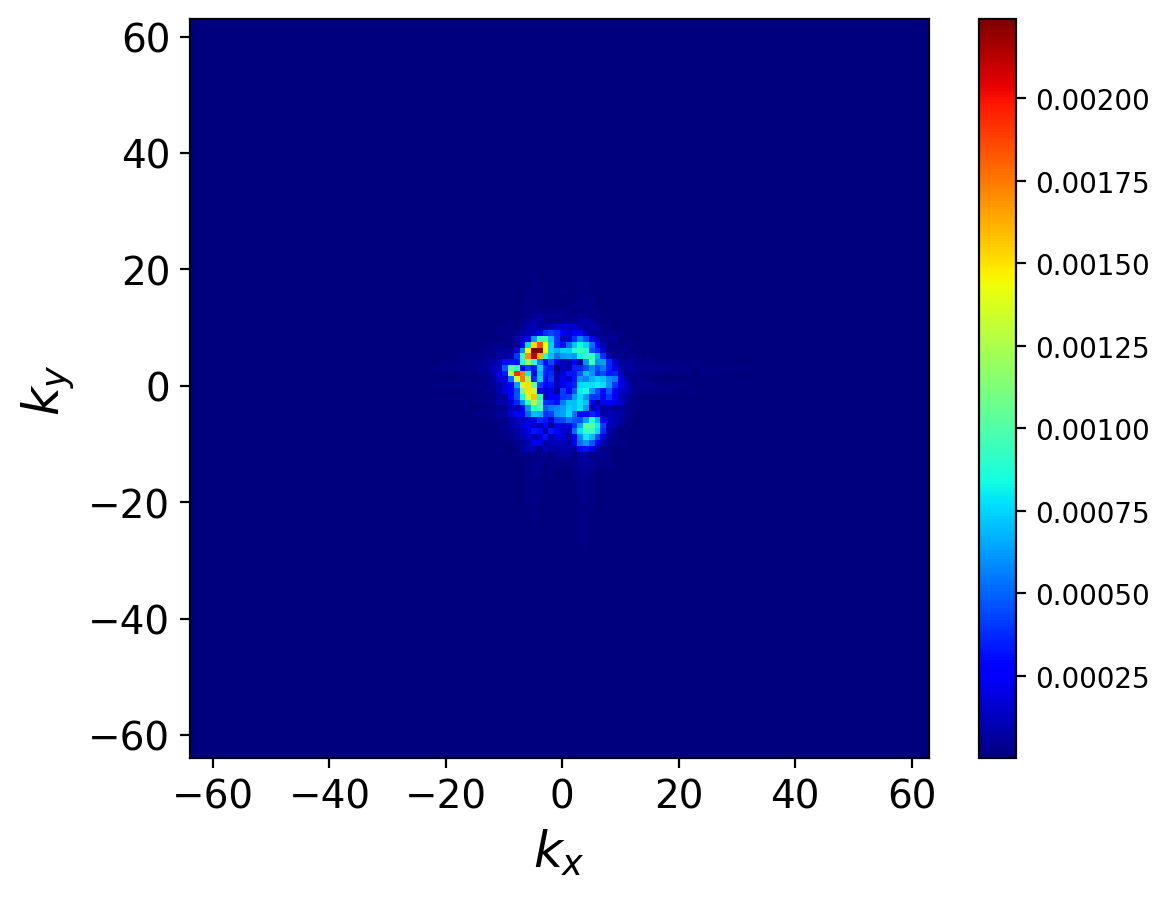}}
     \subfigure[$K=3$]{\includegraphics[width=0.23\textwidth]{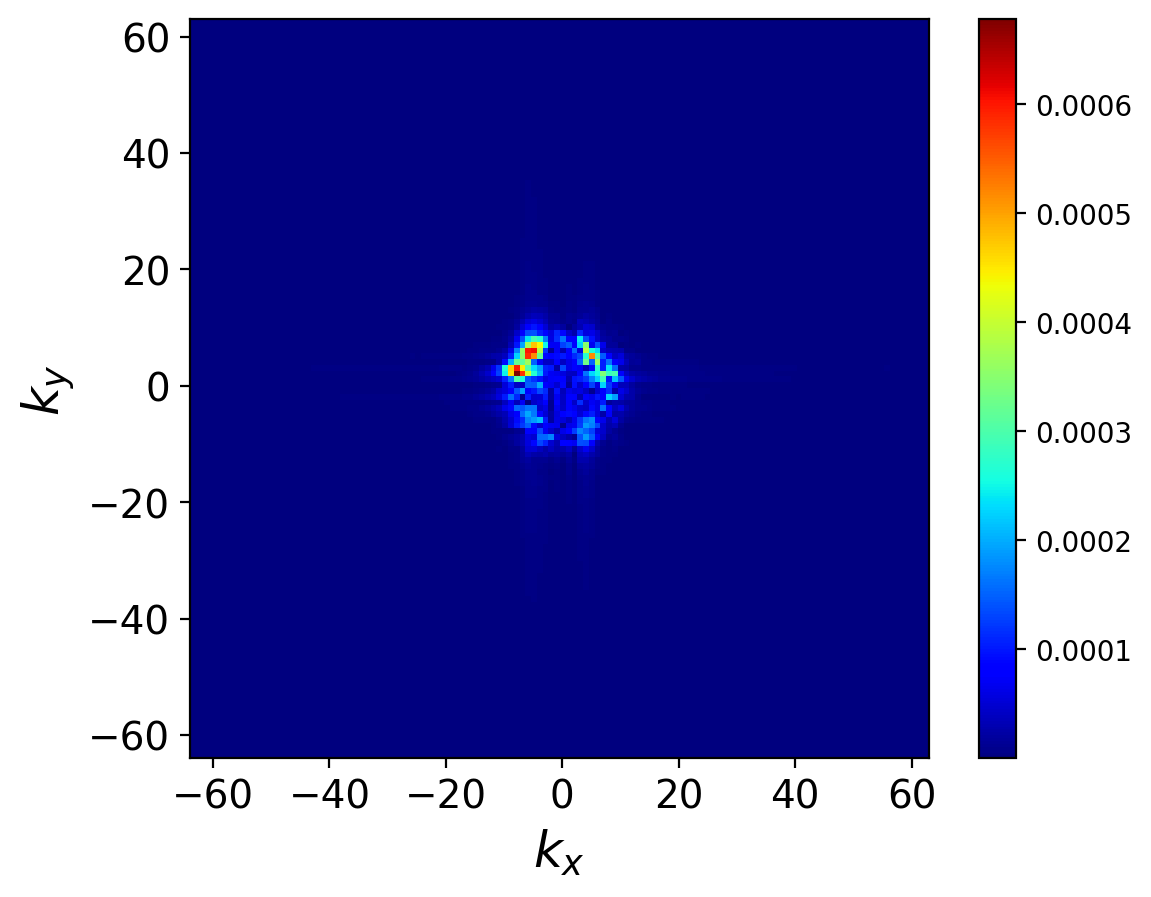}}
     \subfigure[$K=4$]{\includegraphics[width=0.22\textwidth]{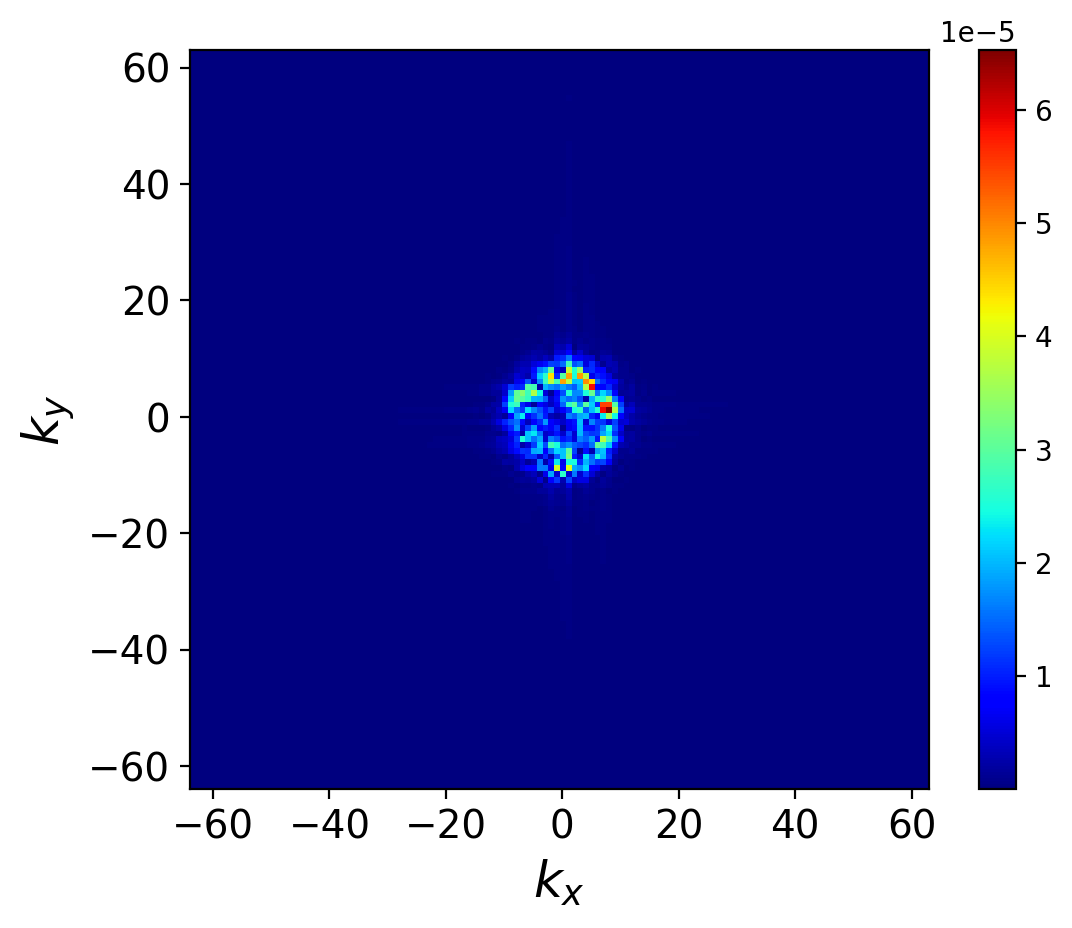}}
     \caption{Comparison of the module of errors in Fourier space of Wave-ADR-NS at different iteration steps when using the fine-tuned $\tau$. Top: Iterative error before correction; Bottom: Error correction obtained through the ADR cycle.}
     \label{fig:fine_tune_error}
 \end{figure}

The reason of the faster convergence of the fine-tuned $\tau$ can be demonstrated through comparing the behavior of Wave-ADR-NS with different $\tau$ values at various iterative steps $K$.
\Cref{fig:eik_error} and \Cref{fig:fine_tune_error} present the results for $\tau$ obtained by solving the eikonal equation and the fine-tuned $\tau$, respectively. For both figures, we compare the quantities representing the iterative errors before and after correction obtained by the ADR cycle in Fourier space.
As demonstrated in the first iteration step, the $\tau$ obtained by solving the eikonal equation can effectively captures the distribution of the characteristic error. However, applying $\tau$ in the iterative method requires additional considerations. Specifically, after the first Wave-ADR cycle, the distribution of the iterative error may deviate from the original characteristic error. In such cases, continuing to use the eikonal-based $\tau$ may no longer be optimal, as evidenced by subsequent iterations where the fine-tuned $\tau$ performs better.
Ideally, $\tau$ should adapt dynamically throughout the iterative process. At present, we regard the fine-tuned $\tau$, trained using the relative residual as the loss function, as an approximation of the best fixed $\tau$ that effectively approximates $A^{-1}$.

\begin{figure}[!htb]
     \centering
     \includegraphics[width=0.4\textwidth]{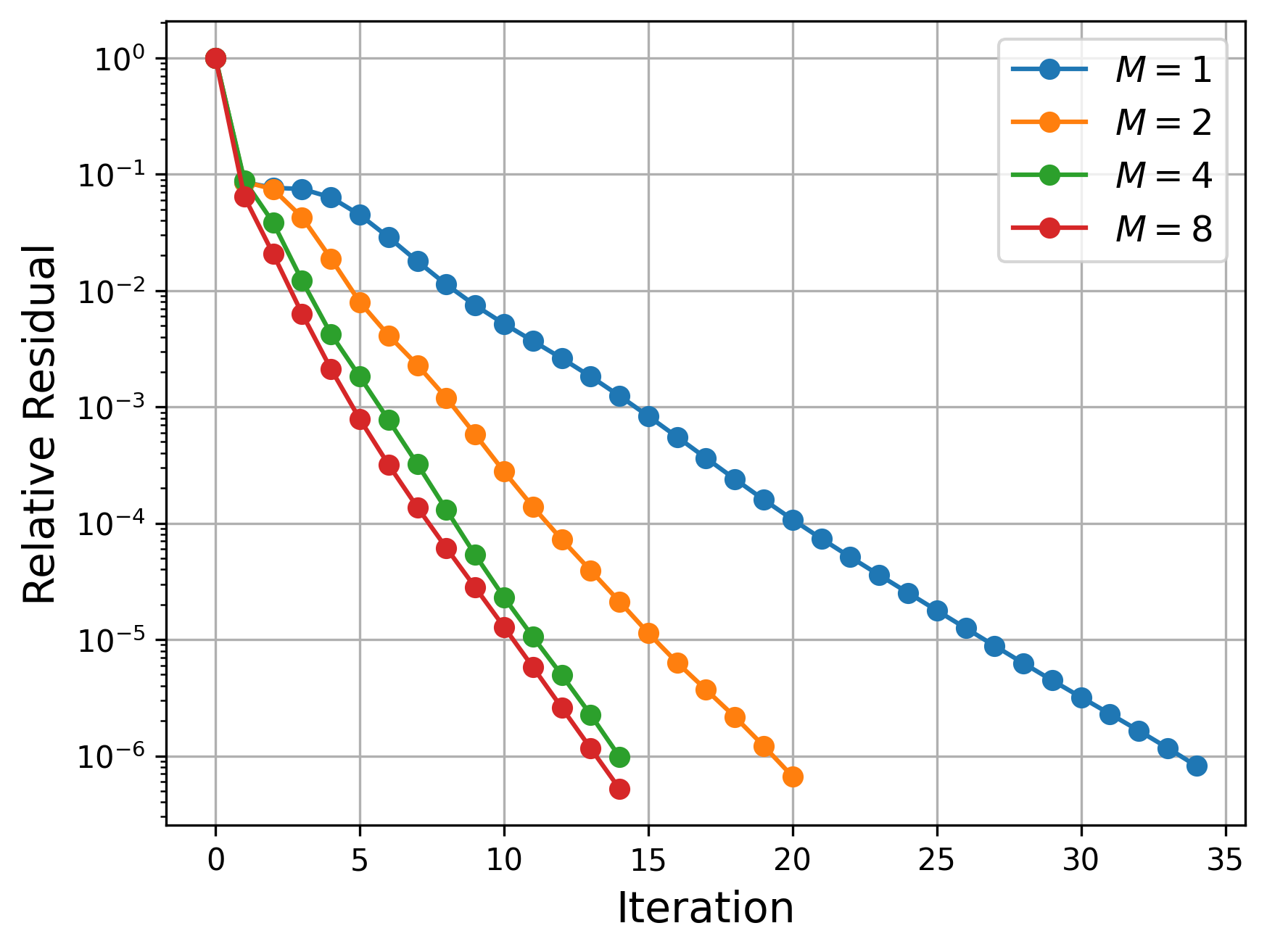}
     \caption{Convergence history of Wave-ADR-NS for different ADR correction steps.}
     \label{fig:diffM}
\end{figure}

We further accelerate Wave-ADR-NS by considering multiple ADR correction steps. 
Regarding the ADR correction step as a post-smoothing relaxation at level $kh \approx 1$, and employ multiple correction steps instead of the previously employed single correction step.
\cref{fig:diffM} illustrates the convergence of Wave-ADR-NS with different ADR correction steps.
It is observed that as $M$ increases, the convergence rate of Wave-ADR-NS significantly improves. 
Therefore, in the following experiments, we use multiple correction steps.

\subsection{Performance of Wave-ADR-NS}\label{sec:053}

\subsubsection{In-distribution generalization}
We first evaluate the generalization capabilities of Wave-ADR-NS on slowness models that share the same distribution as the training dataset, which is converted  from the CIFAR-10 dataset. The right-hand side is a point source located at the center, and the angular frequency $\omega$ is tested up to $640\pi\approx 2010$, with the corresponding discretization scale is 4096. For each scale, we randomly sample 10 different slowness models.
All methods are used as preconditioners for FGMRES with a restart set to 20 for a fair comparison, and the stopping criterion is $\|\mathbf{r}^{(k)}\|/\|\mathbf{g}\| < 10^{-6}$, starting from a zero initial guess. The number of ray functions used in the Wave-Ray method is set to $M=8$, and we apply the same number of ADR correction steps. The CSL preconditioner is configured with $\beta = 0.5\omega^2$ and is inverted by a single MG V-cycle.
\Cref{tab:cifar_iter} and \Cref{tab:cifar_time} present the average iteration counts and computation times required. It can be observed that the Wave-ADR-NS achieves this accuracy with fewer iterations and lower computational time, demonstrating superior efficiency and robustness than classical methods such as Wave-Ray or CSL, as well as the deep learning-enhanced MG preconditioner Encoder-Solver \cite{azulay2022multigrid}.
\begin{table}[!htb]
     \centering
     \footnotesize{
     \caption{Average iteration counts required for different preconditioners when solving Helmholtz with slowness models from CIFAR-10.}
     \label{tab:cifar_iter}
     \begin{tabular}{@{}lllllll@{}}
     \toprule
     $N $              & 128   & 256   & 512    & 1024   & 2048               & 4096               \\ 
     $\omega/2\pi$       & 10    & 20    & 40     & 80     & 160                & 320                \\\midrule
     Wave-ADR-NS     & 8.8   & 15.5  & 28.8   & 65.1               & 137.2              & 344.3              \\
     Wave-Ray        & 23.8  & 46.9  & 93.7   & 191.1              & 404.0              & 874.5              \\
     CSL             & 158.5 & 434.3 & 1163.4 & \textgreater{}2000 & \textgreater{}2000 & \textgreater{}2000 \\
     Encoder-Solver  & 23.5  & 38.3  & 70.9   & 176.5              & \textgreater{}2000 & \textgreater{}2000 \\
     \bottomrule
     \end{tabular}}
\end{table}
\begin{table}[!htb]
     \centering
     \footnotesize{
     \caption{Average computational time (s) required for different preconditioners when solving Helmholtz with slowness models from CIFAR-10.}
     \label{tab:cifar_time}
     \begin{tabular}{@{}lllllll@{}}
     \toprule
     $N $              & 128   & 256   & 512    & 1024   & 2048               & 4096               \\ 
     $\omega/2\pi$       & 10    & 20    & 40     & 80     & 160                & 320                \\\midrule
     Wave-ADR-NS     & 2.31  & 5.60  & 13.72  & 36.96  & 97.42              & 311.27             \\
     Wave-Ray        & 8.35  & 21.24 & 53.61  & 127.67 & 364.91             & \textgreater{}1000 \\
     CSL             & 13.52 & 37.65 & 116.20 & 416.91 & \textgreater{}1000 & \textgreater{}1000 \\
     Encoder-Solver  & 2.03  & 3.43  & 7.27   & 37.56  & \textgreater{}1000 & \textgreater{}1000 \\
     \bottomrule
     \end{tabular}}
\end{table}

\subsubsection{Out-of-distribution transfer}
\begin{figure}[!htb]
     \centering
     \includegraphics[width=0.8\textwidth]{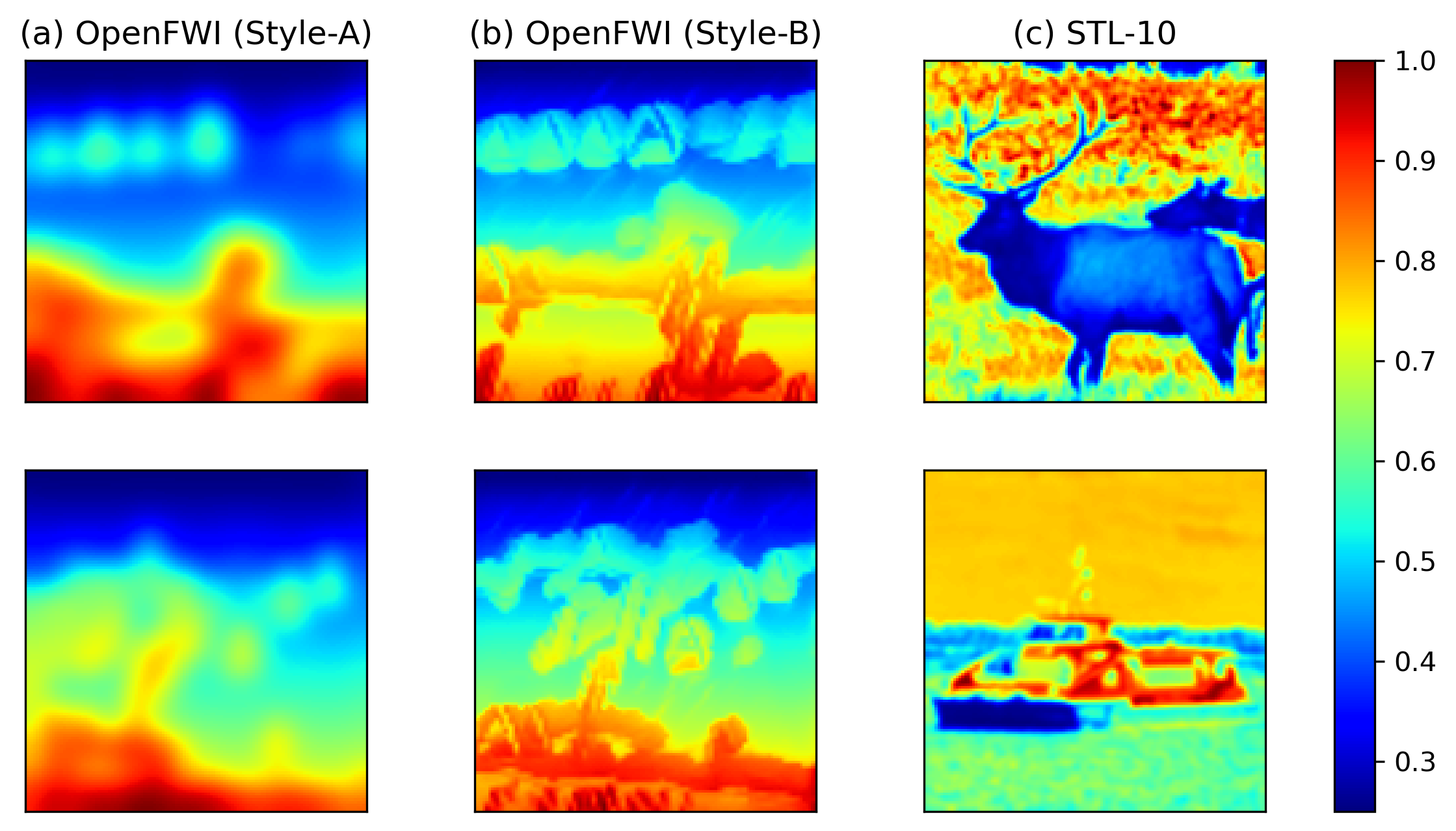}
     \caption{Example slowness models used for out-of-distribution testing. Models from the (a)  OpenFWI (Style-A) dataset; (b) OpenFWI (Style-B) dataset; (c) STL-10 dataset.}
     \label{fig:ood}
\end{figure}

We then evaluate the transfer capabilities of Wave-ADR-NS on slowness models and right-hand sides with distributions different from the training dataset. First, we test the performance of Wave-ADR-NS on slowness models from OpenFWI (Style-A), OpenFWI (Style-B) \cite{deng2022openfwi}, and STL-10 \cite{coates2011analysis}. 
OpenFWI (Style-A) exhibits smoothness comparable to CIFAR-10, while OpenFWI (Style-B) and STL-10 present rougher textures. We scale $s$ values to the range $[0.25, 1]$. \cref{fig:ood} shows the example slowness models.

We use the Wave-ADR-NS trained on the CIFAR-10 dataset to solve Helmholtz equations with a centered point source but  slowness models from the three datasets \textit{without retraining}.
On each scale, we randomly sample 10 different slowness models and compute the average iterations required for the relative residual to drop below $10^{-6}$ as done before. 

\cref{tab:other_square} shows the results. As can be seen, for the slowness models from OpenFWI (Style-A), the performance of Wave-ADR-NS is comparable to that on CIFAR-10. On the more challenging datasets, OpenFWI (Style-B) and STL-10, the number of iterations increases slightly but remains stable, demonstrating strong transfer capabilities.
\begin{table}[!htb]
     \centering
     \footnotesize{
     \caption{Average iteration counts of Wave-ADR-NS when solving Helmholtz equations with slowness models from out-of-distribution datasets.}
     \label{tab:other_square}
     \begin{tabular}{@{}lllllll@{}}
     \toprule
     $N $              & 128   & 256   & 512    & 1024   & 2048               & 4096               \\ 
     $\omega/2\pi$       & 10    & 20    & 40     & 80     & 160                & 320                \\\midrule
     OpenFWI (Style-A)     & 9.3   & 15.3  & 29.9   & 62.3  & 133.7              & 351.9              \\
     OpenFWI (Style-B)        & 9.1  & 15.0  & 30.8   & 70.7  & 144.2           &    407.1          \\
     STL-10            & 11.5 & 27.6 & 56.8 & 109.5 & 224.2 & 451.2 \\ \bottomrule
     \end{tabular}}
\end{table}

\begin{figure}[!htbp]
     \centering
     \subfigure[CIFAR-10]{\includegraphics[width=0.23\textwidth]{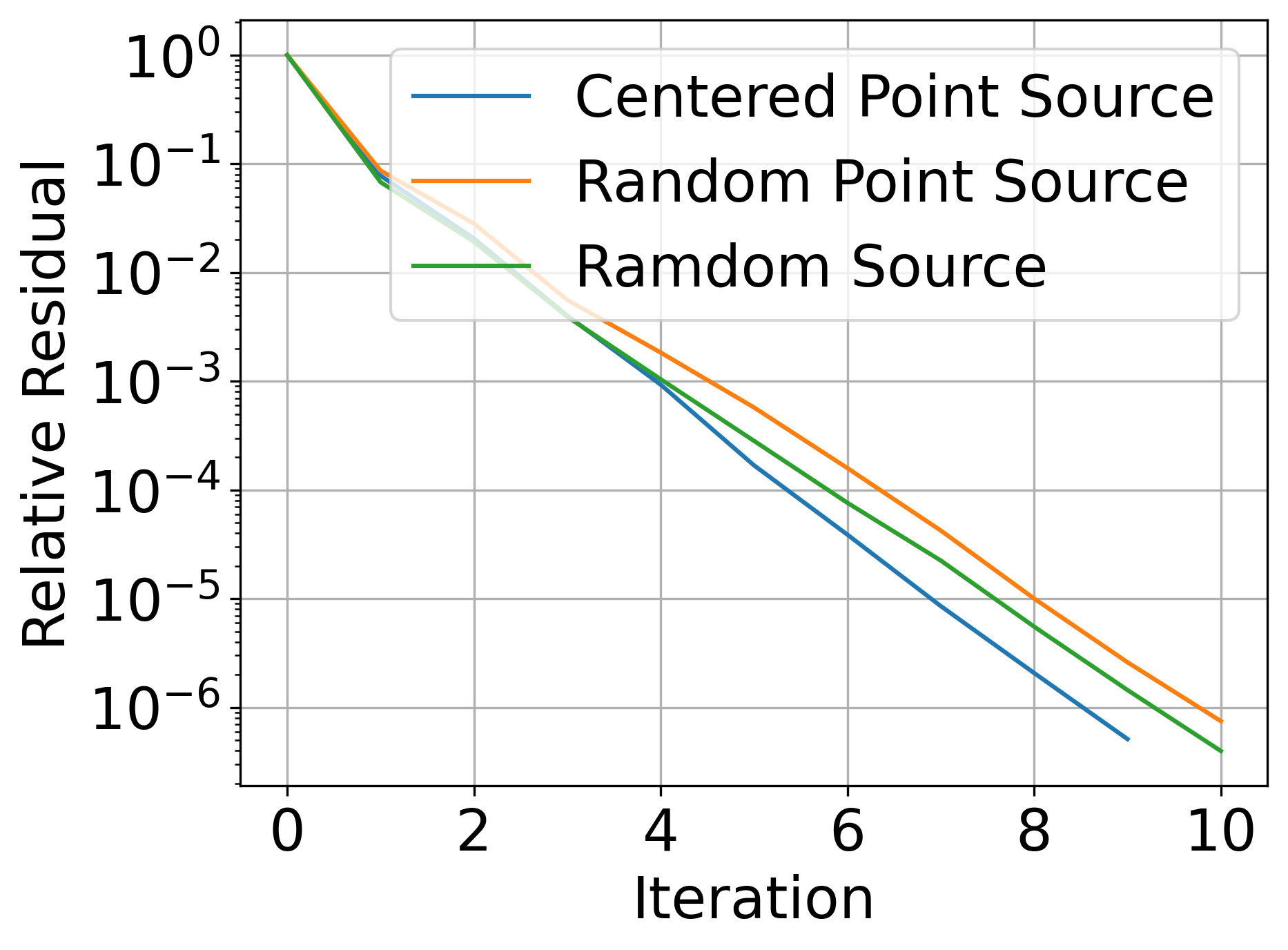}}
     \subfigure[Style-A]{\includegraphics[width=0.23\textwidth]{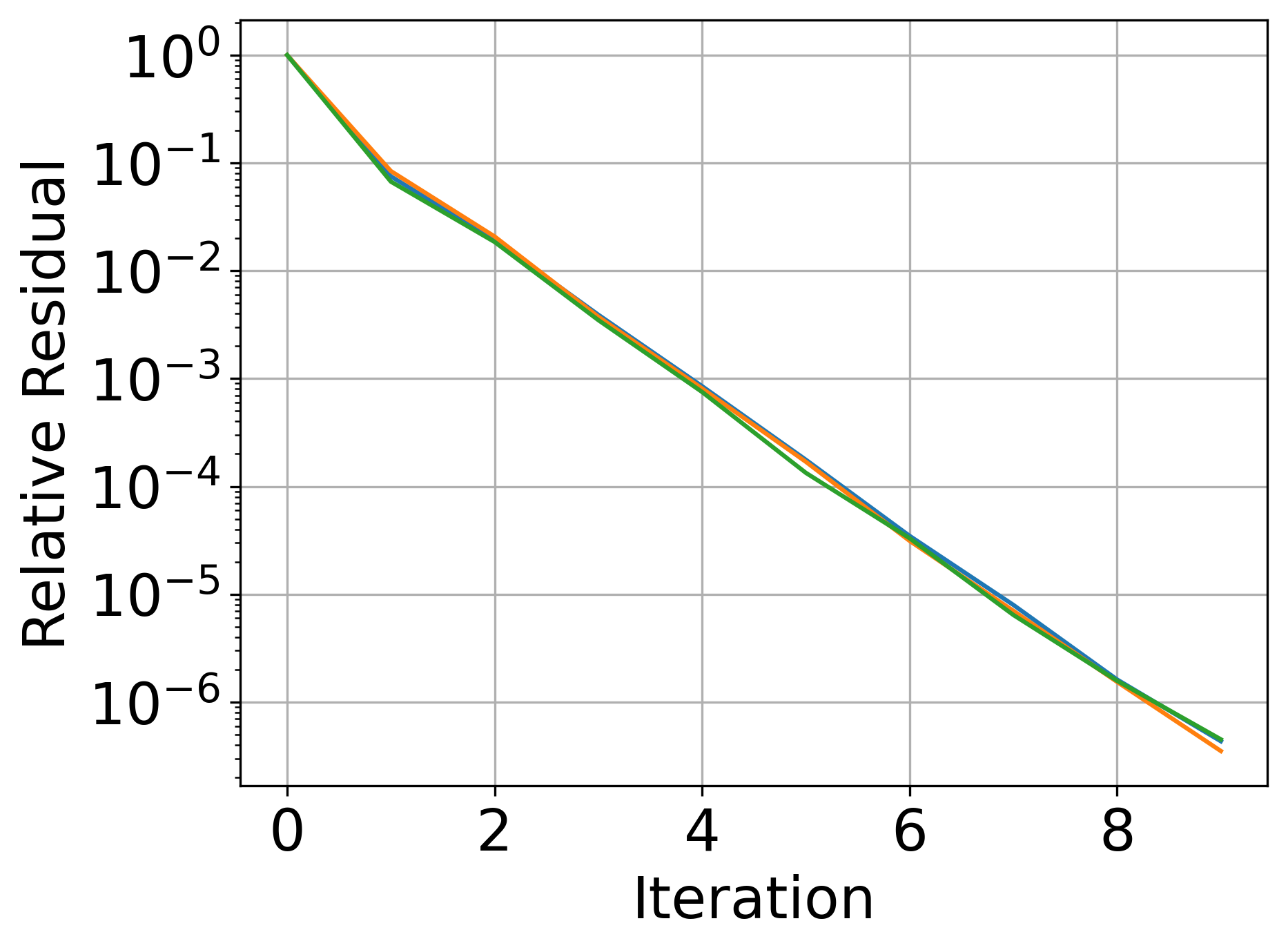}}
     \subfigure[Style-B]{\includegraphics[width=0.23\textwidth]{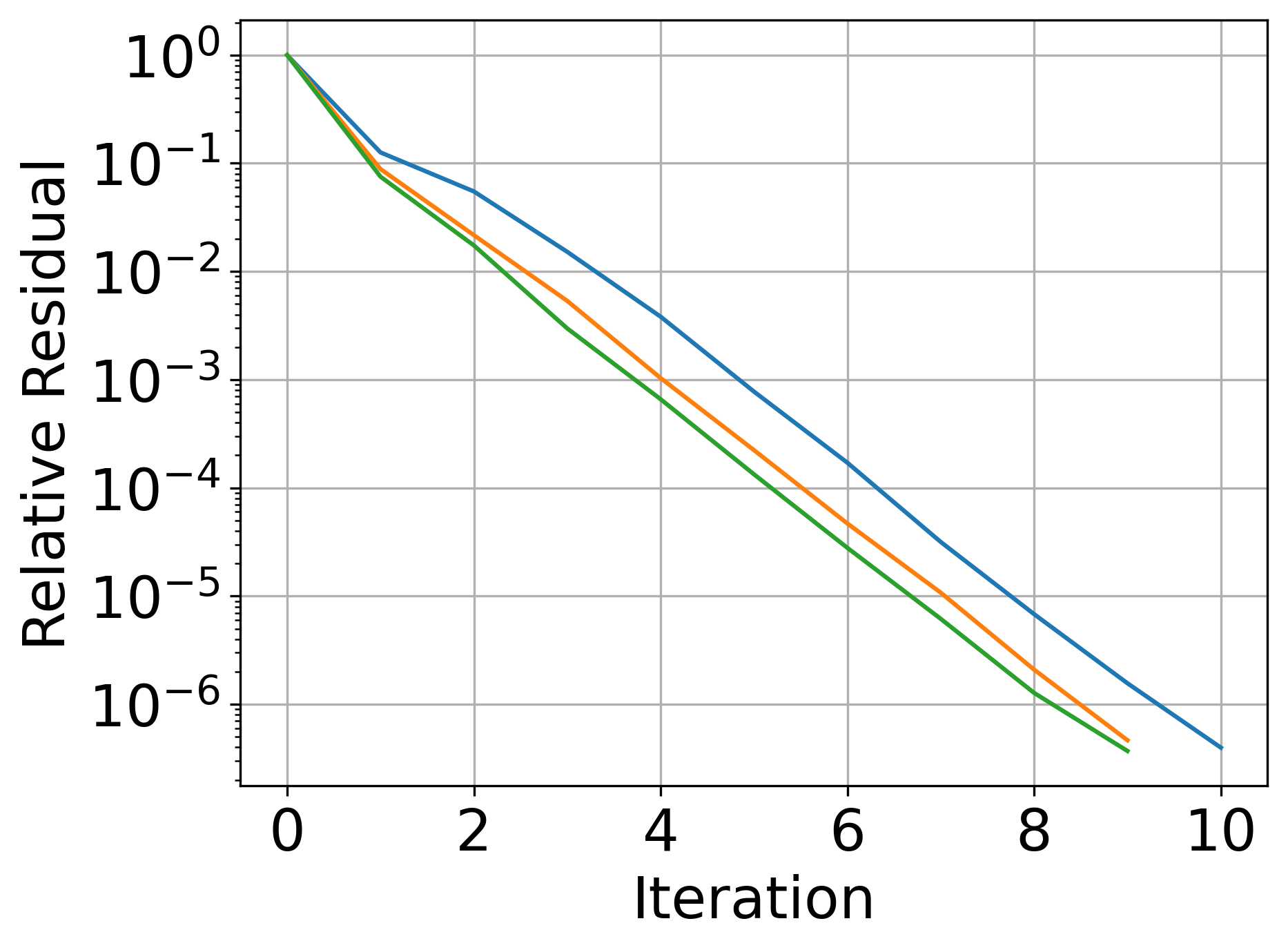}}
     \subfigure[STL-10 ]{\includegraphics[width=0.23\textwidth]{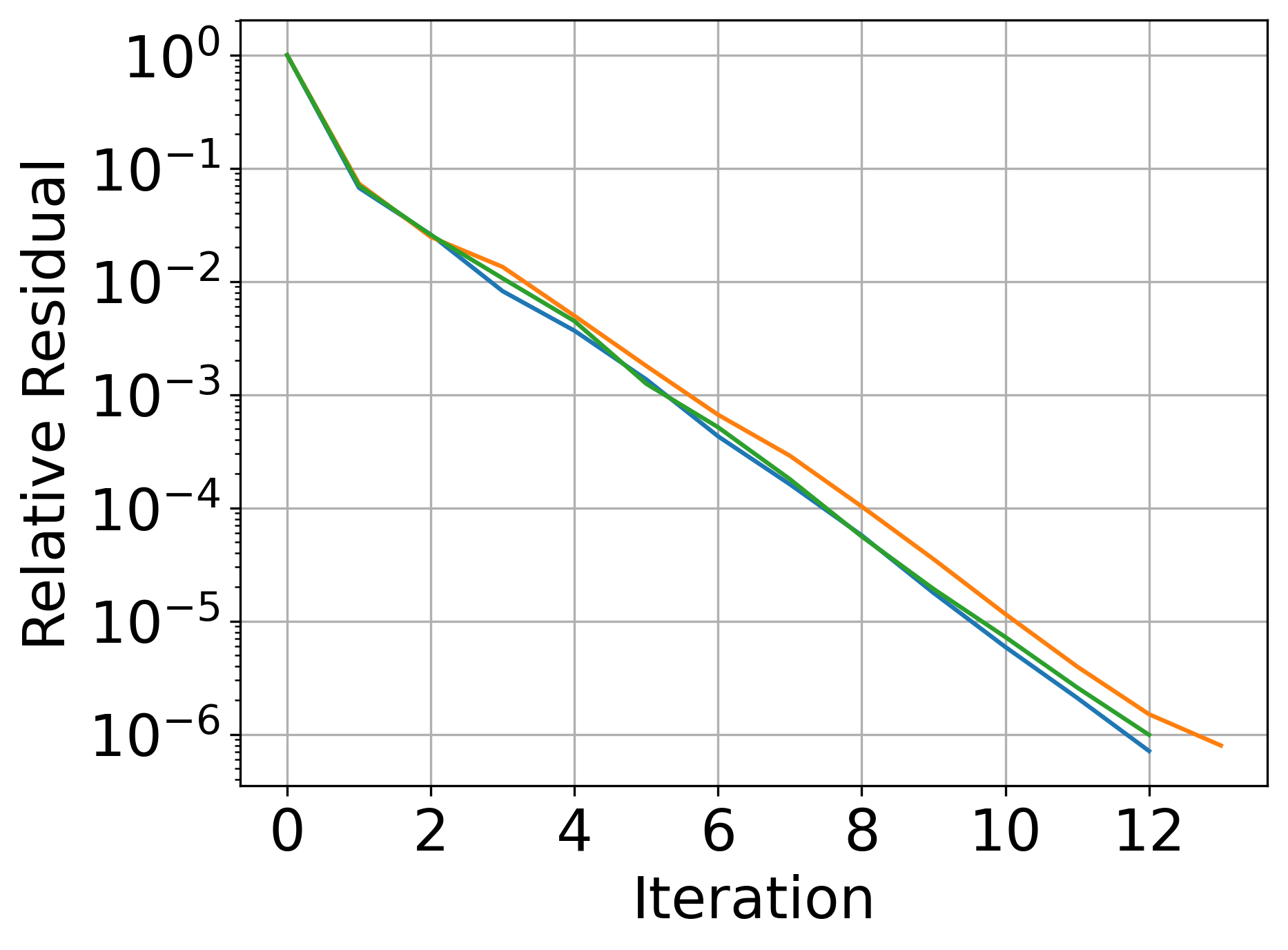}}
     \caption{Convergence history of Wave-ADR-NS for different right-hand sides.}
     \label{fig:diff_rhs}
\end{figure}

Then, we test the performance when solving Helmholtz equations with different right-hand sides. Recall that our training data includes only point sources located at the center of $\Omega$. We now assess the convergence speed for right-hand sides with different source locations, as well as random right-hand sides with normal distribution ($\mathbf{g}\sim \mathcal{N}(\mathbf{0},\mathbf{I})$). 
\cref{fig:diff_rhs} shows the convergence history of Wave-ADR-NS when solving Helmholtz equations with different right-hand sides when $\omega = 20\pi$ and slowness models from the four dataset types. It can be seen that the convergence speed of Wave-ADR-NS is independent of the right-hand sides. The results are similar for other source locations and values of $\omega$.

Next, we evaluate the performance of Wave-ADR-NS in solving the Helmholtz equation with a slight shift
$$
-\Delta u - k^2 u + i\gamma_0 u + i\gamma\frac{\omega}{c^2}u = g, \quad \mathbf{x} \in \Omega,
$$
where $\gamma_0 = 0.01 k^2$. This shift moves the imaginary part of the eigenvalues away from the real axis, mitigating the characteristic errors that are challenging to address and thereby alleviating the computational difficulties associated with solving it. \cite{cocquet2017large}.
We apply Wave-ADR-NS to solve the shifted Helmholtz equation across various slowness models and compare its performance with that of the Implicit Encoder Solver (IES) \cite{lerer2023multigrid}. For consistency, \cref{tab:shift} reports the number of iterations required to achieve a relative residual of $10^{-7}$. The results of IES are directly obtained from \cite{lerer2023multigrid}.
It can be observed that Wave-ADR-NS outperforms IES in terms of iteration counts when solving the shifted Helmholtz equation.
\begin{table}[!htb]
     \centering
     \footnotesize{
     \caption{Average iteration counts of Wave-ADR-NS and IES when solving shifted Helmholtz equations with different slowness models.}
     \label{tab:shift}
     \begin{tabular}{@{}cccccccc@{}}
     \toprule
     Dataset                                                                       & Method      & 128           & 256           & 512           & 1024          & 2048          & 4096                                 \\ \midrule
                                                                                   & Wave-ADR-NS    & \textbf{9.9}  & \textbf{15.7} & \textbf{24.9} & \textbf{38.6} & \textbf{58.4} & \textbf{90.2} \\
     \multirow{-2}{*}{CIFAR-10}                                                    & IES & 16.75         & 28.03         & 43.29         & 68.08         & 85.20         & 117.29                               \\
                                                                                   & Wave-ADR-NS    & \textbf{9.8}  & \textbf{14.6} & \textbf{24.6} & \textbf{38.6} & \textbf{56.7} & \textbf{89.1}                        \\
     \multirow{-2}{*}{\begin{tabular}[c]{@{}c@{}}OpenFWI\\ (Style A)\end{tabular}} & IES & 18.82         & 27.34         & 40.62         & 63.39         & 94.13         & 143.03                               \\
                                                                                   & Wave-ADR-NS    & \textbf{11.7} & \textbf{20.5} & \textbf{30.9} & \textbf{44.9} & \textbf{62.5} & \textbf{85.8}                        \\
     \multirow{-2}{*}{STL-10}                                                      & IES & 26.13         & 33.77         & 47.54         & 63.50        & 130.43        & 189.67                             \\ \bottomrule
     \end{tabular}}
\end{table}

\begin{figure}[!htb]
     \centering
     \includegraphics[width=\textwidth]{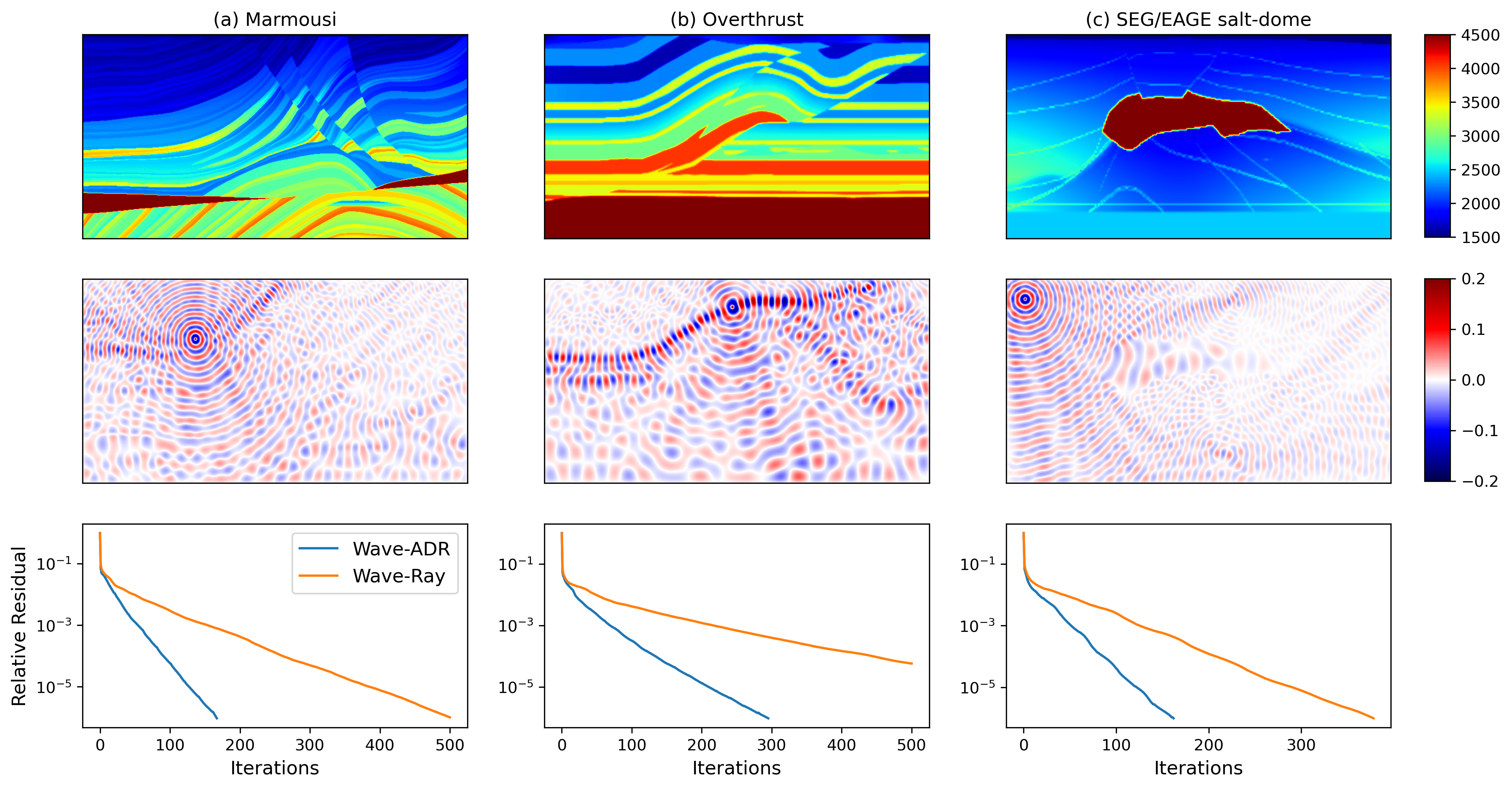}
     \caption{Non-square test. Top: velocity models used for each test. Middle: the corresponding solution of the point source Helmholtz equation. Bottom: convergence history of the Wave-ADR-NS and Wave-Ray preconditioners for each problem.}
     \label{fig:ood_results}
\end{figure}

Finally, we test the performance of Wave-ADR-NS on non-square computational domains, using the Marmousi \cite{brougois1990marmousi}, SEG/EAGE Salt-dome, and Overthrust \cite{aminzadeh1997models} models as test cases. 
\cref{fig:ood_results} presents the comparison with the Wave-Ray method. The results clearly demonstrate that for problems in non-square computational domains, Wave-ADR-NS continues to exhibit clear advantages.

\section{Conclusions and outlook}\label{sec:05}
In this paper, we develop the Wave-ADR-NS for solving high-frequency and heterogeneous Helmholtz equations. 
Wave-ADR-NS comprises two V-cycles. Firstly, the wave cycle, incorporating smoothers selected through spectral analysis, targets the elimination of non-characteristic components. Subsequently, the ADR cycle on the coarse level aims to diminish characteristic components. Parameters in both the smoother and the correction step are learned together by a CNN and a FNO, respectively. Numerical experiments demonstrate that Wave-ADR-NS achieves superior computational efficiency compared to existing multigrid methods and exhibits robust generalization capabilities for both in-distribution and out-of-distribution parameters.
While our development primarily focuses on the 2D scenario, formulating an equivalent model for 3D problems is straightforward.

There is still a lot of work to do in the future.
Firstly, we will enhance Wave-ADR-NS to be able to solve the Helmholtz equation in a wavenumber-independent fashion.
Secondly, we aim to integrate Wave-ADR-NS with an improved discretization method to avoid pollution effects and achieve more accurate numerical solutions.

\section*{Acknowledgments}
We would like to thank the referees for their insightful comments that greatly improved the paper. 
This work was carried out using computing resources at the High Performance Computing Platform of Xiangtan University.

\bibliographystyle{siamplain}
\bibliography{references}
\end{document}